\documentclass[twoside,reqno]{amsart}
\usepackage{amsmath}
\usepackage{amscd}
\usepackage{latexsym}
\usepackage{amsthm}
\usepackage{amssymb}

\theoremstyle{plain}
\newtheorem{thm}{Theorem}[section]

\newtheorem{lem}[thm]{Lemma}
\newtheorem{prop}[thm]{Proposition}

\theoremstyle{definition}
\newtheorem{defn}[thm]{Definition}

\newtheorem{notat}[thm]{Notation}
\newtheorem{hyp}[thm]{Hypothesis}

\theoremstyle{remark}
\newtheorem{rmk}[thm]{Remark}

\newcommand{\mc}{\mathcal}

\newcommand{\QQ}{\mathbb{Q}}
\newcommand{\ZZ}{\mathbb{Z}}
\newcommand{\CC}{\mathbb{C}}

\newcommand{\AAA}{\mathbb{A}}
\newcommand{\PP}{\mathbb{P}}
\newcommand{\GG}{\mathbb{G}}

\newcommand{\SP}{\text{Spec }}

\newcommand{\Kgnb}[2]{\overline{\mc M}_{#1}({#2})}
\newcommand{\Kbm}{\overline{\mc M}}

\newcommand{\EExt}[1]{\mathbb{E}\text{xt}_{\OO_{#1}}}

\newcommand{\lt}{\left}
\newcommand{\rt}{\right}

\newcommand{\OO}{\mathcal O}

\begin{document}

\title[Rational curves on hypersurfaces, II]{Rational curves on
hypersurfaces of low degree, II}  

\author[J. Harris]{Joe Harris} 
 \address{Department of Mathematics \\ 
   Harvard University \\ 
   Cambridge MA 02138}
 \email{harris@math.harvard.edu}

\author[J. Starr]{Jason Starr} 
 \address{Department of Mathematics \\ 
  Massachusetts Institute of Technology \\ 
  Cambridge MA 02139}
 \email{jstarr@math.mit.edu} \date{\today}

\begin{abstract}
  This is a continuation of ~\cite{HRS2} in which we proved
  irreducibility of spaces of rational curves on a general
  hypersurface $X_d\subset \PP^n$ of degree $d<\frac{n+1}{2}$.  In
  this paper, we prove that if $d^2 + d + 2 \leq n$ and if $d\geq 3$,
  then the spaces of rational curves are themselves rationally
  connected.
\end{abstract}

\maketitle

\tableofcontents

\section{Statement of results}~\label{sec-results}

In ~\cite{HRS2}, it is proved that if $X_d \subset \PP^n$ is a general
hypersurface of degree $d < \frac{n+1}{2}$, then each space
$\text{RatCurves}^e(X)$ parametrizing smooth rational curves of degree
$e$ on $X$, is itself an integral, local complete intersection scheme
of the expected dimension $(n+1-d)e+(n-4)$.  More precisely, it is
proved that for every stable $A$-graph $\tau$ and every flag $f\in
\text{Flag}(\tau)$, the Behrend-Manin stack $\Kbm(X,\tau)$ is an
integral, local complete intersection stack of the expected dimension
$\text{dim}(X,\tau)$, and the evaluation morphism
$\text{ev}_f:\Kbm(X,\tau) \rightarrow X$ is flat of the expected fiber
dimension $\text{dim}(X,\tau) - \text{dim}(X)$.  Since
$\text{RatCurves}^e(X)$ is a Zariski open set in the stack
$\Kgnb{0,0}{X,e}$, the result on $\text{RatCurves}^e(X)$ follows.

\

After establishing irreducibility and the dimension of the spaces
$\text{RatCurves}^e(X)$, the next question is to determine the Kodaira
dimension of $\text{RatCurves}^e(X)$.  For a general Fano hypersurface
$X_d \subset \PP^n$ with $d\leq n$, determining the Kodaira dimensions
of $\text{RatCurves}^e(X)$ is subtle.  For instance for smooth cubic
threefolds $X_3 \subset \PP^4$, $X$ has a nontrivial intermediate
Jacobian $J(X)$, and the Abel-Jacobi maps $\text{RatCurves}^e(X)
\rightarrow J(X)$ is dominant for $e \geq 4$.  So the Kodaira
dimension is at least $0$; conjecturally the fibers of the Abel-Jacobi
map are rationally connected so that the Kodaira dimension is exactly
$0$.

\

On the other hand for $d =1,2$, it is a theorem of Kim and
Pandharipande ~\cite[theorem 3]{KP}, that each of the spaces
$\text{RatCurves}^e(X)$ is itself a rational variety, and thus has
Kodaira dimension $-\infty$.  In this paper we present the following
generalization of ~\cite[theorem 3]{KP}:

\begin{thm} \label{thm-thm1}
  Given positive integers $(d,n)$ with $d^2 + d + 2 \leq n$ and $d\geq
  3$, for $X_d \subset \PP^n$ a general hypersurface of degree $d$,
  each of the spaces $\Kgnb{0,0}{X,e}$ is rationally connected.  Thus
  $\text{RatCurves}^e(X)$ has a rationally connected compactification.
\end{thm}

To remind the reader, a variety $V$ is rationally connected if given
two general points $p,q\in V$, there is a map $f:\PP^1 \rightarrow V$
whose image contains $p$ and $q$.  This property is strictly weaker
than rationality, and it is unknown whether this property is the same
as unirationality.  It is a priori a much simpler property to check.
And any rationally connected variety has Kodaira dimension $-\infty$,
therefore each of the schemes $\text{RatCurves}^e(X)$ has Kodaira
dimension $-\infty$.
   
\

The proofs rely on ~\cite[theorem IV.3.7]{K}: given a smooth variety
$V$ and a morphism $f:\PP^1 \rightarrow V$ such that $f^* T_V$ is an
ample vector bundle, then $V$ is rationally connected.  This criterion
also works for smooth Deligne-Mumford stacks $V$.  For readers not
versed in stacks, this is a moot point -- every morphism of $\PP^1$
into a stack constructed in this paper can be deformed to a map
contained in the fine moduli locus of the stack.  For a Behrend-Manin
stack, $\Kbm(X,\tau)$, a morphism $f:\PP^1 \rightarrow \Kbm(X,\tau)$
is equivalent to a fibered surface $\pi:\Sigma \rightarrow \PP^1$,
with some collection of sections $\sigma_1,\dots,\sigma_r:\PP^1
\rightarrow \Sigma$, and a map $f:\Sigma \rightarrow X$.  Assuming
that $f:\PP^1 \rightarrow \Kbm(X,\tau)$ maps into the
\emph{unobstructed locus}, the bundle $f^*T_{\Kbm(X,\tau)}$ can be
computed from a universal construction applied to the datum
$(\pi:\Sigma \rightarrow \PP^1,\sigma_1,\dots,\sigma_r,f)$.  Thus, to
prove $\Kbm(X,\tau)$ is rationally connected, we are reduced to
finding a datum $(\pi:\Sigma \rightarrow
\PP^1,\sigma_1,\dots,\sigma_r,\PP^1)$ satisfying certain properties.

\

In the proof of the induction step, we will use the following
technical hypotheses:

\begin{hyp}\label{hyp-1}
  For each contraction of genus $0$ stable $A$-graphs,
  $\phi:\sigma\rightarrow \tau$, the codimension of the image of the
  corresponding morphism of Behrend-Manin stacks $\Kbm(X,\sigma)
  \rightarrow \Kbm(X,\tau)$ equals $\text{dim}(X,\tau) -
  \text{dim}(X,\sigma)$.
\end{hyp}

In particular, by \cite[proposition 7.4]{HRS2}, for a general
$X_d\subset \PP^n$ with $d<\frac{n+1}{2}$, each stack $\Kbm(X,\sigma)$
has the expected dimension, and hypothesis~\ref{hyp-1} holds for $X$.

\begin{hyp}\label{hyp-1.5}
  The general fiber of the evaluation map $\text{ev}:\Kgnb{0,1}{X,1}
  \rightarrow X$ is irreducible.
\end{hyp}

In particular, for a general pair $([X],[p])$ consisting of a
hypersurface $X\subset \PP^n$ of degree $d$ and a point $p\in X$, the
associated fiber of $\text{ev}$ is a subvariety $Z\subset \PP^{n-1}$
which is a complete intersection of a general sequence of
hypersurfaces $Y_1,\dots, Y_d$ in $\PP^n$ with $\text{deg}(Y_i)=i$.
By the Bertini-Kleiman theorem, a general such complete intersection
is smooth and connected if $d < n-1$.

\begin{hyp}\label{hyp-1.75}
  For each integer $e\geq 0$, the locus in $\Kgnb{0,1}{X,e}$
  parametrizing stable maps with nontrivial automorphism group has
  codimension at least $2$.
\end{hyp}

Of course any stable map with nontrivial automorphism group has an
irreducible component which is a multiple cover of its image.  In
light of \cite[proposition 7.4]{HRS2}, a simple parameter count shows
that for a general hypersurface $X_d \subset \PP^n$ with $d\leq
\frac{n+1}{2}$, hypothesis~\ref{hyp-1.75} is satisfied.

\subsection{Acknowledgments}

We are very grateful to Johan de Jong and Steven Kleiman for many
useful conversations.

\section{Deformation ample}~\label{sec-DA}
Unless stated otherwise, all schemes will be finite type, separated
schemes over $\SP \CC$.  All absolute fiber products will be fiber
products over $\SP \CC$.

\

The following lemma is completely trivial.

\begin{lem}~\label{lem-gend}
  Let $B$ be a connected, proper, prestable curve of arithmetic genus
  $0$ and let $E$ be a locally free sheaf on $E$ which is generated by
  global sections.  Then $H^1(B,E)$ is zero.  Moreover if $p\in B$ is
  any point and ${\mc I}_p\subset \OO_B$ is the corresponding ideal
  sheaf, then $H^1(B,{\mc I}_p\cdot E)$ is zero.
\end{lem}

\begin{proof}
  Since $E$ is generated by global sections, there is a short exact
  sequence of the form:
\begin{equation}
\begin{CD}
0 @>>> K @>>> \OO_B^{\oplus N} @>>> E @>>> 0
\end{CD}
\end{equation}
Since $B$ is a curve, $H^2(B,K)=0$.  Therefore we have a surjection
$H^1(B,\OO_B)^{\oplus N} \rightarrow H^1(B,E)$.  Since $B$ is
connected of arithmetic genus $0$, $H^1(B,\OO_B)$ is zero.  Hence
$H^1(B,E)$ is zero.

\

Now for any point $p\in B$ we have a short exact sequence:
\begin{equation}
\begin{CD}
0 @>>> {\mc I}_p\cdot E @>>> E @>>> E|_p @>>> 0
\end{CD}
\end{equation}
This gives rise to a long exact sequence in cohomology:
\begin{equation}
\begin{CD}
H^0(B,E) @>>> E|_p @>>> H^1(B,{\mc I}_p\cdot E) @>>> H^1(B,E)
\end{CD}
\end{equation}
By the last paragraph, $H^1(B,E)$ is zero.  And since $E$ is generated
by global sections, $H^0(B,E)\rightarrow E|_p$ is surjective.
Therefore $H^1(B,{\mc I}_p \cdot E)$ is zero.
\end{proof}

Suppose given a connected, proper, prestable curve $B$ of arithmetic
genus $0$ and a locally free sheaf $E$ of positive rank on $B$.  Any
pair $({\mc B},{\mc E})$ consisting of a flat family of connected,
proper prestable curves over a DVR, say $\pi:{\mc B} \rightarrow \SP
R$ along with a locally free sheaf ${\mc E}$ on ${\mc B}$ such that
the closed fiber of $\pi$ is isomorphic to $B$, such that the generic
fiber $B_\eta$ of $\pi$ is smooth, and such that the restriction of
${\mc E}$ to $B\subset {\mc B}$ is isomorphic to $E$ will be called a
\emph{smoothing} of the pair $(B,E)$.  We want to know when $(B,E)$
satisfies the condition that for every smoothing $({\mc B},{\mc E})$,
the restriction ${\mc E}_\eta$ of ${\mc E}$ to the generic fiber is an
ample locally free sheaf.  Certainly if $E$ is ample, this is true,
but $E$ need not be ample for this condition to hold: e.g. if $E$ is
an invertible sheaf such that the total degree of $E$ is positive,
then every smoothing of $(B,E)$ will have ample generic fiber.
Although it is not the most general criterion, we find the following
notion to be useful and it is the one we use in the remainder of the
paper.

\

\begin{defn}~\label{def-DA}
  Let $B$ be a connected, proper, prestable curve of arithmetic genus
  $0$.  A locally free sheaf $E$ on $B$ with positive rank is
  \emph{deformation ample} if
  \begin{enumerate}
    \item $E$ is generated by global sections, and
    \item $H^1(B,E(K_B))$ is zero, where $\OO_B(K_B)$ is the dualizing
      sheaf of $B$.
  \end{enumerate}
\end{defn}

\begin{rmk}~\label{rmk-DA}
  \begin{enumerate}
    \item Conditions (1) and (2) above are independent.
    \item If $E$ is invertible, then $E$ is deformation ample iff the
      restriction of $E$ to every irreducible component has
      nonnegative degree and the restriction to at least one
      irreducible component has positive degree.
    \item For a general $E$, one can determine whether $E$ is
      deformation ample in terms of the splitting type of the
      restriction of $E$ to each irreducible component together with
      the patching isomorphisms at the nodes of $B$.
  \end{enumerate}
\end{rmk}

\begin{defn}~\label{defn-fDA}
  Let $T$ be a scheme and let $\pi:B\rightarrow T$ be a family of
  prestable curves of arithmetic genus $0$.  A locally free sheaf $E$
  on $B$ with positive rank is \emph{$\pi$-relatively deformation
    ample} (or simply deformation ample if $\pi$ is understood) if
  \begin{enumerate}
    \item the canonical map $\pi^*\pi_*E\rightarrow E$ is surjective, and
    \item $R^1\pi_*(E(K_\pi))$ is zero, where $\OO_B(K_\pi)$ is the
      relative dualizing sheaf of $\pi$.
  \end{enumerate}
\end{defn}

\begin{lem}~\label{lem-bcDA}
  With notation as in definition~\ref{defn-fDA}, suppose
  $f:T'\rightarrow T$ is a morphism of schemes and let
  $\pi':B'\rightarrow T'$ be the base-change of $\pi$ by $f$.  Let
  $E'$ be the pullback of $E$ by the projection $g:B'\rightarrow B$.
  If $E$ is $\pi$-relatively deformation ample, then $E'$ is
  $\pi'$-relatively deformation ample.  If $f'$ is surjective, the
  converse also holds.
\end{lem}

\begin{proof}
  For the main direction, by ~\cite[section 8.5.2, proposition
  8.9.1]{EGA4}, it suffices to consider the case when $T$ and $T'$ are
  Noetherian affine schemes.

\

There is a canonical map of $\OO_{T'}$-modules $\alpha:f^*\pi_*
E\rightarrow (\pi')_*g^* E$.  There is a commutative diagram:
\begin{equation}
\begin{CD}
(\pi')^*f^*\pi_* E @> = >> g^* \pi^* \pi_* E \\
@V (\pi')^*\alpha VV @VVV \\
(\pi')^*(\pi')_* E' @>>> E'
\end{CD}
\end{equation}
Since $\pi^*\pi_*E \rightarrow E$ is surjective, we conclude that
$g^*\pi^*\pi_* E\rightarrow E'$ is surjective.  Therefore also
$(\pi')^*(\pi')_*E' \rightarrow E'$ is surjective.

\

Now $R^2\pi_* E(K)$ is identically zero.  So by \cite[theorem
III.12.11(b)]{H}, we conclude that for every closed point $t\in T$, we
have $H^1(B_t,E(K)|_{B_t}) = 0$.  By ~\cite[prop. III.9.3]{H}, we
conclude that for every closed point $t'\in T'$, we have
$H^1(B'_{t'},E'(K')|_{B'_{t'}}) = 0$.  So by \cite[theorem
III.12.11(a)]{H} and Nakayama's lemma, we conclude that
$R^1\pi'_*(E'(K'))$ is identically zero.  This shows that $E'$ is
$\pi'$-relatively deformation ample.

\

For the converse in case $T'\rightarrow T$ is surjective, we may
reduce to the case that $T$ is a Noetherian affine scheme.  By the
argument above, we conclude that for each closed point $t'\in T'$, we
have $H^1(B'_{t'},E'(K')|_{B'_{t'}}) = 0$.  Since $T'\rightarrow T$ is
surjective, we conclude by ~\cite[prop. III.9.3]{H} that for each
closed point $t\in T$ we have $H^1(B_t,E(K)|_{B_t}) = 0$.  So by
~\cite[theorem III.12.11(a)]{H} and Nakayama's lemma, we conclude that
$R^1\pi_*(E(K))$ is identically zero.

\

It remains to show that $E$ is $\pi$-relatively generated by global
sections.  For each point $t\in T$, there is some point $t'\in T'$
mapping to $t$.  Since $E'|_{B'_{t'}}$ is generated by global
sections, it follows that $E|_{B_t}$ is generated by global sections.
So, by lemma~\ref{lem-gend}, we conclude that $H^1\lt( B_t, E|_{B_t}
\rt)$ is zero.  So by ~\cite[theorem III.12.11(a)]{H} and Nakayama's
lemma, we conclude that $R^1\pi_*(E)$ is identically zero.  Since $E$
is flat over $T$, it follows by a standard argument that for any
coherent $\OO_T$-module ${\mc F}$, we have that $R^1\pi_*( \pi^*{\mc
  F}\otimes E)$ is also zero.  In particular, applying the long exact
sequence of higher direct images to the short exact sequence:
\begin{equation}
\begin{CD}
0 @>>> \pi^*{\mc I}_t\otimes E @>>> E @>>> E|_{B_t} @>>> 0
\end{CD}
\end{equation}
we conclude that $\pi_*(E) \rightarrow H^0\lt( B_t, E|_{B_t} \rt)$ is
surjective.  Since $E|_{B_t}$ is generated by global sections for each
$t\in T$, we conclude that $E$ is $\pi$-relatively generated by global
sections.
\end{proof}

\begin{lem}\label{lem-smDA}
  With notation as in definition~\ref{defn-fDA}, if $\pi$ is smooth,
  then $E$ is $\pi$-relatively deformation ample iff $E$ is
  $\pi$-relatively ample.
\end{lem}

\begin{proof}
  Both properties are local on $T$ and can be checked after \'{e}tale,
  surjective base-change of $T$.  Thus we may assume that
  $\pi:B\rightarrow T$ is isomorphic to $\pi_1:T\times \PP^1
  \rightarrow T$.  Define $F=(\pi_1)_*(E\otimes
  \pi_2^*\OO_{\PP^1}(-1))$.  There is a natural map $\alpha:
  (\pi_1)^*F \otimes\pi_2^*\OO_{\PP^1}(1)\rightarrow E$.  Suppose that
  $E$ is deformation ample.  Then the claim is that $\alpha$ is
  surjective.  To prove this, it suffices to prove the following:
\begin{enumerate}
\item For each geometric point $t$ of $T$ with residue field $\kappa$,
  we have $H^1(\PP^1_\kappa, E|_{B_t}\otimes \OO_{\PP^1}(-1))$ is
  trivial,
\item $\pi^* F|_{B_t} = H^0(\PP^1_\kappa,E|_{B_t}\otimes
  \OO_{\PP^1}(-1))$, and
\item the map $H^0(\PP^1_\kappa,E|_{B_t}\otimes
  \OO_{\PP^1}(-1))\otimes \OO_{\PP^1}(1) \rightarrow E|_{B_t}$ is
  surjective. 
\end{enumerate}

\

Now by Grothendieck's lemma ~\cite[exercise V.2.6]{H}, $E|_{B_t}$
splits as a direct sum $\OO_{\PP^1}(a_1)\oplus \dots\oplus
\OO_{\PP^1}(a_r)$ for some integers $a_1\leq \dots \leq a_r$.  By
lemma~\ref{lem-bcDA}, we know that $E|_{B_t}$ is deformation ample.
It follows that $a_1\geq 1$.  Thus $H^1(B_t, E|_{B_t}\otimes
\OO_{\PP^1}(-1))=0$ and $(1)$ is established.  Combined
with~\cite[theorem III.12.11(b)]{H}, also $(2)$ follows.  Finally, for
$a_i \geq 1$, we clearly have $H^0(\PP^1_\kappa,
\OO_{\PP^1}(a_i-1))\otimes \OO_{\PP^1}(1)\rightarrow \OO_{\PP^1}(a_i)$
is surjective.  Thus $(3)$ holds and the claim is proved.  Now
$(\pi_1)^*F\otimes (\pi_2)^*\OO_{\PP^1}(1)$ is $\pi_1$-relatively
ample.  Since $E$ is a quotient of $(\pi_1)^*F\otimes
(\pi_2)^*\OO_{\PP^1}(1)$, we conclude that $E$ is also
$\pi_1$-relatively ample.

\

The converse result follows in the same way.
\end{proof}

\begin{lem}~\label{lem-opDA}
  With notation as in lemma~\ref{defn-fDA}, there exists an open
  subscheme $i:U\rightarrow T$ with the following property: for every
  morphism $f:T'\rightarrow T$, $f(T')$ is contained in $U$ iff $E'$
  is $\pi'$-relatively deformation ample.
\end{lem}

\begin{proof}
  By \cite[section 8.5.2,proposition 8.9.1]{EGA4}, we may reduce to
  the case that $T$ and $T'$ are Noetherian affine schemes.

\

Let $Z_1\subset T$ be the closed subset with is the image under $f$ of
the support of $\text{coker}(\pi^*\pi_* E\rightarrow E)$.  Let
$Z_2\subset T$ be the closed subset which is the support of
$R^1\pi_*(E(K_\pi))$.  Let $i:U\rightarrow T$ be the open complement
of $Z_1\cup Z_2$.

\

Suppose $f:T'\rightarrow T$ is a morphism of schemes.  By
~\cite[theorem III.12.11, prop. III.9.3]{H} and Nakayama's lemma,
$R^1\pi'_*(E'(K'))$ is identically zero iff for each $t'\in T'$ with
$t=f(t')$, we have $H^1(B_t,E(K)|_{B_t})$ is zero, i.e. if $t$ is
contained in the complement of $Z_2$.  So $R^1\pi'_*(E'(K'))$ is
identically zero iff $f(T)$ is contained in the complement of $Z_2$.
If $R^1\pi'_*(E'(K'))$ is zero, then also $R^1\pi'_*(E')$ is zero (by
the same argument as in the proof of lemma~\ref{lem-bcDA}).  So using
~\cite[theorem III.12.11, prop. III.9.3]{H} again, $E'$ is
$\pi'$-relatively generated by global sections iff for each $t'\in T'$
with $t=f(t')$, we have $E_t$ is generated by global sections, i.e.
$t$ is in the complement of $Z_1$.  So we conclude that $E'$ is
$\pi'$-relatively deformation ample iff $f(T)$ is contained in the
complement of $Z_1\cup Z_2$.
\end{proof}

\begin{lem}~\label{lem-secDA} 
  With notation as in definition~\ref{defn-fDA}:
\begin{enumerate}
\item If $E_1\rightarrow E_2$ is a morphism of locally free
  sheaves on $B$ whose cokernel is torsion in every fiber
  (in particular, if the
  morphism is surjective), if $E_2$ is nonzero, 
  and if $E_1$ is $\pi$-relatively deformation ample,
  then also $E_2$ is $\pi$-relatively deformation ample.
\item Given a short exact sequence of nonzero locally free sheaves on
  $B$,
\begin{equation}
\begin{CD}
0 @>>> E_1 @>>> E_2 @>>> E_3 @>>> 0 
\end{CD}
\end{equation}
if $E_1$ and $E_3$ are $\pi$-relatively deformation ample, then $E_2$
is also $\pi$-relatively deformation ample.
\item If $E$ is $\pi$-relatively deformation ample, then for each integer
  $n\geq 1$, also $E^{\otimes n}$ is $\pi$-relatively deformation
  ample.
\end{enumerate}
\end{lem}

\begin{proof}
  First we prove $(1)$.  Let $Q$ denote the cokernel, which is assumed
  to be torsion.  Now $R^1\pi_*$ is right exact on coherent sheaves,
  $R^1\pi_*(Q(K))$ is zero (because $Q(K)$ is torsion in fibers), and
  $R^1\pi_*(E_1(K))=0$ by assumption.  So we conclude that also
  $R^1\pi_*(E_2(K))=0$.  Let $I\subset E_2$ denote the image sheaf of
  $E_1\rightarrow E_2$.  The surjective composition map
\begin{equation}
\begin{CD}
\pi^*\pi_*E_1 @>>> E_1 @>>> I
\end{CD}
\end{equation}
factors through $\pi^*\pi_* I\rightarrow I$.  Therefore $I$ is
$\pi$-relatively generated by global sections.  In particular,
$R^1\pi_* I$ is zero (because $R^1\pi_* \OO_B$ is zero).  Since $Q$ is
fiberwise torsion, it is $\pi$-relatively generated by global
sections.  Now we have a short exact sequence of coherent sheaves:
\begin{equation}
\begin{CD}
0 @>>> I @>>> E_2 @>>> Q @>>> 0.
\end{CD}
\end{equation}
Applying the long exact sequence of higher direct images, and using
that $R^1\pi_*I$ is zero, locally on $B$ all global sections of $Q$
lift to global sections of $E_2$.  Since also $I$ is $\pi$-relatively
generated by global sections, we conclude that $E_2$ is
$\pi$-relatively generated by global sections.  Thus $E_2$ is
$\pi$-relatively deformation ample.

\

Next we prove $(2)$.  The long exact sequence of higher direct images
and the vanishings $R^1\pi_*(E_1(K))=R^1\pi_*(E_3(K))=0$ imply that
$R^1\pi_*(E_2(K))$ is also trivial.  Now $E_1$ and $E_3$ are
$\pi$-relatively generated by global sections an so (locally on $T$)
they are quotients of a trivial sheaf $\OO_B^{\oplus N}$.  Since
$R^1\pi_*\OO_B$ is zero, and since $R^1\pi_*$ is right exact for short
exact sequences of coherent sheaves, we conclude that $R^1\pi_*E_1 =
R^1\pi_* E_3 = 0$.  Applying the long exact sequence of higher direct
images to the short exact sequence above, we conclude that
$\pi_*(E_2)\rightarrow \pi_*(E_3)$ is surjective.  Since $E_3$ is
$\pi$-relatively generated by global sections, we conclude that
$\pi^*\pi_*E_2\rightarrow E_2 \rightarrow E_3$ is surjective.
Therefore the map from $E_1$ to the cokernel of
$\pi^*\pi_*E_2\rightarrow E_2$ is surjective.  But since $E_1$ is
$\pi$-relatively generated by global sections, this map is zero and we
conclude $\pi^*\pi_*E_2\rightarrow E_2$ is surjective.  So $E_2$ is
$\pi$-relatively deformation ample.

\

Finally we prove $(3)$ by induction on $n$.  It suffices to consider
the case when $T$ is affine.  For $n=1$ it is trivial.  Suppose $n>1$
and suppose the result is known for all smaller values of $n$.  In
particular, $E^{\otimes(n-1)}$ is $\pi$-relatively generated by global
sections, and we have a natural surjection
\begin{equation}
\pi_*\pi^*(E^{\otimes(n-1)})\otimes_{\OO_B} E \rightarrow E^{\otimes n}.
\end{equation}
Now we can find a surjective map $\OO_T^{\oplus r}\rightarrow
\pi_*(E^{\otimes(n-1)})$.  Thus we have a surjection $E^{\oplus
  r}\rightarrow E^{\otimes n}$.  By $(2)$ above and induction, we have
that $E^{\oplus r}$ is $\pi$-relatively deformation ample.  By $(1)$
above, we conclude that $E^{\otimes n}$ is $\pi$-relatively
deformation ample, and we have proved $(3)$ by induction on $n$.
\end{proof}

\begin{lem}~\label{lem-DAcrit}
  Suppose now that $T=\SP k$ for some algebraically closed field $k$.
  If $E$ satisfies the hypotheses
\begin{enumerate}
\item For every irreducible component $B_i\subset B$ we have
  $E|_{B_i}$ is generated by global sections, and
\item for some nonempty closed subcurve $B'\subset B$, $E|_{B'}$ is
  deformation ample,
\end{enumerate}
then $E$ is deformation ample.
\end{lem}

\begin{proof}
  Let $\delta$ be the number of irreducible components of $B$ which
  aren't contained in $B'$.  We prove the result by induction on
  $\delta$.  If $\delta=0$, then $B=B'$ and there is nothing to prove.
  Thus suppose that $\delta >0$ and suppose the result has been proved
  for all smaller values of $\delta$.

\

Let $B_1\subset B$ be an irreducible component of $B$ which is not in
$B'$, let $B_2\subset B$ denote the union of all irreducible
components other than $B_1$, and suppose that $B_1$ intersects $B_2$
in precisely one node $p$ of $B$.  By the induction assumption
$E|_{B_2}$ is deformation ample.

\

First we prove that $E$ is generated by global sections.  Let
$F\subset E$ be the image of $H^0(B,E)\otimes_k \OO_B \rightarrow E$.
We have a short exact sequence of coherent sheaves:
\begin{equation}
\begin{CD}
0 @>>> E|_{B_1}(-p) @>>> E @>>> E|_{B_2} @>>> 0.
\end{CD}
\end{equation}
Since $E|_{B_1}$ is a locally free sheaf on $\PP^1$ generated by
global sections, Grothendieck's lemma and cohomology of line bundles
on $\PP^1$ imply that $H^1(B,E|_{B_1}(-p))=0$.  Thus all the global
sections of $E|_{B_2}$ lift to global sections of $E$, i.e.
$F\rightarrow E|_{B_2}$ is surjective.  So $E/F$ is supported on $B_1$
and thus is a quotient of $E|_{B_1}$.  Since $E|_{B_1}$ is generated
by global sections, also $E/F$ is generated by global sections.  We
have a short exact sequence
\begin{equation}
\begin{CD}
0 @>>> F @>>> E @>>> E/F @>>> 0.
\end{CD}
\end{equation}
Since $H^1(B,\OO_B)=0$, also $H^1(B,H^0(B,E)\otimes_k \OO_B) = 0$.
Since $H^1(B,*)$ is right exact, $H^1(B,F)=0$.  Thus all the global
sections of $E/F$ lift to global sections of $E$.  This can only hold
if $E/F=0$, i.e. if $E$ is generated by global sections.

\

Next we prove that $H^1(B,E(K_B))=0$.  We have a short exact sequence
of sheaves:
\begin{equation}
\begin{CD}
0 @>>> E(K_B)|_{B_2}(-p) @>>> E(K_B) @>>> E(K_B)|_{B_1} @>>> 0.
\end{CD}
\end{equation}
This gives an exact sequence of vector spaces:
\begin{equation}
\begin{CD}
H^1(B,E(K_B)|_{B_2}(-p)) @>>> H^1(B,E(K_B)) @>>> H^1(B,E(K_B)|_{B_1})
@>>> 0.
\end{CD}
\end{equation}
By standard results on dualizing sheaves and finite morphisms, we have
\begin{equation}
K_{B_2}=\textit{Hom}_{\OO_B}(\OO_{B_2},K_B) = K_B|_{B_2}(-p).
\end{equation}
Therefore $H^1(B,E(K_B)|_{B_2}(-p)) = H^1(B_2,E|_{B_2}(K_{B_2}))$,
which is zero by the induction assumption.  Similarly, $E(K_B)|_{B_1}$
equals $E|_{B_1}(-1)$ (identifying $B_1$ with $\PP^1$).  Since
$E|_{B_1}$ is generated by global sections, it follows by
Grothendieck's lemma and cohomology of line bundles on $\PP^1$ that
$H^1(B_1,E|_{B_1}(-1))$ is trivial.  Therefore we conclude that
$H^1(B,E(K_B))$ is trivial.  So $E$ is deformation ample, and the
result is proved by induction on $\delta$.
\end{proof}

\begin{rmk} \label{rmk-DAcrit}
  A particular case of lemma~\ref{lem-DAcrit} is when $B'$ is one
  irreducible component of $B$.  Then the lemma says that a locally
  free sheaf on $B$ which is \emph{generically ample} in the sense of
  Lazarsfeld ~\cite{F} is deformation ample.
\end{rmk}

\section{Some deformation theory} \label{sec-def}
In the next section we will need some deformation theory of stable
maps.

\begin{notat} \label{not-cpx}
  Suppose $T$ is a scheme and suppose
\begin{equation}
\zeta = ((\pi:B\rightarrow T,\sigma_1,\dots, \sigma_r), g:B
\rightarrow X)
\end{equation}
is a family of marked prestable maps to a smooth scheme $X$.  Denote
by $L_\zeta$ the complex of coherent sheaves on $B$
\begin{equation}
\begin{CD}
-1 & & 0 \\
f^*\Omega_X @> (df)^\dagger >> \Omega_\pi(\sigma_1 + \dots + \sigma_r)
\end{CD}
\end{equation}
Denote by $L_\zeta^\vee$ the object
\begin{equation}
L_\zeta^\vee :=
\mathbb{R}\textit{Hom}_{\OO_B}(L_\zeta,\OO_B)
\end{equation} 
in the derived category of $B$.
\end{notat}

The relevance of the complex $L_\zeta^\vee$ is the following.

\begin{lem} \label{lem-cpx}
  Suppose $X$ is a smooth projective scheme.  Let $\pi:{\mc
    B}\rightarrow \Kgnb{g,r}{X,\beta}$ denote the universal curve, let
  $\sigma_i: \Kgnb{g,r}{X,\beta} \rightarrow {\mc B}$ denote the
  universal sections, and let $f:{\mc B} \rightarrow X$ denote the
  universal map, i.e.
\begin{equation}
\zeta = \lt( \lt( \pi:{\mc B} \rightarrow \Kgnb{g,r}{X,\beta},
\sigma_1, \dots, \sigma_r \rt), f:{\mc B} \rightarrow X \rt)
\end{equation}
is the universal family of stable maps.  There is an \emph{obstruction
  theory} for $\Kgnb{g,r}{X,\beta}$ in the sense of \cite[definition
4.4]{BF} of the form
\begin{equation}
\phi:\lt( \mathbb{R}\pi_*(L_\zeta^\vee) [1] \rt)^\vee \rightarrow
L_{\Kgnb{g,r}{X,\beta}}.
\end{equation}
\end{lem}

\begin{proof} Essentially this follows from ~\cite{BF} and
  ~\cite{B}.
\end{proof}

\begin{rmk} \label{rmk-cpx}
  Explicitly, if $\zeta = ((B,p_1,\dots,p_r),f:B\rightarrow X)$ is a
  stable map (i.e $T = \SP C$), then the space of first order
  deformations of $\zeta$ is given by $\EExt{B}^1(L_\zeta,\OO_B)$ and
  the obstruction group is a subgroup of $\EExt{B}^2(L_\zeta,\OO_B)$.
  In particular, if $\EExt{B}^2(L_\zeta,\OO_B)$ vanishes, then
  $\Kgnb{g,r}{X,\beta}$ is smooth at the point $[\zeta]$.
\end{rmk}

\subsection{Contracting unstable components}\label{subsec-unstable}

In this subsection we wish to investigate the relationship between
$\EExt{B}(L_\zeta, \OO_B)$ and $\EExt{B'}(L_\zeta',\OO_B')$ where
\begin{equation}
\zeta = ((B,p_1,\dots,p_r),f:B \rightarrow X)
\end{equation}
is a prestable map, where
\begin{equation}
(B',p'_1,\dots, p'_r,q'_1,\dots,q'_s)
\end{equation}
is a prestable curve, where
\begin{equation}
h:(B',p'_1,\dots, p'_r) \rightarrow (B,p_1,\dots, p_r)
\end{equation} 
is a map which contracts some of the unstable components of
$(B',p_1,\dots,p_r)$, where $f' = f\circ h$ and where $\zeta'$ is the
prestable map
\begin{equation}
\zeta' = ((B',p'_1,\dots, p'_r,q'_1, \dots, q'_s), f':B' \rightarrow
X).
\end{equation}

\

Any morphism $h:B'\rightarrow B$ as above can be factored as a
sequence of \emph{elementary} maps.  We begin by analyzing the case
$\zeta = (B,f:B\rightarrow X)$ is a prestable map without marked
points, $h:B'\rightarrow B$ contracts a single unstable component, and
$\zeta' = (B', f' = f\circ h:B' \rightarrow X)$ (again without marked
points).  We call such maps either \emph{type I} or \emph{type II},
depending on whether the image under $h$ of the contracted component
is a smooth point of $B$ or a node of $B$.  After this, we analyze the
case where $\zeta = (B,(p_1,\dots,p_n),f:B\rightarrow X)$ is a marked
prestable map, where $B'$ is the same prestable map, but with one
extra marked point, and where $h:B'\rightarrow B$ is the identity map.
We call such a map \emph{type III}.

\

To simplify calculations, we replace $L_\zeta$ by a quasi-isomorphic
complex as follows.  Choose a regular embedding $i:B\rightarrow S$ of
$B$ into a smooth surface $S$.  Then the morphism $(f,i):B \rightarrow
X \times S$ is a regular embedding.  Let $N_{(f,i)}$ denote the normal
bundle of the embedding.  There are induced morphisms $\alpha_f:
N_{(f,i)}^\vee \rightarrow f^* \Omega_X$ and $\beta_i: N_{(f,i)}^\vee
\rightarrow i^* \Omega_S$ which are the components of the canonical
morphism $N_{(f,i)}^\vee \rightarrow (f,i)^* \Omega_{X \times S}$.
Define the complex $L_{(f,i)}$ to be
\begin{equation}
\begin{CD}
-1 & & 0 \\
N_{(f,i)}^\vee @> \beta_i >> i^* \Omega_S
\end{CD}
\end{equation}
There is a map of complexes $\gamma_{(f,i)}:L_{(f,i)} \rightarrow
L_\zeta$ defined by the commutative diagram
\begin{equation}
\begin{CD}
N_{(f,i)}^\vee @> \beta_i >> i^* \Omega_S \\
@V \alpha_f VV  @VV (di)^\dagger V \\
f^* \Omega_X @> (df)^\dagger >> \Omega_B
\end{CD}
\end{equation}

\begin{lem} \label{lem-qism}
  The morphism $\gamma_{(f,i)}: L_{(f,i)} \rightarrow L_\zeta$ is a
  quasi-isomorphism of complexes.
\end{lem}
 
\begin{proof} 
  This is an easy local argument which is left as an exercise for the
  reader.
\end{proof}

As a corollary of the lemma, we see that $L_\zeta^\vee$ is represented
in $D(\OO_{B})$ by the complex $L_{(f,i)}^\vee$ defined to be
\begin{equation}
\begin{CD}
0 & & 1 \\
i^* T_S @> (\beta_i)^\vee >> N_{(f,i)}
\end{CD}
\end{equation}

\

Now suppose that $p\in B \subset S$ is a point, define
$g:S'\rightarrow S$ to be the blow up of $S$ at $p$ with exceptional
divisor $E$, and define $i':B' \rightarrow S'$ to be the reduced total
transform of $B$.

\

We break up our analysis according to the type of behavior of $p\in
B$.  If $p\in B$ is a smooth point, we call this \emph{type (I)}.  If
$p\in B$ is a node, we call this \emph{type (II)}.  We further
decompose each type as follows.  If $p\in B$ is a smooth point which
lies on a stable component, we say this is \emph{type (Ia)}.  If $p\in
B$ is a smooth point which lies on an unstable component, we say this
is \emph{type (Ib)}.  If $p\in B$ is a node, and the first order
deformations of $\zeta$ smooth the node, we say this is \emph{type
  (IIa)}.  If the first order deformations of $\zeta$ don't smooth the
node, we say this is \emph{type (IIb)}.

\

For type (I), we have $g^* B$ equals $B'$ as Cartier divisors.  For
type (II), we have $g^* B$ equals $B' + E$ as Cartier divisors.  For
both types, we define $h:B' \rightarrow B$ to be the restriction of
$g$ and we define $f':B' \rightarrow B$ to be $f' = f\circ h$.

\

Of course we have $g_*\OO_{S'} = \OO_S$ and $R^{k>0}g_* \OO_{S'}$ is
zero.  For type (I), we have $\OO_{S'}(-B') = g^* \OO_S(-B)$.  So by
the projection formula we have $g_*\OO_{S'}(-B') = \OO_S(-B)$ and
$R^{k>0}g_*\OO_{S'}(-B')$ is zero.  Also we have that $g_*\OO_{S'}(E)
= \OO_S$ and $R^{k>0}g_* \OO_{S'}(E)$ is zero.  So also for type (II)
we have $g_*\OO_{S'}(-B') = \OO_S(-B)$ and $R^{k>0}g_* \OO_{S'}(-B')$
is zero.

\

Using the resolution of $\OO_{B'}$
\begin{equation}
\begin{CD}
0 @>>> \OO_{S'}(-B') @>>> \OO_{S'} @>>> \OO_{B'} @>>> 0
\end{CD}
\end{equation}
we conclude that $h_* \OO_{B'} = \OO_B$ and $R^{k>0}h_* \OO_{B'}$ is
zero.  In other words, the canonical morphism $\OO_B \rightarrow
\mathbb{R}h_* \OO_{B'}$ is a quasi-isomorphism.  From this it follows
by the projection formula that the canonical morphism
\begin{equation}
L_{(f,i)}^\vee \rightarrow 
\mathbb{R}h_* \mathbb{L}h^* (L_{(f,i)}^\vee)
\end{equation}
is a quasi-isomorphism.  Therefore the pullback morphisms
\begin{equation}
\mathbb{H}^k(B,L_\zeta^\vee) \rightarrow
\mathbb{H}^k(B',g^*L_\zeta^\vee)
\end{equation}
are isomorphisms.

\

Now there is a canonical morphism $\OO_{S'}(B') \rightarrow g^*
\OO_S(B)$.  For type (I), this morphism is an isomorphism.  For type
(II), this morphism is injective and the cokernel is $g^*\OO_S(B)|_E$,
i.e. $\OO_E \otimes_\CC M$ where $M= \OO_S(B)|_p$ is a one-dimensional
vector space.  For type (I), we conclude that the canonical morphism
$N_{(f',i')} \rightarrow h^* N_{(f,i)}$ is an isomorphism.  For type
(II), we conclude that there is an exact sequence:
\begin{equation}
0 \rightarrow \textit{Tor}_1^{\OO_{S'}}(\OO_{B'},\OO_E)\otimes_\CC M
\rightarrow N_{(f',i')} \rightarrow h^* N_{(f,i)} \rightarrow \OO_E
\otimes_\CC M \rightarrow 0
\end{equation}
For both types, denote by $\Gamma \subset B'$ the subcurve which is
the union of all components other than $E$, and let $D = E\cap
\Gamma$.  From the resolution of $\OO_E$
\begin{equation}
\begin{CD}
0 @>>> \OO_{S'}(-E) @>>> \OO_S @>>> \OO_E @>>> 0
\end{CD}
\end{equation}
we have the relation
\begin{equation}
\textit{Tor}_1^{\OO_{S'}}(\OO_{B'},\OO_E) = {\mc
  I}_{\Gamma}\otimes_{\OO_{S'}}\OO_{S'}(-E) = \OO_E(-E)(-D).
\end{equation}
where ${\mc I}_{\Gamma}$ is the ideal sheaf of $\OO_{B'}$ defining
$\Gamma \subset B'$, i.e. $\OO_E(-D)$.  So in case $(2)$, we have an
exact sequence:
\begin{equation}
0 \rightarrow \OO_E(-E)(-D)\otimes_\CC M \rightarrow N_{(f',i')}
\rightarrow h^*N_{(f,i)} \rightarrow \OO_E \otimes_\CC M \rightarrow
0.
\end{equation}

\

For both types, we have a short exact sequence of $\OO_{S'}$-modules
\begin{equation}
\begin{CD}
0 @>>> g^* \Omega_S @> (dg)^\dagger >> \Omega_{S'} @>>> \Omega_E @>>>
0.
\end{CD}
\end{equation}
Applying $\mathbb{R}\textit{Hom}_{\OO_{S'}}(*,\OO_{S'})$, we have an
exact sequence
\begin{equation}
\begin{CD}
0 @>>> T_{S'} @>>> g^* T_S @>>>
\textit{Ext}^1_{\OO_{S'}}(\OO_E,\OO_{S'}) @>>> 0 
\end{CD}
\end{equation}
Using the resolution of $\OO_E$ in the last paragraph, we have that
\begin{equation}
\textit{Ext}^1_{\OO_{S'}}(\OO_E,\OO_{S'}) = \OO_E(E).
\end{equation}
So we have an exact sequence
\begin{equation}
\begin{CD}
0 @>>> T_{S'} @>>> g^* T_S @>>> T_E(E) @>>> 0
\end{CD}
\end{equation}
Using our \emph{Tor} result from the last paragraph, we have a short
exact sequence
\begin{equation}
0 \rightarrow T_E(-D) \rightarrow (i')^*T_{S'} \xrightarrow{dg} h^*
i^* T_S \rightarrow T_E(E) \rightarrow 0.
\end{equation}
Notice this holds in both cases.

\

The maps $N_{(f',i')} \rightarrow h^* N_{(f,i)}$ and $(i')^*T_{S'}
\rightarrow h^* i^* T_S$ considered in the last two paragraphs are
compatible with $\alpha_{i'}$ and $h^*\alpha_i$.  So we have an
induced map of complexes $dg: L_{\zeta'}^\vee \rightarrow h^*
L_{\zeta}^\vee$.  Define $I\hookrightarrow h^* L_\zeta^\vee$ to be the
image complex of $dg:L_{\zeta'}^\vee \rightarrow h^* L_{\zeta}^\vee$.
For type (I), we define two complexes of coherent sheaves on $B'$,
$K_I$ and $Q_I$, by $K_I = T_E(-D)[0]$ and $Q_I = T_E(E)[0]$.  For
type (II), we define complexes of coherent sheaves on $B'$, $K_{II}$
and $Q_{II}$ by
\begin{eqnarray}
K_{II} = T_E(-D)[0] \oplus\lt( \OO_E(-E)(-D)\otimes_\CC M\rt)[-1] \\
Q_{II} = T_E(E)[0] \oplus \lt(\OO_E \otimes_\CC M\rt)[-1]
\end{eqnarray}
Then we have exact sequences of complexes
\begin{equation}
\begin{CD}
0 @>>> K @>>> L_{\zeta'}^\vee @> dg >> I @>>> 0 \\
0 @>>> I @>>> h^* L_{\zeta}^\vee @>>> Q @>>> 0
\end{CD}
\end{equation}

\

For type (I), we have that $\mathbb{H}^0(B',K_I) = H^0(E,T_E(-D))$ is
$2$-dimensional, because $T_E(-D) \cong \OO_E(1)$.  And $\mathbb{H}^{k
  > 0}(B',K_1)$ is zero.  Similarly $\mathbb{H}^0(B',Q_I) =
H^0(E,T_E(E))$ is $2$-dimensional and $\mathbb{H}^{k>0}(B',Q_I)$ is
zero.  Therefore we have a long exact sequence of hypercohomology
groups:
\begin{eqnarray*}
0 \rightarrow H^0(E,T_E(-D)) \rightarrow
\mathbb{H}^0(B',L_{(f',i')}^\vee) \rightarrow 
\mathbb{H}^0(B, L_{(f,i)}^\vee) \rightarrow \dots \\ 
\dots  \rightarrow  H^0(E, T_E(E)) \rightarrow
\mathbb{H}^1(B',L_{(f',i')}^\vee) \rightarrow
 \mathbb{H}^1(B, L_{(f,i)}^\vee ) \rightarrow 0 \\
  0 \rightarrow \mathbb{H}^2(B',L_{(f',i')}^\vee) \rightarrow \mathbb{H}^2(B,
 L_{(f,i)}^\vee ) \rightarrow 0
\end{eqnarray*}

\

For type (Ia), the map
\begin{equation}
\mathbb{H}^0(B,L_{(f,i)}^\vee) \rightarrow H^0(E,T_E(E))
\end{equation}
is the zero map.  So we have proved the following lemma.

\begin{lem} \label{lem-def1a}
  Suppose that $h:B'\rightarrow B$ is type (Ia).  Then we have exact
  sequences:
\begin{eqnarray}
0 \rightarrow H^0(E,T_E(-D)) \rightarrow
\mathbb{H}^0(B',L_{(f',i')}^\vee) \rightarrow 
\mathbb{H}^0(B, L_{(f,i)}^\vee) \rightarrow 0 \\ 
0 \rightarrow  H^0(E, T_E(E)) \rightarrow
\mathbb{H}^1(B',L_{(f',i')}^\vee) \rightarrow
 \mathbb{H}^1(B, L_{(f,i)}^\vee ) \rightarrow 0 \\
  0 \rightarrow \mathbb{H}^2(B',L_{(f',i')}^\vee) \rightarrow \mathbb{H}^2(B,
 L_{(f,i)}^\vee ) \rightarrow 0
\end{eqnarray}
So the canonical map from the Lie algebra of infinitesimal
automorphisms of $\zeta'$ to the Lie algebra of infinitesimal
automorphisms of $\zeta$ is surjective with $2$-dimensional kernel,
the canonical map from the space of first order deformations of
$\zeta'$ to the space of first order deformations of $\zeta$ is
surjective with $2$-dimensional kernel, and the obstruction space of
$\zeta'$ equals the obstruction space of $\zeta$.
\end{lem}

\

For type (Ib), the map
\begin{equation}
\mathbb{H}^0(B,L_{(f,i)}^\vee) \rightarrow H^0(E,T_E(E))
\end{equation}
has a $1$-dimensional image; we will call it $N$.  We have proved the
following lemma.

\begin{lem} \label{lem-def1b}
  Suppose that $h:B'\rightarrow B$ is type (Ib).  Then there is a
  $1$-dimensional subspace $N\subset H^0(E,T_E(E))$ such that we have
  exact sequences:
\begin{eqnarray}
0 \rightarrow H^0(E,T_E(-D)) \rightarrow
\mathbb{H}^0(B',L_{(f',i')}^\vee) \rightarrow 
\mathbb{H}^0(B, L_{(f,i)}^\vee) \rightarrow N \rightarrow 0 \\ 
0 \rightarrow  H^0(E, T_E(E))/N \rightarrow
\mathbb{H}^1(B',L_{(f',i')}^\vee) \rightarrow
 \mathbb{H}^1(B, L_{(f,i)}^\vee ) \rightarrow 0 \\
  0 \rightarrow \mathbb{H}^2(B',L_{(f',i')}^\vee) \rightarrow 
\mathbb{H}^2(B, L_{(f,i)}^\vee ) \rightarrow 0
\end{eqnarray}
So the canonical map from the Lie algebra of infinitesimal
automorphisms of $\zeta'$ to the Lie algebra of infinitesimal
automorphisms of $\zeta$ has both $1$-dimensional kernel and cokernel,
the canonical map from the space of first order deformations of
$\zeta'$ to the space of first order deformations of $\zeta$ is
surjective with $1$-dimensional kernel, and the obstruction space of
$\zeta'$ equals the obstruction space of $\zeta$.
\end{lem}

\

Next we consider type (II).  Then $\OO_E(-E)(-D)$ is isomorphic to
$\OO_E(-1)$.  Since $H^0(E,\OO_E(-1)) = H^1(E,\OO_E(-1)) = 0$, the
terms $\OO_E(-E)(-D)\otimes_\CC M[-1]$ do not contribute to the
hypercohomology of $K_{II}$.  And $T_E(-D)$ is isomorphic to $\OO_E$.
So we have $\mathbb{H}^0(B',K_{II}) = H^0(E,T_E(-D))$ is
$1$-dimensional and $\mathbb{H}^{k>0}(B',K_{II})$ is zero.

\

For $Q_{II}$ both terms contribute to the cohomology.  We have
$T_E(E)$ is isomorphic to $\OO_E(1)$ so that $\mathbb{H}^0(B',Q_{II})
= H^0(E,T_E(E))$ is $2$-dimensional, and $\mathbb{H}^1(B',Q_{II}) = M
= \OO_S(B)|_p$ is $1$-dimensional.  Therefore we have a long exact
sequence in hypercohomology:
\begin{eqnarray*}
0 \rightarrow H^0(E,T_E(-D)) \rightarrow
\mathbb{H}^0(B',L_{(f',i')}^\vee) \rightarrow
\mathbb{H}^0(B,L_{(f,i)}^\vee) \rightarrow \dots \\
\dots \rightarrow H^0(E,T_E(E)) \rightarrow
\mathbb{H}^1(B',L_{(f',i')}^\vee) \rightarrow
\mathbb{H}^1(B,L_{(f,i)}^\vee) \rightarrow \dots \\
\dots \rightarrow M \rightarrow \mathbb{H}^2(B',L_{(f',i')}^\vee)
\rightarrow \mathbb{H}^2(B,L_{(f,i)}^\vee) \rightarrow 0
\end{eqnarray*}

\

Geometrically, every infinitesimal automorphism of $\zeta$ lifts to an
infinitesimal automorphism of $\zeta'$, so that
$\mathbb{H}^0(B,L_{(f,i)}^\vee)\rightarrow H^0(E,T_E(E))$ is the zero
map.  Similarly, the map $\mathbb{H}^1(B,L_{(f,i)}^\vee) \rightarrow
\OO_S(B)|_p$ is nonzero iff there are deformations of $\zeta$ which
smooth the node $p$ to first order.  So we have the following lemma:

\begin{lem} \label{lem-def2a}
  Suppose that $h:B'\rightarrow B$ is type (IIa).  Then we have exact
  sequences:
\begin{eqnarray}
0 \rightarrow H^0(E,T_E(-D)) \rightarrow
\mathbb{H}^0(B',L_{(f',i')}^\vee) \rightarrow
\mathbb{H}^0(B,L_{(f,i)}^\vee) \rightarrow 0 \\
0 \rightarrow H^0(E,T_E(E)) \rightarrow
\mathbb{H}^1(B',L_{(f',i')}^\vee) \rightarrow
\mathbb{H}^1(B,L_{(f,i)}^\vee) \rightarrow \OO_S(B)|_p \rightarrow 0 \\
0 \rightarrow
\mathbb{H}^2(B',L_{(f',i')}^\vee) \rightarrow \mathbb{H}^2(B,
L_{(f,i)}^\vee) \rightarrow 0 
\end{eqnarray}
So the canonical map from the space of first order deformations of
$\zeta'$ to the space of first order deformations of $\zeta$ has both
$1$-dimensional kernel and cokernel, and the obstruction space to
$\zeta'$ equals the obstruction space to $\zeta$.
\end{lem}

\begin{lem} \label{lem-def2b}
  Suppose that $h:B'\rightarrow B$ is type (IIb).  Then we have exact
  sequences:
\begin{eqnarray}
0 \rightarrow H^0(E,T_E(-D)) \rightarrow
\mathbb{H}^0(B',L_{(f',i')}^\vee) \rightarrow
\mathbb{H}^0(B,L_{(f,i)}^\vee) \rightarrow 0 \\
0 \rightarrow H^0(E,T_E(E)) \rightarrow
\mathbb{H}^1(B',L_{(f',i')}^\vee) \rightarrow
\mathbb{H}^1(B,L_{(f,i)}^\vee) \rightarrow 0 \\
0 \rightarrow \OO_S(B)|_p \rightarrow
\mathbb{H}^2(B',L_{(f',i')}^\vee) \rightarrow \mathbb{H}^2(B,
L_{(f,i)}^\vee) \rightarrow 0 
\end{eqnarray}
So the canonical map from the space of first order deformations of
$\zeta'$ to the space of first order deformations of $\zeta$ is
surjective with $1$-dimensional kernel, and the map from the
obstruction space of $\zeta'$ to the obstruction space of $\zeta$ is
surjective and has a $1$-dimensional kernel.
\end{lem}

\

Finally, we consider the case when $h:B'\rightarrow B$ is the identity
map, but there is one extra marked point $q\in B'$ which is not in
$B$; we call this \emph{type (III)}.  We further break this up as
follows.  If $q\in B$ lies on an unstable component, we call this
\emph{type (IIIa)}.  If $q\in B$ lies on a stable component, we call
this \emph{type (IIIb)}.

\
  
For type (III), there is a canonical short exact sequence of
complexes:
\begin{equation}
\begin{CD}
0 @>>> L_\zeta @>>> L_\zeta' @>>> \Omega_B(q)|_q[0] @>>> 0
\end{CD}
\end{equation}
Of course we have:
\begin{eqnarray}
\EExt{B}^0(\Omega_B(q)|_q[0],\OO_B) = 0, \\
\EExt{B}^1(\Omega_B(q)|_q[0],\OO_B) = T_B|_q, \\
\EExt{B}^2(\Omega_B(q)|_q[0],\OO_B) = 0.
\end{eqnarray}
The induced map $\mathbb{E}\text{xt}^0_{\OO_B}(L_\zeta,\OO_B)
\rightarrow T_B|_q$ is nonzero iff $q$ lies on an unstable component
of $\zeta$.  Thus we have the following lemma:

\begin{lem} \label{lem-def3a}
  Suppose that $h:B'\rightarrow B$ is type (IIIa).  Then we have exact
  sequences:
\begin{eqnarray}
0 \rightarrow \EExt{B'}^0(L_{\zeta'},\OO_{B'}) \rightarrow
\EExt{B}^0(L_\zeta,\OO_B) \rightarrow T_B|_q \rightarrow 0 \\
0 \rightarrow \EExt{B'}^1(L_{\zeta'},\OO_{B'}) \rightarrow
\EExt{B}^1(L_\zeta,\OO_B) \rightarrow 0 \\
0 \rightarrow \EExt{B'}^2(L_{\zeta'},\OO_{B'}) \rightarrow
\EExt{B}^2(L_\zeta,\OO_B) \rightarrow 0 
\end{eqnarray}
So the Lie algebra of infinitesimal automorphisms of $\zeta'$ has
codimension $1$ in the Lie algebra of infinitesimal automorphisms of
$\zeta$, the space of first order deformations of $\zeta'$ equals the
space of first order deformations of $\zeta$, and the obstruction
space of $\zeta'$ equals the obstruction space of $\zeta$.
\end{lem}

\begin{lem} \label{lem-def3b}
  Suppose that $h:B'\rightarrow B$ is type (IIIb).  Then we have exact
  sequences:
\begin{eqnarray}
0 \rightarrow \EExt{B'}^0(L_{\zeta'},\OO_{B'}) \rightarrow
\EExt{B}^0(L_\zeta,\OO_B) \rightarrow 0 \\
0 \rightarrow T_B|_q \rightarrow \EExt{B'}^1(L_{\zeta'},\OO_{B'}) 
\rightarrow \EExt{B}^1(L_\zeta,\OO_B) \rightarrow 0 \\
0 \rightarrow \EExt{B'}^2(L_{\zeta'},\OO_{B'}) \rightarrow
\EExt{B}^2(L_\zeta,\OO_B) \rightarrow 0 
\end{eqnarray}
So the Lie algebra of infinitesimal automorphisms of $\zeta'$ equals
the Lie algebra of infinitesimal automorphisms of $\zeta$, the
canonical map from the space of first order deformations of $\zeta'$
to the space of first order deformations of $\zeta$ is surjective with
$1$-dimensional kernel, and the obstruction space of $\zeta'$ equals
the obstruction space of $\zeta'$.
\end{lem}

Combining lemma~\ref{lem-def1a} through lemma~\ref{lem-def3b}, one can
analyze the associated maps of vector spaces
$\EExt{B'}^k(L_{\zeta'},\OO_{B'}) \rightarrow
\EExt{B}^k(L_\zeta,\OO_B)$ for any morphism $h:B'\rightarrow B$ which
removes some subset of marked points from $B'$ and then contracts some
subset of the unstable components.

\subsection{Gluing stable curves}\label{subsec-glue}

Just as one has an obstruction theory for $\Kgnb{g,r}{X,\beta}$ of the
form $\lt( {\mathbb R}\pi_*\lt(L_\zeta^\vee\rt)[1]\rt)^\vee$, also for
any stable $A$-graph $\tau$, one has an analogous obstruction theory
for each of the Behrend-Manin stacks $\Kbm(X,\tau)$ (c.f.~\cite{BM}
for the definition of $\Kbm(X,\tau)$).  We will not describe this
obstruction theory here.  In \cite{B}, a \emph{relative obstruction
  theory} for the morphism $\Kbm(X,\tau) \rightarrow {\mathfrak
  M}(\tau)$ is given from which an absolute obstruction theory for
$\Kbm(X,\tau)$ can be deduced.

\

Suppose that $\tau$ is a stable $A$-graph, and suppose that
$\{f_1,f_2\}$ is a disconnecting edge of $\tau$.  Let $\tau_1\subset
\tau$ be the maximal connected subgraph which contains $f_1$ and not
$f_2$, and let $\tau_2\subset \tau$ be the maximal connected subgraph
which contains $f_2$ and not $f_1$.  So we have forgetful
$1$-morphisms $\Kbm(X,\tau) \rightarrow \Kbm(X,\tau_i)$ for $i=1,2$.

\begin{lem}\label{lem-glue1}
  Suppose that $\zeta:T \rightarrow \Kbm(X,\tau)$ is a $1$-morphism
  and let $\zeta_i:T \rightarrow \Kbm(X,\tau_i)$, $i=1,2$ be the
  composition of $\zeta$ with the forgetful $1$-morphism above.
  Suppose that for each point $t\in T$, the stable maps $\zeta_1(t)\in
  \Kbm(X,\tau_1)$ and $\zeta_2(t) \in \Kbm(X,\tau_2)$ are unobstructed
  (in the sense that the obstruction groups described above are zero)
  and the evaluation morphism $\text{ev}_{f_1}:\Kbm(X,\tau_1)
  \rightarrow X$ is smooth at $\zeta_1(t)$.  Let $T_{\text{ev}_{f_1}}$
  denote the dual of the sheaf of relative differentials of
  $\text{ev}_{f_1}$.  Then also $\zeta(t)\in \Kbm(X,\tau)$ is
  unobstructed, and there is a short exact sequence:
\begin{equation}
\begin{CD}
0 @>>> \zeta_1^*T_{\text{ev}_{f_1}} @>>> \zeta^* T_{\Kbm(X,\tau)} @>>>
\zeta_2^* T_{\Kbm(X,\tau_2)} @>>> 0
\end{CD}
\end{equation}
\end{lem}

\begin{proof}
  The proof essentially follows from the fact that $\Kbm(X,\tau)$ is
  an open substack of the $2$-fiber product:
\begin{equation}
\Kbm(X,\tau_1) \times_{\text{ev}_{f_1},X,\text{ev}_{f_2}}
\Kbm(X,\tau_2).
\end{equation}
The details are left to the reader.
\end{proof}

Now suppose that $\phi:\tau \rightarrow \sigma$ is the contraction of
stable $A$-graphs which contracts the edge $\{f_1,f_2\}$.  The induced
$1$-morphism $\Kbm(X,\tau) \rightarrow \Kbm(X,\sigma)$ is unramified
of codimension at most $1$.  In some circumstances, it is the
normalization of a Cartier divisor.

\begin{lem}\label{lem-glue2}
  With the same notation as in lemma~\ref{lem-glue1}, suppose that
  $\tau$ is a genus $0$ tree.  For $i=1,2$, let the domain of
  $\zeta_i$ be given by $\pi_i:C_i \rightarrow T$, let $g_i:C_i
  \rightarrow X$ be the map of $\zeta_i$, and let $s_i:T \rightarrow
  C_i$ be the section corresponding to the flag $f_i$ of $\tau_i$.
  Denote by $T_{\pi_i}$ the dual of the sheaf of relative
  differentials of $\pi_i$.
  
  Suppose that for every point $t\in T$, $\Kbm(X,\tau_2)$ is
  unobstructed at $\zeta_2(t)$, and suppose that $g_1^*T_X$ is
  $\pi_1$-relatively generated by global sections.  Then for each
  point $t\in T$, $\Kbm(X,\sigma)$ is unobstructed at $\zeta(t)$, the
  morphism $\Kbm(X,\tau) \rightarrow \Kbm(X,\sigma)$ is a regular
  embedding of codimension $1$ at $\zeta(t)$, and we have a short
  exact sequence:
\begin{equation}
\begin{CD}
0 @>>> \zeta^*T_{\Kbm(X,\tau)} @>>> \zeta^* T_{\Kbm(X,\sigma)} @>>>
s_1^*T_{\pi_1} \otimes s_2^* T_{\pi_2} @>>> 0
\end{CD}
\end{equation}
\end{lem}

\begin{proof}
  Let $\pi:C\rightarrow T$ be the family of curves obtained by
  identifying the section $s_1$ of $\pi_1$ and the section $s_2$ of
  $\pi_2$.  Let $s:T \rightarrow C$ be the section corresponding to
  $s_1$ and $s_2$.  Let $g:C \rightarrow X$ be the map obtained by
  gluing $g_1$ and $g_2$.

\

First we prove that for any point $t\in T$ we have that
$\Kbm(X,\sigma)$ is unobstructed at $\zeta(t)$ and the morphism
$\Kbm(X,\tau) \rightarrow \Kbm(X,\sigma)$ is a regular embedding of
codimension $1$.  The two statements together are equivalent to the
statement that there are first order deformations of the map
$\zeta(t)\in \Kbm(X,\sigma)$ which smooth the node $s(t)\in C_t$.

\

Now there is an exact sequence for $\zeta(t)\in \Kbm(X,\sigma)$:
\begin{equation}
T_{\Kbm(X,\sigma)}|_{\zeta(t)} \rightarrow \text{Ext}^1\lt(
\Omega_{C_t},\OO_{C_t} \rt) \rightarrow H^1\lt( C_t,g^*T_X \rt) \rightarrow
\text{Obs}(\zeta(t))  
\rightarrow 0
\end{equation}
Here $\text{Obs}(\zeta)$ is the obstruction group to $\Kbm(X,\sigma)$
at $\zeta$.  Similarly, we have an exact sequence for $\zeta_2(t)$:
\begin{equation}
\begin{CD}
T_{\Kbm(X,\tau_2)}|_{\zeta_2(t)} @>>> \text{Ext}^1\lt(
\Omega_{(C_2)_t} \rt) @>>> H^1\lt( (C_2)_t, g_2^*T_X \rt) @>>> 0
\end{CD}
\end{equation}

\

We have a short exact sequence of sheaves on $C_t$:
\begin{equation}
\begin{CD}
  0 @>>> g_1^*T_X\lt(-s(t)\rt) @>>> g^*T_X @>>> g_2^*T_X @>>> 0
\end{CD}
\end{equation}
By assumption, $g_1^*T_X$ is generated by global sections.  Thus by
lemma~\ref{lem-gend}, we conclude that $H^1\lt(C_t, g_1^*T_X\lt(-s(t)
\rt) \rt)$ is zero.  Thus we have an identification of $H^1\lt(
C_t,g^*T_X \rt)$ and $H^1\lt( (C_2)_t, g_2^*T_X \rt)$.  Also,
$\text{Ext}^1(\Omega_{C_t},\OO_{C_t})$ is canonically isomorphic to
the product over all nodes $q\in C_t$ of $T'_q\otimes T''_q$ where
$T'_q$ and $T''_q$ are the tangent spaces to the two branches of $C_t$
at $q$.  We have the analogous result for $(C_2)_t$.  Via these
identifications, we have a commutative diagram:
\begin{equation}
\begin{CD}
\text{Ext}^1\lt( \Omega_{C_t},\OO_{C_t} \rt) @>>> H^1\lt( C_t, g^* T_X
\rt) \\
@V p VV @VV = V \\
\text{Ext}^1\lt( \Omega_{(C_2)_t}, \OO_{(C_2)_t} \rt) @>>> H^1\lt(
(C_2)_t, g_2^* T_X \rt)
\end{CD}
\end{equation}
where $p$ is the canonical projection.  Choose any section
\begin{equation}
\phi:T'_{s(t)}\otimes T''_{s(t)} \rightarrow \text{Ext}^1\lt(
\Omega_{C_t}, \OO_{C_t} \rt)
\end{equation}
and any section
\begin{equation}
\psi:\text{Ext}^1\lt( \Omega_{(C_2)_t}, \OO_{(C_2)_t} \rt) \rightarrow
\text{Ext}^1\lt( \Omega_{C_t},\OO_{C_t} \rt)
\end{equation}
of the canonical projections.  Choose any element $u\in
T'_{s(t)}\otimes T''_{s(t)}$ and consider the image
$\overline{\phi(u)}$ of $\phi(u)$ in $H^1\lt( C_t, g^* T_X \rt)$.
Since the map
\begin{equation}
\begin{CD}
\text{Ext}^1\lt( \Omega_{(C_2)_t}, \OO_{(C_2)_t} \rt) @>>> H^1\lt(
(C_2)_t, g_2^* T_X \rt)
\end{CD}
\end{equation}
is surjective, we can find some element $v\in \text{Ext}^1\lt(
\Omega_{(C_2)_t}, \OO_{(C_2)_t} \rt)$ such that $\overline{\psi(v)}$
equals $\overline{\phi(u)}$.  Consider $\phi(u)-\psi(v)$.  This has
image $0$ in $H^1\lt( C_t, g^*T_X \rt)$.  Therefore we conclude that
it is in the image of $T_{\Kbm(X,\sigma)}|_{\zeta(t)}$.  So we
conclude that $u$ is in the image of the projection map
$T_{\Kbm(X,\sigma)} \rightarrow T'_{s(t)}\otimes T''_{s(t)}$, i.e.
this projection map is surjective.  So the deformations of $\zeta(t)$
smooth the node $s(t)$.  Therefore $\Kbm(X,\sigma)$ is smooth at
$\zeta(t)$ and the morphism $\Kbm(X,\tau) \rightarrow \Kbm(X,\sigma)$
is a regular embedding of codimension $1$ at $\zeta(t)$.

\

Finally, the short exact sequence above is just the globalized version
of the projection map $T_{\Kbm(X,\sigma} \rightarrow T'_{s(t)}\otimes
T''_{s(t)}$ appearing in the last paragraph.
\end{proof}

\section{Properties of families of stable maps} \label{sec-props}

In this section we introduce some definitions and lemmas regarding
properties of families of stable maps.  Recall, to prove that
$\Kgnb{0,0}{X,e}$ is rationally connected, we have to find a
\emph{very free} $1$-morphism $\zeta:\PP^1 \rightarrow
\Kgnb{0,0}{X,e}$, i.e. a $1$-morphism whose image is contained in the
smooth locus and such that $\zeta^*T_{\Kgnb{0,0}{X,e}}$ is an ample
vector bundle.  Our proof that such a $1$-morphism exists is by
induction, where the induction step consists of attaching a
$1$-parameter family of lines to our $1$-parameter family of degree
$e$ stable maps.  To make the induction argument work, we need a bit
more than a very free $1$-morphism $\zeta:\PP^1 \rightarrow
\Kgnb{0,0}{X,e}$.  The property we need is what we call a \emph{very
  positive} $1$-morphism.  Additionally, we need that our
$1$-parameter family of lines has a property which we call \emph{very
  twisting}.  Finally, because of an operation we perform on
$1$-morphisms which we call \emph{modification}, and which we
introduce in the next section, we need to consider the case of
$1$-morphisms with reducible domain, i.e. $\zeta:B \rightarrow
\Kgnb{0,1}{X,e}$ where $B$ is a connected, prestable curve of
arithmetic genus $0$.  This is the level of generality in which we
make all our definitions.

\begin{defn}  Given a genus $0$ stable map
  $\zeta = ((B,p_1,\dots,p_n),f:B\rightarrow X)$, we say $\zeta$ is
  \emph{very stable} if the unmarked prestable map $(B,f:B\rightarrow
  X)$ is stable.
\end{defn}

\begin{notat} \label{not-1mor}
  Given a closed subscheme $X\subset \PP^n$ and a scheme $B$, a
  $1$-morphism $\zeta:B\rightarrow \Kgnb{g,r}{X,e}$ is equivalent to a
  datum:
\begin{equation}
\zeta = \lt( \lt(p_\zeta:\Sigma_\zeta \rightarrow B,
\sigma_{\zeta,1},\dots, \sigma_{\zeta,r} \rt), g_\zeta \rt).
\end{equation}
Here $p_\zeta:\Sigma_\zeta\rightarrow B$ is a family of prestable
curves, $\sigma_{\zeta,i}:B\rightarrow \Sigma_\zeta$ is a collection
of sections, $g_\zeta:\Sigma_\zeta\rightarrow X$ is a morphism of
schemes and we denote $h_{\zeta,i}=g_\zeta\circ \sigma_{\zeta,i}$.
When there is no risk of confusion, we will suppress the $\zeta$
subscripts.
\end{notat}

\begin{defn}~\label{defn-twisting}
  Suppose $\pi:B\rightarrow T$ is a family of prestable, geometrically
  connected curves of arithmetic genus $0$.  Suppose given a
  $1$-morphism $\zeta:B \rightarrow \Kgnb{0,1}{X,1}$, i.e. a datum
\begin{equation}
\zeta = \lt( p:\Sigma \rightarrow B, \sigma: B \rightarrow \Sigma,
g:\Sigma \rightarrow X \rt)
\end{equation}
such that $X$ is smooth along $g(\Sigma)$.  The $1$-morphism
$\zeta:B\rightarrow \Kgnb{0,1}{X,1}$ is \emph{twisting} if
\begin{enumerate}
\item The data $(\pi:B\rightarrow T, h:B \rightarrow X)$ is a
  family of stable maps to $X$, i.e. a $1$-morphism $\xi:T\rightarrow
  \Kgnb{0,0}{X,e}$ for some $e\geq 0$.
\item The image of $\xi:T\rightarrow \Kgnb{0,0}{X,e}$ is contained in
  the \emph{very unobstructed} locus of $\Kgnb{0,0}{X}$.
\item The image of $\zeta:T \rightarrow \Kgnb{0,1}{X,1}$ is contained
  in the \emph{very unobstructed} locus of the evaluation morphism
  $\text{ev}:\Kgnb{0,1}{X,1} \rightarrow X$.
\item Denoting by $T_{\text{ev}}$ the dual of the sheaf of relative
  differentials $\Omega_{\text{ev}}$, the pullback bundle
  $\zeta^*T_{\text{ev}}$ is
  $\pi$-relatively generated by global sections.
\item Denoting by $\text{pr}:\Kgnb{0,1}{X,1}\rightarrow \Kgnb{0,0}{X,1}$
  the projection map, and by $T_{\text{pr}}$ the dual of the sheaf of
  relative differentials $\Omega_{\text{pr}}$, the pullback bundle
  $\zeta^*T_{\text{pr}}$, i.e. the line bundle
  $\sigma^*\OO_\Sigma(\sigma)$, is $\pi$-relatively generated by
  global sections.
\end{enumerate}
\end{defn}

\

\begin{defn}~\label{defn-verytwisting}
  With notation as in definition~\ref{defn-twisting}, a morphism
  $\zeta:T\rightarrow \Kgnb{0,1}{X,1}$ is \emph{very twisting} if it
  is twisting and if $\zeta^*T_{\text{ev}}$ is $\pi$-relatively
  deformation ample.
\end{defn}

\

\begin{rmk}~\label{rmk-twisting} 
  Regarding the definitions above:
\begin{enumerate}
  
\item In $(2)$ and $(3)$ of defintion~\ref{defn-twisting}, \emph{very
    unobstructed} means that the naive obstruction group vanishes.
  For $(2)$ this means that for each $t\in T$ and the corresponding
  stable map $(h_t:B_t\rightarrow X)$, the following group vanishes:
\begin{equation}
\begin{CD}
  & -1 &  & 0 \\
\mathbb{E}\text{xt}^1_{B_t}( & h_t^*\Omega_X @>>> \Omega_{B_t},\ & \OO_{B_t})
\end{CD}
\end{equation}

\item It is easy to see that $\zeta^* T_{\text{pr}}$ is just
  $\sigma^*\OO_{\Sigma}(\sigma)$.  
  
\item Observe the product morphism $(p,g):\Sigma\rightarrow B \times
  X$ is a regular embedding.  Denote by ${\mc N}$ the normal bundle of
  this regular embedding.  Then $(3)$ is equivalent to the condition
  that $R^1p_*\lt({\mc N}(-\sigma)\rt)$ is trivial.  In this case
  $\zeta^*T_{\text{ev}}$ is the locally free sheaf $p_*\lt({\mc
    N}(-\sigma)\rt)$.

\item Since the prestable family of maps $(\pi:B \rightarrow T, \xi:B
  \rightarrow X)$  is stable, clearly also $(\pi:B \rightarrow T,
  \zeta: B \rightarrow \Kgnb{0,1}{X,1})$ is stable.
  
\item There are some degree conditions implicit in these definitions.
  The total degree of $\sigma^*\OO_{\Sigma}(\sigma)$ is simply
  $s=2e-e'$ where $e$ is the degree of $h:B \rightarrow \PP^N$ and
  $e'$ is the degree of $g:\Sigma \rightarrow \PP^N$ (both degrees
  with respect to $\OO_{\PP^N}(1)$).  So if $\zeta$ is twisting, we
  have that $2e \geq e'$ and if $\zeta$ is very twisting, we have that
  $2e > e'$.
  
\item Additionally, given a twisting family $\zeta$ a point $b\in B$,
  and a deformation of the line $g_b:\Sigma_b \rightarrow X$ which
  continues to contain $h(b)$, there must be a deformation of the
  whole family $\zeta$ giving rise to the deformation of $\Sigma_b$
  and which does not deform $h:B \rightarrow X$.  In particular, if
  $h:B \rightarrow X$ is also an embedded line, then the map of the
  surface $g:\Sigma \rightarrow X$ must have degree $1$ or $2$ and
  must deform along with a line which intersects $h(B)$ in a fixed
  point.  If $X \subset \PP^n$ is a hypersurface of low degree $d > 1$
  such that to a general line $h:B \rightarrow X$ there is a
  corresponding twisting family $\zeta$ (what we refer to as a
  \emph{twistable line} below), then we must have that $h:B
  \rightarrow X$ is an embedding of a smooth quadric surface.  The
  condition on such $X$ that a general line is \emph{twistable} is
  essentially that, given two general intersecting lines $B$ and $L$
  in $X$, there is a smooth quadric surface $\Sigma$ in $X$ which
  contains both $B$ and $L$.  In a later section we will see that this
  condition does hold for a general hypersurface $X\subset \PP^n$ of
  degree $d$ when $d^2 \leq n+1$.

\end{enumerate}
\end{rmk}

\begin{lem} \label{lem-twisttogether}
  Suppose $B = B_1 \cup B_2$ is a prestable, geometrically connected
  curve of arithmetic genus $0$ where $B_1$ and $B_2$ are connected
  subcurves such that $B_1 \cap B_2$ is a single node of $B$.  Suppose
  given $\zeta:B \rightarrow \Kgnb{0,1}{X,1}$ such that $\zeta|_{B_i}
  : B_i \rightarrow \Kgnb{0,1}{X,1}$ is twisting for $i=1,2$.  Then
  $\zeta$ is twisting.  If, in addition, at least one of $\zeta_i$ is
  very twisting, then $\zeta$ is very twisting.
\end{lem}

\begin{proof}
  This is an easy consequence of lemma~\ref{lem-DAcrit}.
\end{proof}

\begin{lem} \label{lem-twistopen}
  Let $\pi:B\rightarrow T$ be a family of prestable, geometrically
  connected curves of arithmetic genus $0$ and let $\zeta:B
  \rightarrow \Kgnb{0,1}{X,1}$ be a morphism.  There is an open
  subscheme $U_{\text{twist}}\subset T$ (resp. $U_{vtwist}\subset T$)
  with the following property: for any morphism of schemes
  $f:T'\rightarrow T$, the pullback family $f^*\pi:f^*B \rightarrow
  T'$ and $f^*\zeta: f^*B \rightarrow \Kgnb{0,1}{X,1}$ is twisting
  (resp. very twisting) iff $f(T')\subset U_{\text{twist}}$ (resp.
  $f(T')\subset U_{\text{vtwist}}$).
\end{lem}

\begin{proof}
  By \cite[lemma 1]{B} there is an universal open subscheme
  $U_1\subset T$ over which $(\pi:B\rightarrow T, h:B\rightarrow X)$
  is a family of stable maps.  It is clear that $U$, if it exists,
  must also be contained in the complement of the support of
\begin{equation}
\mathbb{R}^1\lt(\pi_*\textit{Hom}_{\OO_B}\rt)\lt(h^*\Omega_X
\rightarrow \Omega_\pi, \OO_B\rt),
\end{equation}
and in the complement of the image under $\pi$ of the supports of the
sheaves:
\begin{eqnarray}
R^1p_*\lt({\mc N}(-\sigma)\rt), \\
\text{coker}\lt( \pi^*\pi_* \sigma^*\OO_{\Sigma}(\sigma) \rightarrow
\sigma^*\OO_{\Sigma}(\sigma) \rt).
\end{eqnarray}
Let $U_2$ denote the complement of these sets in $U_1$.  On $U_2$ all
of the conditions to be twisting (resp. very twisting) are satisfied
except the condition that $\zeta^*T_{\text{ev}}$ is $\pi$-relatively
generated by global sections (resp. $\pi$-relatively deformation
ample).  So we define $U_{\text{twist}}$ to be the complement in $U_2$
of the image under $\pi$ of the cokernel of the morphism
\begin{equation}
\pi^*\pi_* \zeta^*T_{\text{ev}}\rightarrow \zeta^*T_{\text{ev}}.
\end{equation}
And, using lemma~\ref{lem-opDA}, we define $U_{\text{vtwist}}$ to be
the universal open subscheme of $U_2$ over which
$\zeta^*T_{\text{ev}}$ is $\pi$-relatively deformation ample.  It
follows immediately from the construction that $U_{\text{twist}}$ and
$U_{\text{vtwist}}$ have the desired universal properties.
\end{proof}

\begin{defn} \label{defn-twistable}
  Suppose $(\pi:B\rightarrow T, h:B\rightarrow X)$ is a family of
  genus $0$ stable maps, i.e. a $1$-morphism $\xi:T\rightarrow
  \Kgnb{0,0}{X,e}$ for some $e \geq 0$.  We say $\xi:T \rightarrow
  \Kgnb{0,0}{X,e}$ is \emph{twistable} (resp. \emph{very twistable})
  if there exists a surjective \'etale morphism $u:T'\rightarrow T$
  and a morphism $\zeta:u^*B \rightarrow \Kgnb{0,1}{X,1}$ with
  $h_\zeta = u^* h$ such that $\zeta$ is twisting (resp. very
  twisting).
\end{defn}

\begin{prop} \label{prop-twistopen}
  Let $\xi = (\pi:B\rightarrow T, h:B\rightarrow X)$ be a $1$-morphism
  $\xi:T \rightarrow \Kgnb{0,0}{X,e}$.  There is an open subscheme
  $U_{t-able}\subset T$ (resp $U_{vt-able}\subset T$) such that for
  each morphism of schemes $f:T'\rightarrow T$, the pullback $(f^*\pi:
  f^*B \rightarrow T', f^*h: f^*B \rightarrow X)$ is twistable (resp.
  very twistable) iff $f(T') \subset U_{t-able}$ (resp. $f(T')\subset
  U_{vt-able}$).
\end{prop}

\begin{proof}
  It suffices to check that if $t_0\in T$ is a geometric point such
  that $h_{t_0}:B_{t_0} \rightarrow X$ is twistable (resp. very
  twistable), then there is an \'etale neighborhood of $t_0\in T$ over
  which $\xi$ is twistable (resp. very twistable).  Denote by
  $\zeta_0:B_{t_0} \rightarrow \Kgnb{0,1}{X,1}$ the twisting morphism.
  We consider $\Kgnb{0,1}{X,1}$ as a projective scheme via the
  Pl\"ucker and Segr\'e embeddings of $\GG(1,n) \times \PP^n
  \hookrightarrow \PP^{\frac{n(n+1)^2}{2}-1}$.  Let $\beta$ denote the
  degree of the stable map $\zeta_0$.

\

Define ${\mc M} = T \times \Kgnb{0,0}{\Kgnb{0,1}{X,1}, \beta}$, i.e.
${\mc M}$ parametrizes pairs $(t,\zeta)$ where $t\in T$ is a point and
where $\zeta:B \rightarrow \Kgnb{0,1}{X,1}$ is a genus $0$ stable map
of degree $\beta$.  Denote the universal stable map by
\begin{eqnarray}
\rho:{\mc B} \rightarrow \Kgnb{0,0}{\Kgnb{0,1}{X,1},\beta}, \\
\zeta: {\mc B} \rightarrow \Kgnb{0,1}{X,1}.
\end{eqnarray}
As in notation~\ref{not-1mor}, let $p:\Sigma \rightarrow {\mc B}$ be
the pullback by $\zeta$ of the universal curve over $\Kgnb{0,1}{X,1}$,
let $\sigma:{\mc B} \rightarrow \Sigma$ be the pullback of the
universal section, let $g:\Sigma \rightarrow X$ be the pullback of the
universal map, and let $h = g \circ \sigma$.  So we have a family of
prestable maps 
\begin{equation}
\widetilde{\xi} = \lt( \rho:{\mc B} \rightarrow \Kgnb{0,0}{\Kgnb{0,1}{X,1},
  \beta}, h: {\mc B} \rightarrow X \rt).
\end{equation}
By \cite[lemma 1]{B} there is a maximal open substack ${\mc U}_e \subset
\Kgnb{0,0}{\Kgnb{0,1}{X,1}, \beta}$ over which $\widetilde{\xi}$ is
stable of degree $e$.  By assumption, $(t_0,\zeta_0)$ is in $T \times
{\mc U}_e$.  

\

In the last paragraph we constructed a $1$-morphism
\begin{equation}
(1_T,\widetilde{\xi}): T \times {\mc U}_e \rightarrow T \times
\Kgnb{0,0}{X,e}.
\end{equation}
We also saw that $(t_0,\zeta_0)$ is in the domain of this
$1$-morphism.  The claim is that $(1_T,\widetilde{\xi})$ is smooth on
a neighborhood of $(t_0,\zeta_0)$, i.e. that $\widetilde{\xi}: {\mc
  U}_e \rightarrow \Kgnb{0,0}{X,e}$ is smooth at $\zeta_0$.  First we
will show that ${\mc U}_e$ is smooth at $\zeta_0$.  The space of first
order deformations and the obstruction space of ${\mc U}_e$ at
$\zeta_0$ are given by
\begin{equation}
\mathbb{E}\text{xt}^i_{B_{t_0}}(L^\cdot_{\zeta_0},\OO_{B_{t_0}}),
\end{equation}
for $i=1,2$ respectively, where $L^\cdot_{\zeta_0}$ is the complex
\begin{equation}
\begin{CD}
-1 & & 0 \\
\zeta_0^* \Omega_{\Kgnb{0,1}{X,1}} @> d(\zeta_0)^\dagger >> \Omega_{B_{t_0}}.
\end{CD}
\end{equation}
Now the induced morphism $\xi_0:B_{t_0} \rightarrow X$ by $\xi_0 =
\text{ev} \circ \zeta_0$ also has an associated complex
$L^\cdot_{\xi_0}$:
\begin{equation}
\begin{CD}
-1 & & 0 \\
\xi_0^* \Omega_X @> d(\xi_0)^\dagger >> \Omega_{B_{t_0}}.
\end{CD}
\end{equation}
There is a morphism of complexes:
\begin{eqnarray}
\gamma: L^\cdot_{\xi_0} \rightarrow L^\cdot_{\zeta_0} \\
\gamma^0 = \text{id} : \Omega_B \rightarrow \Omega_B \\
\gamma^{-1} = \zeta_0^*\lt( d(\text{ev})^\dagger \rt): \zeta_0^*
\text{ev}^* \Omega_X \rightarrow \zeta_0^* \Omega_{\Kgnb{0,1}{X,1}}.
\end{eqnarray}
There is also a morphism of complexes:
\begin{equation}
\delta: L^\cdot_{\zeta_0} \rightarrow \zeta_0^* \Omega_{\text{ev}}[1],
\end{equation}
where $\delta^{-1}:\zeta_0^* \Omega_{\Kgnb{0,1}{X,1}} \rightarrow
\zeta_0^*\Omega_{\text{ev}}$ is the pullback of the canonical
surjection.  And the triple:
\begin{equation}
\begin{CD}
L^\cdot_{\xi_0} @> \gamma >> L^\cdot_{\zeta_0} @> \delta >> \zeta_0^*
\Omega_{\text{ev}}[1] 
\end{CD}
\end{equation}
is an exact triangle.  Thus there is a corresponding long exact
sequence of $\mathbb{E}\text{xt}$'s.  Condition $(2)$ of
definition~\ref{defn-twisting} says that
$\mathbb{E}\text{xt}^1_{B_{t_0}}(L^\cdot_{\xi_0},\OO_{B_{t_0}})$ is
zero.  By condition $(4)$ of the definition, $\zeta_0^* T_{\text{ev}}$
is generated by global sections.  Since $B_{t_0}$ is connected of
arithmetic genus $0$, we have that $H^1(B_{t_0},\OO_{B_{t_0}})$ is
zero.  So for any trivial bundle, $H^1$ is zero.  Since $H^2$ vanishes
on all coherent sheaves, we conclude that for any sheaf generated by
global sections, $H^1$ is zero.  Thus we have that the group
\begin{equation}
\mathbb{E}\text{xt}^2(\zeta_0^*\Omega_{\text{ev}}, \OO_{B_{t_0}}) =
H^1(B_{t_0}, T_{\text{ev}}),
\end{equation}
is also zero.  By the long exact sequence, we conclude that
$\mathbb{E}\text{xt}^1_{B_{t_0}}( L^\cdot_{\zeta_0}, \OO_{B_{t_0}} )$
is also zero.  So the obstruction group vanishes and ${\mc U}_e$ is
smooth at $\zeta_0$.

\

By condition $(2)$ the image point $\xi_0\in \Kgnb{0,0}{X,e}$ is a
smooth point of $\Kgnb{0,0}{X,e}$.  Thus to prove that
$\widetilde{\xi}: {\mc U}_e \rightarrow \Kgnb{0,0}{X,e}$ is smooth, it
suffices to prove that derivative map $d(\widetilde{\xi})$ is
surjective on the space of first order deformations.  The map
\begin{equation}
d(\widetilde{\xi}):
\mathbb{E}\text{xt}^1_{B_{t_0}}(L^{\cdot}_{\zeta_0}, \OO_{B_{t_0}})
\rightarrow \mathbb{E}\text{xt}^1_{B_{t_0}}(L^\cdot_{\xi_0},
\OO_{B_{t_0}}),
\end{equation}
is precisely the map occurring in the long exact sequence of
$\mathbb{E}\text{xt}$'s from the paragraph above.  By the long exact
sequence, the cokernel of this map is a subgroup of
$H^1(B_{t_0},T_{\text{ev}})$, and this is zero as we have seen.
Therefore $\widetilde{\xi}$ is smooth at $\zeta_0\in {\mc U}_e$.  

\

Consider the morphism $(1_T,\xi): T\rightarrow T \times
\Kgnb{0,0}{X,e}$.  We can form the fiber product ${\mc M}$ of
$(1_T,\xi)$ with the morphism $(1_T,\widetilde{\xi}): T\times {\mc
  U}_e \rightarrow T \times \Kgnb{0,0}{X,e}$.  The fiber product ${\mc
  M}$ exactly parametrizes triples $(t,\zeta,\theta)$ where $t\in T$
is a point, $\zeta:B \rightarrow \Kgnb{0,1}{X,1}$ is a point in ${\mc
  U}_e$, and $\theta:\xi_t \rightarrow \widetilde{\zeta}$ is an
equivalence of objects in the groupoid $\Kgnb{0,0}{X,e}(\SP
\kappa(t))$.  The projection map $\text{pr}_1:{\mc M} \rightarrow T$
is smooth at $(t_0,\zeta_0)$ by the last paragraph.  So we can find an
\'etale morphism $f:M \rightarrow {\mc M}$ of a scheme to ${\mc M}$
whose image contains $(t_0,\zeta_0)$, and such that $M\rightarrow T$
is smooth.  Thus there is an \'etale morphism $u:T' \rightarrow T$ and
a section $z:T' \rightarrow M$.  Define $\zeta:T'\rightarrow {\mc
  U}_e$ to be the composition $\text{pr}_2\circ g \circ z$.

\

We also denote by $\zeta:B' \rightarrow \Kgnb{0,1}{X,1}$ the pullback
by $\zeta:T'\rightarrow {\mc U}$ of the universal stable map.  As
$\widetilde{\xi}(\zeta): B' \rightarrow \Kgnb{0,0}{X,e}$ is equivalent
to $u^* \xi: u^*B \rightarrow \Kgnb{0,0}{X,e}$, after replacing $T'$
by an \'etale, cover, we may suppose that $B'=u^*B$ as $T'$-schemes,
and $\widetilde{\xi}(\zeta) = u^*\xi$.  Now the fiber of
$\zeta:u^*B\rightarrow \Kgnb{0,1}{X,1}$ over any preimage of
$(t_0,\zeta_0)$ is twisting.  So by lemma~\ref{lem-twistopen}, up to
replacing $T'$ by a Zariski open subscheme, we may suppose that
$\zeta:u^*B \rightarrow \Kgnb{0,1}{X,1}$ is twisting.  Similarly, if
$(t_0,\zeta_0)$ is very twisting, we may suppose that $\zeta$ is very
twisting.  So we conclude that on the Zariski open subscheme of $T$
which is the image of $u:T'\rightarrow T$, the family $\xi:B
\rightarrow \Kgnb{0,0}{X,e}$ is twistable (resp. very twistable).
Since this holds for every point $t_0\in T$ where $\xi_0$ is
twistable, the lemma is proved.
\end{proof}

\begin{lem} \label{lem-twisttogether2}
  Suppose given two families
\begin{equation}
\xi_i = \lt( \lt( \pi_i:B_i \rightarrow T, \sigma_i:T \rightarrow B_i \rt),
h_i:B_i \rightarrow X \rt), i=1,2.
\end{equation}
such that for each $(\pi_i:B_i \rightarrow T, h_i: B_i \rightarrow X)$
is twistable, and such that $h_1\circ \sigma_1 = h_2 \circ \sigma_2$.
For each $t\in T$, assume that the locus of free lines in $X$ passing
through $h_1\circ\sigma_1(t) = h_2\circ \sigma_2(t)$ is irreducible.

\

Let us denote by
\begin{equation}
\xi =  \lt( \pi:B \rightarrow T, h:B \rightarrow X \rt)
\end{equation}
the family obtained by taking $B$ to be the connected sum of $B_1$ and
$B_2$ where the section $\sigma_1$ is identified with the section
$\sigma_2$.  Then $\xi$ is a twistable family.  Moreover, if at least
one of $\xi_1, \xi_2$ is very twistable, then $\xi$ is very twistable.
\end{lem}

\begin{proof}
  This follows essentially by lemma~\ref{lem-twisttogether}.  First of
  all, using proposition~\ref{prop-twistopen}, it suffices to prove
  the result when $T=\SP k$ for some algebraically closed field $k$.
  We suppose that we are in this case.

\

For each of $i=1,2$, let ${\mc M}_1$ denote the fiber product
constructed in the proof of proposition~\ref{prop-twistopen}, i.e.
${\mc M}_1$ parametrizes pairs $(\zeta_i,\theta_i)$ where $\zeta_i:B_i
\rightarrow \Kgnb{0,1}{X,1}$ is a twisting family (resp. very twisting
family) such that the induced map
\begin{equation}
\widetilde{\zeta}_i = \lt( (B_i,\sigma_i),g_i\circ \rho_i:B_i
\rightarrow X \rt)
\end{equation}
is stable, and where $\theta_i:\xi_i \rightarrow
\widetilde{\zeta}_i$ is an equivalence of objects.
Since each of $\xi_i$ is twistable, we see that each of ${\mc M}_i$ is
nonempty.  

\

By the proof of proposition~\ref{prop-twistopen}, each of ${\mc M}_i$
is smooth.  By the definition of twisting families, for each $i=1,2$
the morphism
\begin{equation}
e_i: {\mc M}_i \rightarrow \Kgnb{0,1}{X,1},\ \zeta_i\mapsto
\zeta_i(\sigma_i) 
\end{equation}
has image
contained in the unobstructed locus of $\text{ev}:
\Kgnb{0,1}{X,1}\rightarrow X$.  Let $P\subset \Kgnb{0,1}{X,1}$ be the
preimage under $\text{ev}$ of the point
$p=h_1(\sigma_1)=h_2(\sigma_2)$.  The image of $e_i$ is contained in
the smooth locus of $P$.  The claim is that $e_i:{\mc M}_i \rightarrow
P$ is smooth.  By \cite[proposition I.2.14.2]{K}, the obstruction space at a
point $\zeta_i$ is contained in the cohomology group $H^1\lt( B_i,
\zeta_i^* T_{\text{ev}}(-\sigma_i) \rt)$.  By the definition of a
twisting family, $\zeta_i^* T_{\text{ev}}$ is generated by global
sections.  Thus, by lemma~\ref{lem-gend}, the cohomology group above
is zero.  Since the obstruction space vanishes, we conclude that $e_i$
is smooth.  

\

Since both $e_1:{\mc M}_1 \rightarrow P$ and $e_2:{\mc M}_2
\rightarrow P$ are smooth, both have nonempty, open image.  And $P$ is
irreducible by assumption.  Therefore the image of $e_1$ and the image
of $e_2$ intersect.  If we choose a family $\zeta_1\in {\mc M}_1$ and
$\zeta_2\in {\mc M}_2$ such that $e_1(\zeta_1)=e_2(\zeta_2)$, then we
can glue $\zeta_1$ and $\zeta_2$ to obtain a morphism $\zeta:B
\rightarrow \Kgnb{0,1}{X,1}$ such that $\zeta|_{B_1}=\zeta_1$ and
$\zeta|_{B_2} = \zeta_2$.  By lemma~\ref{lem-twisttogether}, we
conclude that $\zeta$ is twisting.  Moreover, if at least one of
$\zeta_i, i=1,2$ is very twisting, then $\zeta$ is very twisting.  And
$\widetilde{\zeta} = \xi$.  This shows that $\xi$ is twistable, and it
is very twistable if at least one of $\xi_i, i=1,2$ is very twistable.
\end{proof}

\begin{hyp} \label{hyp-2}
  Let $U \subset \Kgnb{0,1}{X,1}$ denote the preimage of
  $U_{\text{t-able}} \subset \Kgnb{0,0}{X,1}$ under $\text{pr}$.  The
  evaluation morphism $\text{ev}:U \rightarrow X$ has Zariski dense
  image.
\end{hyp}

\begin{defn} \label{defn-pos}
  Suppose $\pi:B \rightarrow T$ is a family of prestable,
  geometrically connected curves of arithmetic genus $0$.  A
  $1$-morphism $\zeta:B \rightarrow \Kgnb{0,1}{X,e}$ is
  \emph{positive} (resp. \emph{very positive}) if:
\begin{enumerate}
\item The data $(\pi:B \rightarrow T, h:B \rightarrow X)$ is a family
  of stable maps to $X$, i.e. a $1$-morphism $\xi:T \rightarrow
  \Kgnb{0,0}{X,\epsilon}$ for some $\epsilon\geq 0$.
\item The image of $\xi:T \rightarrow \Kgnb{0,0}{X,\epsilon}$ is contained in
  the \emph{very unobstructed} locus of $\Kgnb{0,0}{X,\epsilon}$.
\item The image of $\text{pr}\circ \zeta: T\rightarrow \Kgnb{0,0}{X,e}$
  is contained in the \emph{very unobstructed} locus of
  $\Kgnb{0,0}{X,e}$.
\item The pullback bundle $(\text{pr}\circ \zeta)^*
  T_{\Kgnb{0,0}{X,e}}$ is $\pi$-relatively deformation ample.
\item The pullback line bundle $\sigma^*\OO_{\Sigma}(\sigma)$ is
  $\pi$-relatively generated by global sections
  (resp. $\pi$-relatively ample).  
\end{enumerate}
\end{defn}

\begin{rmk} \label{rmk-pos}
  Regarding the definition above:
\begin{enumerate}
  
\item This definition is very similar to
  definition~\ref{defn-twisting}.  It differs in that $e$ need not
  equal $1$ and that we only require $\text{pr}\circ \zeta$ to have
  image in the very unobstructed locus, instead of requiring $\zeta$
  to have image in the very unobstructed locus of $\text{ev}$.
  
\item Consider the case when $T = \SP \kappa$ for some field $\kappa$,
  and suppose that $B$ is smooth, i.e. $B \cong \PP^1_\kappa$.  If
  $\zeta:B \rightarrow \Kgnb{0,1}{X,e}$ is positive, then the morphism
  $\text{pr} \circ \zeta: B \rightarrow \Kgnb{0,0}{X,e}$ is \emph{very
    free} in the sense of Debarre \cite[p. 86]{De}.

\end{enumerate}
\end{rmk}

\begin{lem}\label{lem-posopen}
  Let $\pi:B \rightarrow T$ be a family of prestable, geometrically
  connected curves of arithmetic genus $0$ and let $\zeta:B
  \rightarrow \Kgnb{0,1}{X,1}$ be a $1$-morphism.  There is an open
  subscheme $U_{\text{pos}} \subset T$ (resp. $U_{\text{v-pos}}\subset
  T$) with the following property: for any morphism of schemes
  $f:T'\rightarrow T$, the pullback family $f^*\pi: f^*B \rightarrow
  T'$ and $f^*\zeta: f^*B \rightarrow \Kgnb{0,1}{X,1}$ is positive
  (resp. very positive) iff $f(T')\subset U_{\text{free}}$ (resp.
  $f(T')\subset U_{\text{v-free}}$).
\end{lem}

\begin{proof}
  The proof is almost identical to the proof of
  lemma~\ref{lem-twistopen}.
\end{proof}

\begin{lem} \label{lem-posdef}
  Suppose that $T=\SP k$ is a point and $\zeta:B\rightarrow
  \Kgnb{0,1}{X,e}$ is a positive $1$-morphism whose image is contained
  in the locus of very stable maps.
\begin{enumerate}

\item If $B$ is smooth, then $\zeta:B \rightarrow \Kgnb{0,1}{X,e}$ is
  \emph{free} in the sense of Koll\'ar ~\cite[definition II.3.11]{K}.
  If $\zeta$ is very positive, then $\zeta$ is \emph{very free}.

\item In any case, the $1$-morphism is \emph{unobstructed} in the
  sense of Koll\'ar \cite[definition I.2.6]{K}; in particular it is
  the specialization of a positive $1$-morphism $\zeta_\eta:B_\eta
  \rightarrow \Kgnb{0,1}{X,e}$ with $B_\eta$ geometrically connected
  and smooth, and whose image is contained in the locus of very stable
  maps.
\end{enumerate}
\end{lem}

\begin{proof}
  Suppose that the image of $\zeta$ lies in the locus of very stable
  maps.  By the relative version of lemma~\ref{lem-def3a}, we have
  that the image of $\zeta$ is in the smooth locus of
  $\Kgnb{0,1}{X,e}$.  And we have a short exact sequence:
\begin{equation}
\begin{CD}
0 @>>> \sigma^* \OO_{\Sigma}(\sigma) @>>> \zeta^* T_{\Kgnb{0,1}{X,e}}
@>>> (\text{pr}\circ \zeta)^* T_{\Kgnb{0,0}{X,e}} @>>> 0
\end{CD}
\end{equation}
By condition $(4)$ of definition ~\ref{defn-pos}, $(\text{pr}\circ
\zeta)^* T_{\Kgnb{0,0}{X,e}}$ is deformation ample.  And by condition
$(5)$ of definition ~\ref{defn-pos}, $\sigma^*\OO_{\Sigma}(\sigma)$ is
generated by global sections (resp. deformation ample).  Therefore
$\zeta^* T_{\Kgnb{0,1}{X,e}}$ is generated by global sections.  And if
$\zeta$ is very positive, then $\zeta^* T_{\Kgnb{0,1}{X,e}}$ is
deformation ample by $(2)$ of lemma ~\ref{lem-secDA}.  So if $B$ is
smooth, then $\zeta$ is free, and it is very free if $\zeta$ is very
positive.  This proves $(1)$.

\

Since $\zeta^* T_{\Kgnb{0,1}{X,e}}$ is generated by global sections,
by lemma ~\ref{lem-gend}, $H^1(B,\zeta^* T_{\Kgnb{0,1}{X,e}})$ is
zero.  Therefore $\zeta$ is unobstructed in the sense of Koll\'ar.
Now let $\pi: {\mc B} \rightarrow \SP R$ be a smoothing of $B$, i.e. a
flat family of proper, geometrically connected, prestable curves of
arithmetic genus $0$ over a DVR such that the special fiber is
isomorphic to $B$ and such that the general fiber $B_\eta$ is smooth.
By ~\cite[theorem I.2.10]{K}, the projection of the relative
Hom-scheme,
\begin{equation}
\textit{Hom}_{\SP R}({\mc B}, \SP R \times \Kgnb{0,1}{X,e})
\rightarrow \SP R,
\end{equation}
is smooth at $[\zeta]$.  Therefore, after making some finite, flat
base change $\SP R' \rightarrow \SP R$ we may suppose that $\zeta$ is
the specialization of a $1$-morphism $\zeta_R:{\mc B} \rightarrow
\Kgnb{0,1}{X,e}$.  By lemma ~\ref{lem-posopen}, we have that $\zeta_R$
is positive, in particular $\zeta_\eta:B_\eta \rightarrow
\Kgnb{0,1}{X,e}$ is positive.  Since the locus of very stable maps in
$\Kgnb{0,1}{X,e}$ is open, we conclude that the image of $\zeta_R$ is
contained in this locus.
\end{proof}

Now we come to the main notion of this section.

\begin{defn} \label{defn-inducts}
  Suppose $\pi:B \rightarrow T$ is a family of prestable,
  geometrically connected curves of arithmetic genus $0$.  An
  \emph{inducting pair} of degree $e$ is a pair
\begin{equation}
\lt( \zeta_1:B\rightarrow \Kgnb{0,1}{X,1}, \overline{\zeta}_e:B
\rightarrow \Kgnb{0,1}{X,e} \rt),
\end{equation}
such that:
\begin{enumerate}

\item $\zeta_1$ is very twisting,

\item $\overline{\zeta}_e$ is very positive and the image of
  $\overline{\zeta}_e$ is contained in the locus of very stable maps,
  and

\item the two morphisms $h_{\zeta_1}:B \rightarrow X$ and
  $h_{\overline{\zeta}_e}: B \rightarrow X$ are equal.
\end{enumerate}
\end{defn}

\begin{lem} \label{lem-inducting}
  Let $\pi:B \rightarrow T$ be a family of prestable, geometrically
  connected curves of arithmetic genus $0$, and let
\begin{equation}
\lt( \zeta_1:B
\rightarrow \Kgnb{0,1}{X,1}, \overline{\zeta}_e:B \rightarrow
\Kgnb{0,1}{X,e} \rt),
\end{equation}
be a pair of $1$-morphisms such that $h_{\zeta_1} =
h_{\overline{\zeta}_e}$.  Then there is an open subscheme
$U_{\text{induct}}\subset T$ with the following property: for any
morphism of schemes $f:T'\rightarrow T$, the pullback of
$(\zeta_1,\overline{\zeta}_e)$ is inducting iff $f(T')\subset U$.
\end{lem}

\begin{proof}
  We just define $U_{\text{induct}}$ to be the intersection of the
  open subset $U_{\text{vtwist}}\subset T$ as in
  lemma~\ref{lem-twistopen} for $\zeta_1$ and the open subset
  $U_{\text{v-pos}}\subset T$ as in lemma~\ref{lem-posopen} for
  $\overline{\zeta}_e$.
\end{proof}

We finally come to our last definition.

\begin{defn} \label{defn-inductable}
  Suppose $(\pi:B\rightarrow T, \overline{\zeta}_e:B \rightarrow
  \Kgnb{0,1}{X,e})$ is a very positive family whose image is contained
  in the locus of very stable maps.  We say $\overline{\zeta}_e$ is
  \emph{inductable} if there is a surjective \'etale morphism $u:T'
  \rightarrow T$ and a morphism $\zeta_1:u^*B \rightarrow
  \Kgnb{0,1}{X,1}$ with $h_{\zeta_1} = u^* h_{\overline{\zeta}_e}$
  such that $(\zeta_1,\overline{\zeta}_e)$ is an inducting pair.
\end{defn}

\begin{lem} \label{lem-inductable}
  Let $(\pi:B \rightarrow T, \overline{\zeta}_e:B \rightarrow
  \Kgnb{0,1}{X,e})$ be a very free family.  There is an open subscheme
  $U_{i-able} \subset T$ such that for each morphism of schemes
  $f:T'\rightarrow T$, the pullback $(f^*\pi: f^*B \rightarrow T',
  f^*\overline{\zeta}_e: f^* B \rightarrow \Kgnb{0,1}{X,e})$ is
  inductable iff $f(T')\subset U_{i-able}$.
\end{lem}

\begin{proof}
  We simply apply proposition~\ref{prop-twistopen} to
\begin{equation}
\xi := (\pi:B
\rightarrow T, h_{\overline{\zeta}_e}: B \rightarrow X).
\end{equation}
\end{proof}

\section{The induction argument} \label{sec-induct}

In this section we will show that given an inductable $1$-morphism
$\overline{\zeta}_e: B \rightarrow \Kgnb{0,1}{X,e}$, this gives rise
to an inductable $1$-morphism $\overline{\zeta}_{e+1}: B \rightarrow
\Kgnb{0,1}{X,e+1}$.  The basic idea is, given an inducting pair
$(\zeta_1,\overline{\zeta}_{e+1})$ to form the family of connected
sums.  This isn't quite an inductable $1$-morphism, but after
deforming and then performing a simple operation which we call a
\emph{modification}, we do obtain an inductable $1$-morphism
$\overline{\zeta}_{e+1}$.

\begin{notat} \label{not-divs}
  We will follow \cite{QDiv} in our notation of the tautological
  divisors on $\Kgnb{0,1}{\PP^N,e}$. Specifically, in
  $A^1(\Kgnb{0,1}{\PP^N,e})_\QQ$ we denote by $\Delta_{(e_1,e_2)}$ the
  $\QQ$-divisor whose general point parametrizes a reducible embedded
  curve with one irreducible component of degree $e_1$, one
  irreducible component of degree $e_2$ and where the marked point is
  on the first irreducible component.  We denote by ${\mc L}$ the
  divisor class $\text{ev}^*\OO_{\PP^N}(1)$.  And we denote by ${\mc
    H}$ the divisor which parametrizes stable maps whose image in
  $\PP^r$ intersects a given codimension $2$ linear space.  Given a
  closed subscheme $X\subset \PP^N$, we also denote by
  $\Delta_{(e_1,e_2)}$, ${\mc L}$ and ${\mc H}$ the pullbacks of the
  divisors given above by the induced $1$-morphism $\Kgnb{0,1}{X,e}
  \rightarrow \Kgnb{0,1}{\PP^N,e}$.
\end{notat}

Before proceeding to the main result of this section, we describe an
operation which we will perform repeatedly in the proof.  Suppose that
$\zeta:B \rightarrow \Kgnb{g,r}{X,e}$ is a family of stable maps as in
notation ~\ref{not-1mor}.  Suppose $b\in B$ is a point whose image
$\zeta(b)$ is a stable map
\begin{equation}
((\Sigma_b, p_1, \dots, p_r), g_b:\Sigma_b
\rightarrow X).
\end{equation}  
Suppose that $L\subset \Sigma_b$ is an irreducible component which is
not contracted by $g_b$ and $p_i\in L$ is one of the marked points.
For simplicity assume that $L$ contains no nodes, in particular this
is the case when $L$ has genus $0$.  Let $M\subset \Sigma_b$ denote
all the irreducible components of $\Sigma_b$ other than $L$.  Let
$R=(r_1,\dots, r_l)$ denote the set of intersection points of $L$ and
$M$ with some ordering.  Let $D = (p_{j_1},\dots, p_{j_m})$ denote the
set of marked points which lie on $L$ other than $p_i$ and let
$E=(p_{k_1},\dots, p_{k_n})$ denote the set of marked points which lie
on $M$.

\

Form the product surface $L\times L$ with diagonal $\Delta:L
\rightarrow L\times L$ and let $u:\Lambda \rightarrow L$ denote the
blowing up of $L\times L$ along the set of points $\Delta(R\cup D)$.
For each point $p \in R\cup D$, let $F_{p} \subset \Lambda$ be the
proper transform of $L \times \{ p \}$.  Let $F_{p_i}$ denote the
proper transform of the diagonal $\Delta(L) \subset L \times L$.
Consider $\text{pr}_1\circ u: \Lambda \rightarrow L$ as a family of
prestable curves parametrized by $L$.  Then the data
\begin{equation}
\widetilde{\zeta}_{\Lambda} = 
(\text{pr}_1 \circ u: \Lambda
\rightarrow L, (\sigma_{p_i}, \sigma_p | p\in R\cup D) ),
\end{equation}
where $\sigma_p:L \rightarrow F_p$ is the unique isomorphism such that
$\text{pr}_1\circ \sigma_p$ is the identity, gives a family of
prestable marked curves parametrizing by $L$, it is essentially the
constant family $L\times (L,\{p_i\}\cup R \cup D) \rightarrow L$
except that we are allowing $p_i$ to vary among all points in $L$ and
then blowing up to obtain a stable family.

\

Next we form the constant family of prestable marked curves
parametrized by $L$:
\begin{equation}
\widetilde{\zeta}_{L\times M} =
(\text{pr}_1: L \times M 
\rightarrow L, (s_p | p \in R \cup E)
\end{equation}
where $s_p: L \rightarrow L\times M$ is simply $s_p(t)=(t,p)$.  We can
glue $\widetilde{\zeta}_{\Lambda}$ and $\widetilde{\zeta}_{L\times M}$
as follows.  For each $p\in R$, we identify the section $\sigma_p$ of
$\widetilde{\zeta}_{\Lambda}$ with the section $s_p$ of
$\widetilde{\zeta}_{L\times M}$.  Here the identification is the
unique one compatible with projection to $L$.  Let us denote the new
family of prestable marked curves by
\begin{equation}
(\rho:\Pi \rightarrow L, (\phi_j: L
\rightarrow \Pi | j=1,\dots, r)
\end{equation}
where $\Pi$ is the surface obtained by gluing $\Lambda$ and $L\times
M$ as above, and where
\begin{equation}
\phi_j = \lt\{ \begin{array}{ll}
                \sigma_{p_j} & ,p_j\in L \\
                s_{p_j}      & ,p_j\in M
               \end{array} \rt.
\end{equation}
Notice that there is a unique morphism $\text{pr}_2:\Pi \rightarrow
\Sigma_b$ whose restriction to $\Lambda$ is $\text{pr}_2\circ u:
\Lambda \rightarrow L \subset \Sigma_b$ and whose restriction to $L
\times M$ is $\text{pr}_2: L \times M \rightarrow M\subset \Sigma_b$.
We form a family of stable maps parametrized by $L$,
\begin{equation}
\widetilde{\zeta}_\Pi = \lt( \lt( \rho: \Pi \rightarrow L, \phi_1,\dots,
\phi_r \rt), g_b 
\circ \text{pr}_2: \Pi \rightarrow X \rt).
\end{equation}
Notice that the family $\widetilde{\zeta}_\Pi$ is stable, and if we
remove the section $\phi_i$ and stabilize, we just get the constant
family parametrized by $L$ whose image is the stabilization of
$\zeta(b)$ upon removing $p_i$.  Also, the image
$\widetilde{\zeta}_\Pi(p_i)$ is precisely $\zeta(b)$.  Let $\tilde{B}$
be the connected sum of $B$ and $L$ where $b\in B$ is identified with
$p_i\in L$.  Since $\widetilde{\zeta}_\Pi(p_i) = \zeta(b)$, we may
form a $1$-morphism $\widetilde{\zeta}: \widetilde{B} \rightarrow
\Kgnb{g,r}{X,e}$ such that $\widetilde{\zeta}$ restricted to $B$ is
$\zeta$, and $\widetilde{\zeta}$ restricted to $L$ is
$\widetilde{\zeta}_\Pi$.

\begin{notat} \label{not-modif}
  Given a $1$-morphism $\zeta:B \rightarrow \Kgnb{g,r}{X,e}$, a point
  $b\in B$, an irreducible component $L\subset \Sigma_b$, and a marked
  point $p_i\in L$ as above, we call the $1$-morphism
  $\widetilde{\zeta}:\widetilde{B} \rightarrow \Kgnb{g,r}{X,e}$
  constructed in the last paragraph the \emph{modification} of $\zeta$
  determined by $b\in B$, by $L$ and by $p_i$.
\end{notat}

\begin{lem} \label{lem-modif}
  Suppose given a $1$-morphism $\zeta:B \rightarrow \Kgnb{0,r}{X,e}$,
  a point $b\in B$, an irreducible component $L\subset \Sigma_b$ which
  is not contracted, and a marked point $p_i\in L$ such that when we
  remove $p_i$, the resulting stable map is a smooth point of
  $\Kgnb{0,r-1}{X,e}$ (this condition is equivalent to the condition
  that the image of $\widetilde{\zeta}_\Pi$ is contained in the smooth
  locus of $\Kgnb{0,r}{X,e}$).  Then we have the vanishing
\begin{equation}
H^1 \lt( L,\widetilde{\zeta}_\Pi^*T_{\Kgnb{0,r}{X,e}} \rt) = 0.
\end{equation}
In particular, if $\zeta:B \rightarrow \Kgnb{0,r}{X,e}$ is a free
morphism of a rational curve into the smooth locus, then there are
deformations of $\widetilde{\zeta}:\widetilde{B} \rightarrow
\Kgnb{0,r}{X,e}$ which smooth the node of $\widetilde{B}$.
\end{lem}

\begin{proof}
  This is an application of the deformation theory of section
  ~\ref{sec-def}.  Let $\tau$ denote the dual graph of $\zeta(b)$ and
  let $\psi:\tau \rightarrow \tau'$ be the combinatorial morphism of
  graphs which removes the tail associated to $p_i$.  The morphism
  $\Kbm(X,\psi):\Kbm(X,\tau) \rightarrow \Kbm(X,\tau')$ is smooth
  along the image of $\widetilde{\zeta}_\Pi$ and the pullback of the
  vertical tangent bundle is simply $T_L$.  So
  $\widetilde{\zeta}_\Pi^* T_{\Kbm(X,\tau)}$ is generated by global
  sections.  The pullback of the normal sheaf ${\mc N}$ of
  $\Kbm(X,\tau) \rightarrow \Kgnb{0,r}{X,e}$ is the direct sum over
  all $r_j\in R$ of $N_{F_j/\Lambda} \otimes_{\CC} T_{r_j} M$.  As the
  normal bundle of $N_{F_j/\Lambda}$ is just $\OO_L(-1)$, we conclude
  that $H^1(L,{\mc N}) = 0$.  We have a short exact sequence:
\begin{equation}
\begin{CD}
0 @>>> \widetilde{\zeta}_\Pi^*T_{\Kbm(X,\tau)} @>>> 
\widetilde{\zeta}_\Pi^*T_{\Kgnb{0,r}{X,e}}  @>>> {\mc N} @>>> 0
\end{CD}
\end{equation}
In the corresponding long exact sequence of cohomology, $H^1$ of the
first and third terms vanishes.  Therefore we have the vanishing
result.

\

Suppose that $\zeta:B \rightarrow \Kgnb{0,r}{X,e}$ is a free
$1$-morphism of a rational curve into the smooth locus of
$\Kgnb{0,r}{X,e}$.  Then $\zeta^*T_{\Kgnb{0,r}{X,e}}$ is generated by
global sections, so $H^1(B,\zeta^* T_{\Kgnb{0,r}{X,e}}(-b))$ is zero
by lemma~\ref{lem-gend} We have a short exact sequence:
\begin{equation}
\begin{CD}
0 @>>> \zeta^* T_{\Kgnb{0,r}{X,e}}(-b) @>>> \widetilde{\zeta}^*
T_{\Kgnb{0,r}{X,e}} @>>> \widetilde{\zeta}_\Pi^* T_{\Kgnb{0,r}{X,e}}
@>>> 0
\end{CD}
\end{equation}
In the corresponding long exact sequence of cohomology, $H^1$ of the
first and third terms vanishes.  Therefore
$H^1(\widetilde{B},\widetilde{\zeta}^* T_{\Kgnb{0,r}{X,e}} )$
vanishes.  This cohomology group is the obstruction to smoothing the
node, therefore there are deformations of
$\widetilde{\zeta}:\widetilde{B}\rightarrow \Kgnb{0,r}{X,e}$ which
smooth the node of $\widetilde{B}$.
\end{proof}

\begin{rmk} \label{rmk-modif}
  In case the line bundle $\sigma^*\OO_{\Sigma}(\sigma)$ is generated
  by global sections and $L$ contains only one node of $\Sigma_b$, we
  have a simpler proof of the deformation result.  We have a short
  exact sequence of sheaves on $\Sigma$:
\begin{equation}
\begin{CD}
0 @>>> \OO_{\Sigma}(\sigma) @>>> \OO_{\Sigma}(\sigma + L) @>>>
N_{L/\Sigma}(p_i) @>>> 0 
\end{CD}
\end{equation}
Since $L$ contains only one node, $N_{L/\Sigma}\cong \OO_L(-1)$, so
the last term is isomorphic to $\OO_L$.  So $\OO_{\Sigma}(\sigma+L)$
is generated by global sections.  A small deformation $\sigma'$ of
$\sigma + L$ in the linear series $|\sigma + L|$ will be a section of
$\pi:\Sigma \rightarrow B$, and the stabilization of the $1$-morphism
$B \rightarrow \Kgnb{0,r}{X,e}$ which removes the section $\sigma$
from $\zeta$ and replaces it by $\sigma'$ will be a small deformation
of $\widetilde{\zeta}$ which smooths the node of $\widetilde{B}$.
\end{rmk}

\

Now we come to the main theorem of this section, which we use for the
induction step in the proof of theorem~\ref{thm-thm1}.

\begin{thm} \label{thm-induction}
  Suppose that $X$ satisfies hypothesis~\ref{hyp-1},
  hypothesis~\ref{hyp-1.5}, hypothesis~\ref{hyp-1.75}, and
  hypothesis~\ref{hyp-2}.  For each integer $e\geq 1$, if there exists
  an inductable map $\overline{\zeta}_e: B_e \rightarrow
  \Kgnb{0,1}{X,e}$, then there exists an inductable map
  $\overline{\xi}_{e+1}: B_{e+1} \rightarrow \Kgnb{0,1}{X,e+1}$.
  
  More precisely, suppose that $(\zeta_1,\overline{\zeta}_e)$ is an
  inducting pair.  Let us denote:
\begin{eqnarray}
s=\text{deg}(\zeta_1)^*(2{\mc
L} - {\mc H}), \\
\overline{s}=\text{deg}(\overline{\zeta}_e)^*(2{\mc
L} - {\mc H})
\end{eqnarray}  
Then for each $k=1,\dots,\overline{s}$, there is an inducting pair
$(\xi^k_1,\overline{\xi}^k_{e+1})$ satisfying the following.
\begin{enumerate}
  
\item We have
\begin{equation}
\text{deg}\lt((\xi^k_1)^*{\mc H}\rt) = \text{deg}\lt(\zeta_1^*{\mc H}
\rt) + 2k. 
\end{equation}

\item We have
\begin{equation}
\text{deg}\lt( (\xi^k_1)^*{\mc L} \rt) = \text{deg}\lt(
(\overline{\xi}^k_{e+1})^*{\mc L} \rt) = \text{deg}\lt( \zeta_1^*{\mc L}
\rt) + k.
\end{equation}

\item We have
\begin{equation}
\text{deg}\lt( (\overline{\xi}^k_{e+1})^*{\mc H} \rt) = \text{deg} \lt(
\zeta_1^*{\mc H} \rt) + \text{deg} \lt( \overline{\zeta}_e^*{\mc H}
\rt).
\end{equation}

\item For each $e_1 + e_2 = e$ with both $e_1, e_2 \geq 2$, we have
\begin{equation}
\text{deg}\lt( (\overline{\xi}^k_{e+1})^* \Delta_{(e_1+1,e_2)} \rt) = 
\text{deg}\lt( \overline{\zeta}_{e}^* \Delta_{(e_1,e_2)} \rt).
\end{equation}

\item If $e > 1$, we have
\begin{equation}
\text{deg}\lt( (\overline{\xi}^k_{e+1})^* \Delta_{(e,1)} \rt)
=\text{deg}\lt( \overline{\zeta}_{e}^* \Delta_{(e-1,1)} \rt)+ s+ k,
\end{equation}
and we have
\begin{equation}
\text{deg}\lt( (\overline{\xi}^k_{e+1})^* \Delta_{(1,e)} \rt) = \overline{s}-k.
\end{equation}

\item If $e=1$, we have
\begin{equation}
\text{deg}\lt( (\overline{\xi}^k_{2})^* \Delta_{(1,1)} \rt) = s+\overline{s}.
\end{equation}
\end{enumerate}

\end{thm}

\begin{proof}
  By lemma ~\ref{lem-posdef} and lemma ~\ref{lem-inductable}, we may
  suppose that $B_e$ is smooth.  Moreover, in this case
  $\overline{\zeta}_e:B_e \rightarrow \Kgnb{0,1}{X,e}$ is free (in
  fact very free).  Therefore, we may suppose that
  $\overline{\zeta}_e(B_e)$ is in general position: for any finite
  collection of codimension $2$ subvarieties $(Z_\alpha |
  \alpha=1,\dots,M)$ and any finite collection of divisors $(D_\beta |
  \beta=1,\dots, N)$, we may suppose that $\overline{\zeta}_e(B_e)$ is
  disjoint from each $Z_\alpha$ and has $0$-dimensional intersection
  with each $D_\beta$.

\
 
Let us denote the family of stable maps $\overline{\zeta}_e$ by:
\begin{equation}
\lt( \overline{p}:\overline{\Sigma} \rightarrow B,
\overline{\sigma}:B \rightarrow \overline{\Sigma}, \overline{g}:
\overline{\Sigma} \rightarrow X \rt).
\end{equation}
And let us denote the family of stable maps $\zeta_1$ by
\begin{equation}
\lt( p:\Sigma \rightarrow B, \sigma:B \rightarrow \Sigma, g: \Sigma
\rightarrow X \rt).
\end{equation}
The basic idea is to form the connected sum of the surfaces $\Sigma$
and $\overline{\Sigma}$ glued along the sections $\sigma$ and
$\overline{\sigma}$.  The actual family $\overline{\zeta}_{e+1}$ is a
bit more complicated.

\

Define $\pi':\Sigma'\rightarrow B$ to be the family of curves obtained
by taking the connected sum of $\Sigma$ and $\overline{\Sigma}$ glued
along the sections $\sigma$ and $\overline{\sigma}$.  Here $\pi'$ is
the unique morphism such that $\pi'|_{\Sigma} = \pi$ and
$\pi'|_{\overline{\Sigma}} = \overline{\pi}$.  Define $g':\Sigma'
\rightarrow X$ to be the unique morphism such that $g'|_{\Sigma} = g$
and $g'|_{\overline{\Sigma}} = \overline{g}$.  Then $\zeta' =
(\pi':\Sigma'\rightarrow B, g':\Sigma' \rightarrow X)$ is a family of
stable maps in the boundary divisor $\Delta_{e,1}$ of
$\Kgnb{0,0}{X,e+1}$.  Moreover, $\zeta'$ clearly factors through the
Behrend-Manin stack $\Kbm(X,\tau) \rightarrow X$ where $\tau$ is the
genus $0$ stable $A$-graph with two vertices $v_1, v_2$ with
$\beta(v_1)=1$ and $\beta(v_2)=e$.  By lemma~\ref{lem-glue1}, we have
a short exact sequence:
\begin{equation}
\begin{CD}
0 @>>> \zeta_1^* T_{\text{ev}} @>>> (\zeta')^* T_{\Kbm(X,\tau)} @>>>
\overline{\zeta}_e^* T_{\Kgnb{0,1}{X,e}} @>>> 0
\end{CD}
\end{equation}
Since $\zeta_1$ is very twisting, by definition
~\ref{defn-verytwisting} we have that $\zeta_1^* T_{\text{ev}}$ is
ample.  Since $\overline{\zeta}_e$ is very positive, by lemma
~\ref{lem-posdef}, we have that $\overline{\zeta}_e^*
T_{\Kgnb{0,1}{X,e}}$ is ample.  Therefore $(\zeta')^*
T_{\Kbm(X,\tau)}$ is ample.  By lemma~\ref{lem-glue2}, we have a short
exact sequence:
\begin{equation}
0 \rightarrow (\zeta')^* T_{\Kbm(X,\tau)} \rightarrow (\zeta')^*
T_{\Kgnb{0,0}{X,e+1}} \rightarrow \sigma^* \OO_{\Sigma}(\sigma) \otimes
\overline{\sigma}^* \OO_{\overline{\Sigma}}(\overline{\sigma}) \rightarrow 0
\end{equation}
By definition ~\ref{defn-verytwisting} and definition ~\ref{defn-pos},
we have that both $\sigma^*\OO_{\Sigma}(\sigma)$ and
$\overline{\sigma}^* \OO_{\overline{\Sigma}}(\overline{\sigma})$ are
ample.  Therefore their tensor product is ample, and we conclude that
$(\zeta')^* T_{\Kgnb{0,0}{X,e+1}}$ is ample.

\

Denote by $\overline{s}$ the self-intersection of
$\overline{\sigma}\subset \overline{\Sigma}$, i.e. the degree of the
invertible sheaf
$\overline{\sigma}^*\OO_{\overline{\Sigma}}(\overline{\sigma})$.
Notice that we have
\begin{equation}
\overline{s} = \text{deg}\lt( 2\overline{\zeta}_e^*{\mc L} -
\overline{\zeta}_e^*{\mc H} \rt).
\end{equation}
Let $\sigma':B \rightarrow \overline{\Sigma}$ be a general member of
the linear series of $|\overline{\sigma}|$.  Since
$\overline{\sigma}^*\OO_{\overline{\Sigma}}(\overline{\sigma})$ is
generated by global sections, we can find such a $\sigma'$ which has
transverse intersections with $\overline{\sigma}$ at points
$p_1,\dots,p_{\overline{s}} \in \overline{\Sigma}$.  Define
$b:\widetilde{\overline{\Sigma}}\rightarrow \overline{\Sigma}$ to be
the blowing up of $\overline{\Sigma}$ at the points
$p_1,\dots,p_{\overline{s}}$.  Let
$\widetilde{\overline{\pi}}:\widetilde{\overline{\Sigma}} \rightarrow
\overline{\Sigma}$ denote the projection $\overline{\pi}\circ b$.  Let
$\widetilde{\overline{g}}$ denote $\overline{g}\circ b$.  Let
$\widetilde{\overline{\sigma}}:B \rightarrow
\widetilde{\overline{\Sigma}}$ and $\widetilde{\sigma}:B \rightarrow
\widetilde{\overline{\Sigma}}$ denote the proper transforms of
$\overline{\sigma}$ and $\sigma'$ respectively.  Notice that
$\widetilde{\overline{\sigma}}$ and $\widetilde{\sigma}$ are disjoint
sections.  So the data
\begin{equation}
\lt( \lt(\widetilde{\overline{p}}:\widetilde{\overline{\Sigma}}
\rightarrow B, \widetilde{\overline{\sigma}}, \widetilde{\sigma} \rt),
\widetilde{\overline{f}}: 
\widetilde{\overline{\Sigma}} \rightarrow X \rt)
\end{equation}
is a family of stable pointed maps, i.e. a $1$-morphism
$\widetilde{\overline{\zeta}}_e: B \rightarrow \Kgnb{0,2}{X,e}$. By
the deformation theory in subsection~\ref{subsec-unstable}, we have a
short exact sequence:
\begin{equation}
\begin{CD}
0 @>>> (\widetilde{\sigma})^*
\OO_{\widetilde{\overline{\Sigma}}}(\widetilde{\sigma}) @>>> \lt( 
\widetilde{\overline{\zeta}}_e \rt)^* T_{\Kgnb{0,2}{X,e}} @>>> \lt(
\overline{\zeta} \rt)_e^* T_{\Kgnb{0,1}{X,e}} @>>> 0
\end{CD}
\end{equation}
Of course we have that $(\widetilde{\sigma})^*
\OO_{\widetilde{\overline{\Sigma}}}(\widetilde{\sigma})$ is the
trivial invertible sheaf $\OO_B$.  In particular, we have that $\lt(
\widetilde{\overline{\zeta}}_e \rt)^* T_{\Kgnb{0,2}{X,e}}$ is
generated by global sections.

\

Define $\widetilde{\pi}:\widetilde{\Sigma} \rightarrow B$ to be the
family of curves obtained by taking the connected sum of $\Sigma$ and
$\widetilde{\overline{\Sigma}}$ glued along the sections $\sigma$ and
$\tau$ respectively.  Here $\widetilde{\pi}$ is the unique morphism
such that $\widetilde{\pi}|_{\widetilde{\overline{\Sigma}}} =
\widetilde{\overline{\pi}}$ and $\widetilde{\pi}|_\Sigma = \pi$.
Define $\widetilde{g}:\widetilde{\Sigma}\rightarrow X$ to be the
unique morphism such that
$\widetilde{g}|_{\widetilde{\overline{\Sigma}}} =
\widetilde{\overline{g}}$ and such that $\widetilde{g}|_{\Sigma} = g$.
Define $\widetilde{\sigma}:B \rightarrow \widetilde{\Sigma}$ to be the
section from the last paragraph.  This gives a family of stable maps:
\begin{equation}
\widetilde{\zeta} = \lt( \lt( \widetilde{\pi}:\widetilde{\Sigma}
\rightarrow B, \widetilde{\sigma}:B \rightarrow \widetilde{\Sigma}
\rt), \widetilde{g}: \widetilde{\Sigma} \rightarrow X \rt).
\end{equation}
Define $\tau$ to be the unique genus $0$ stable $A$-graph with two
vertices $v_1,v_2$, with one edge joining $v_1$ and $v_2$, with one
flag attached to $v_1$, with $\beta(v_1)=e$, and with $\beta(v_2)=1$.
In other words, $\tau$ is the dual graph of a stable map with one
marked point, with reducible domain consisting of two irreducible
components, where the component with the marked point has degree $e$
and where the other component has degree $1$.  Then
$\widetilde{\zeta}:B \rightarrow \Kgnb{0,1}{X,e+1}$ is a $1$-morphism
which factors through the Behrend-Manin stack $\Kbm(X,\tau)
\rightarrow \Kgnb{0,1}{X,e+1}$.  By lemma~\ref{lem-glue1}, we have a
short exact sequence:
\begin{equation}
\begin{CD}
0 @>>> \zeta_1^* T_{\text{ev}} @>>> \lt( \widetilde{\zeta} \rt)^*
T_{\Kbm(X,\tau)} @>>>
\lt( \widetilde{\overline{\zeta}}_e \rt)^* T_{\Kgnb{0,2}{X,e}} @>>> 0
\end{CD}
\end{equation}
As the first and third term in this exact sequence are generated by
global sections, also $\lt( \widetilde{\zeta} \rt)^*
T_{\Kbm(X,\gamma)}$ is generated by global sections.  Finally, by
lemma~\ref{lem-glue2} we have a short exact sequence:
\begin{equation}
0 \rightarrow \lt( \widetilde{\zeta} \rt)^* T_{\Kbm(X,\gamma)} \rightarrow \lt(
\widetilde{\zeta} \rt)^* T_{\Kgnb{0,1}{X,e+1}} \rightarrow \sigma^*
\OO_{\Sigma}(\sigma) \otimes_{\OO_B} \lt(
\widetilde{\overline{\sigma}} \rt)^*
\OO_{\widetilde{\overline{\Sigma}}}(\widetilde{\overline{\sigma}})
\rightarrow 0
\end{equation}
Of course $\lt( \widetilde{\overline{\sigma}} \rt)^*
\OO_{\widetilde{\overline{\Sigma}}}(\widetilde{\overline{\sigma}})$ is
the trivial invertible sheaf.  And $\sigma^* \OO_{\Sigma}(\sigma)$ is
just $\OO_B(s)$.  Thus the last term in the sequence is an ample
invertible sheaf.  In particular, $\lt( \widetilde{\zeta} \rt)^*
T_{\Kgnb{0,1}{X,e+1}}$ is generated by global sections.  So the
$1$-morphism $\widetilde{\zeta}: B \rightarrow \Kgnb{0,1}{X,e+1}$ is
free.  Also, the pullback by $\widetilde{\zeta}$ of the normal sheaf
of the unramified $1$-morphism, $\Kbm(X,\tau) \rightarrow
\Kgnb{0,1}{X,e+1}$ is the last term in the sequence, and so has
positive degree.

\
  
One final remark, when $e>1$, the image of $\widetilde{\zeta}$
intersects the divisor $\Delta_{1,e}$ transversely precisely at the
images of the points $p_1,\dots, p_{\overline{s}} \in B$, in
particular the degree of the $\QQ$-Cartier divisor class
$\widetilde{\sigma}^*\OO_{\Kgnb{0,1}{X,e+1}}(\Delta_{1,e})$ is
positive.  In the special case $e=1$, we have that
$\Delta_{1,e}=\Delta_{e,1} =\Delta_{1,1}$.  In this case
$\Kbm(X,\tau)$ is the normalization of $\Delta_{1,1}$ (at least in a
neighborhood of the image of $\widetilde{\sigma}$).  So the degree of
$\widetilde{\sigma}^*\OO_{\Kgnb{-,1}{X,e+1}}(\Delta_{1,1})$ is the sum
of the degree of the pullback of the normal sheaf of
$\Kbm(X,\tau)\rightarrow \Kgnb{0,1}{X,e+1}$ and the divisor $p_1 +
\dots + p_{\overline{s}}$, i.e. $s+ \overline{s}$.  So in this case,
the degree is again positive.

\

Since $\widetilde{\zeta}:B \rightarrow \Kgnb{0,1}{X,e+1}$ is free we
can find deformations of $\widetilde{\zeta}$ which are in general
position.  Now by hypothesis ~\ref{hyp-1}, the locus parametrizing
stable maps with at least three irreducible components in their domain
has codimension $2$.  By hypothesis~\ref{hyp-1.75} the locus
parametrizing stable maps with automorphisms has codimension at least
$2$.  Therefore we can find a deformation $\overline{\xi}_{e+1}^0:B
\rightarrow \Kgnb{0,1}{X,e+1}$ of $\widetilde{\zeta}$ such that the
image of $\overline{\xi}_{e+1}^0:B \rightarrow \Kgnb{0,1}{X,e+1}$
misses the locus parametrizing stable maps with at least three
irreducible components in their domain and misses the locus
parametrizing stable maps with automorphisms.  As well, we can assume
that the pullback of $T_{\Kgnb{0,0}{X,e+1}}$ by $\text{pr}\circ
\overline{\xi}_{e+1}^0: B \rightarrow \Kgnb{0,0}{X,e+1}$ is an ample
vector bundle since the pullback of $T_{\Kgnb{0,0}{X,e+1}}$ by
$\text{pr}\circ \widetilde{\zeta}$, i.e. by $\zeta':B \rightarrow
\Kgnb{0,0}{X,e+1}$, is an ample bundle.  Let us denote
$\overline{\xi}_{e+1}^0$ by
\begin{equation}
\overline{\xi}_{e+1}^0 = \lt( \lt( \varpi^0:\Xi^0 \rightarrow B, \lambda^0:B
\rightarrow \Xi^0 \rt), g^0:\Xi^0 \rightarrow B \rt).
\end{equation}
Define $h^0:B \rightarrow X$ to be $g^0\circ \lambda^0$.  By
assumption, $\overline{h}:B \rightarrow X$ is very twistable, and
$h^0:B \rightarrow X$ is a small deformation of $\overline{h}$.  So by
proposition ~\ref{prop-twistopen}, $h^0$ is very twistable.

\

Since
$\widetilde{\sigma}^*\OO_{\widetilde{\Sigma}}(\widetilde{\sigma})$ is
trivial, we may also assume that $\lt( \lambda^0 \rt)^*
\OO_{\Xi^0}(\lambda^0)$ is trivial.  Finally, we may assume that all
intersections of the image of $\overline{\xi}_{e+1}^0$ and the divisor
$\Delta_{1,e}$ are transverse and occur at general points of
$\Delta_{1,e}$.  In particular, if $b\in B$ is such a point, we may
assume that the corresponding stable map $\overline{\xi}_{e+1}^0(b)$
is of the form $((\Xi^0_b,\lambda^0_b),h^0_b:\Xi^0_b \rightarrow X)$
where $\Xi^0_b$ is a reducible curve $L \cup M$, with $L\cap M$
consisting of one node of $\Xi^0_b$ which is a general point of $X$,
with $\lambda^0_q \in L$, and such that $L \rightarrow X$ is an
embedding of a line which is twistable.

\

Now we define $\widetilde{\overline{\xi}}^0_{e+1}$ to be the
modification of $\overline{\xi}^0_{e+1}$ associated to the point $b\in
B$, to $L\subset \Xi^0_b$ and $\lambda^0_q \in L$ (c.f.
notation~\ref{not-modif}).  We saw above that $h^0:B \rightarrow X$ is
very twistable.  And the evaluation map $\widetilde{h}_L:L \rightarrow
X$ associated to the restriction $\widetilde{\overline{\xi}}^0_{e+1}:
L \rightarrow \Kgnb{0,1}{X,e+1}$ is just the embedding $L\subset X$,
which is twistable.  By hypothesis~\ref{hyp-1.5} and the assumption
that $L\cap M$ is a general point of $X$, we see that the hypotheses
of lemma~\ref{lem-twisttogether2} are satisfied.  So by
lemma~\ref{lem-twisttogether2} we have that the evaluation map
$\widetilde{h}:\widetilde{B}\rightarrow X$ associated to
$\widetilde{\overline{\xi}}^0_{e+1}$ is also very twistable.  Also,
notice that the image of $\widetilde{\overline{\xi}}^0_{e+1}$ is
contained in the locus of very stable maps.

\

The $1$-morphism $\widetilde{\overline{\xi}}^0_{e+1}$ satifies the
criterion of lemma~\ref{lem-modif}, in fact the criterion of
remark~\ref{rmk-modif}, thus we can smooth the node of
$\widetilde{B}$.  Let $\overline{\xi}^1_{e+1}:B \rightarrow $ denote a
small deformation which smooths the node of $\widetilde{B}$.  We will
denote this by
\begin{equation}
\overline{\xi}^1_{e+1} = \lt( \lt( \pi^1:\Xi^1\rightarrow B, \sigma^1:
B \rightarrow \Xi^1 \rt), g^1: \Xi^1 \rightarrow X \rt).
\end{equation}
The claim is that $\overline{\xi}^1_{e+1}$ is inductable.  Denote $h^1
= g^1\circ \sigma^1$.  Since $\widetilde{h}$ is very twistable and
since $h^1$ is a small deformation of $\widetilde{h}$, it follows by
proposition~\ref{prop-twistopen} that $h^1$ is very twistable.  Now
the image of $\widetilde{\overline{\xi}}^0_{e+1}$ is not contained in
the locus of very stable maps, precisely because of the point $r\in L$
where $L\cap M = \{r\}$: the stable map
$\widetilde{\overline{\xi}}^0_{e+1}(r)$ is not very stable.  But from
the description of the smoothing of
$\widetilde{\overline{\xi}}^0_{e+1}$ given in remark~\ref{rmk-modif},
we can find small deformations whose image is contained in the very
stable locus.  So we may assume the image of $\overline{\xi}^1_{e+1}$
is contained in the locus of very stable maps.  It remains only to
show that $\overline{\xi}^1_{e+1}$ is positive.

\

That $h^1:B \rightarrow T$ is a stable map follows from the fact that
$B$ is smooth and $h^1$ is non-constant: in fact the degree of
$h^1(B)$ equals the degree of $\widetilde{h}(\widetilde{B})$, which is
just $\text{deg}(h(B)) +1$.  This proves ($1$) of
definition~\ref{defn-pos}.  To show that $h^1:B\rightarrow X$ is
unobstructed, it suffices to prove that
$\widetilde{h}:\widetilde{B}\rightarrow X$ is unobstructed.  By
assumption, $h^0:B \rightarrow X$ is unobstructed (being a small
deformation of $h:B \rightarrow X$).  And the restriction to $L$ of
$T_X$ is generated by global sections (since $L$ is general).  Thus
$\widetilde{h}:\widetilde{B} \rightarrow X$ is unobstructed.  This
proves ($2$) of definition~\ref{defn-pos}. The map
\begin{equation}
\text{pr}\circ \widetilde{\overline{\xi}}^0_{e+1}: \widetilde{B}
\rightarrow \Kgnb{0,0}{X,e}
\end{equation}
is constant on $L\subset \widetilde{B}$, and on $B\subset
\widetilde{B}$, it is just $\zeta'$.  Thus the image is contained in
the unobstructed locus.  So we have ($3$) of definition
~\ref{defn-pos}.  The pullback of $T_{\Kgnb{0,0}{X,e}}$ to
$\widetilde{B}$ is trivial on $L$ and equals
$(\zeta')^*T_{\Kgnb{0,0}{X,e}}$ on $B$.  Since
$(\zeta')^*T_{\Kgnb{0,0}{X,e}}$ is ample, we conclude from
lemma~\ref{lem-DAcrit} that the pullback to $\widetilde{B}$ of
$T_{\Kgnb{0,0}{X,e}}$ is deformation ample.  Since
$\overline{\xi}^1_{e+1}$ is a small deformation of
$\widetilde{\overline{\xi}}^0_{e+1}$, the pullback of
$T_{\Kgnb{0,0}{X,e}}$ via $\overline{\xi}^1_{e+1}$ is deformation
ample by lemma~\ref{lem-opDA}.  This proves ($4$) of definiton
~\ref{defn-pos}.  Now for $\widetilde{\overline{\zeta}}^0_{e+1}$, the
pullback of $\widetilde{\sigma}^*\OO_{\widetilde{\Sigma}}(\sigma)$ is
trivial when restricted to $B\subset \widetilde{B}$, but on $L\subset
\widetilde{B}$ it is $\OO_L(1)$.  So by lemma~\ref{lem-DAcrit}, it is
deformation ample.  Since $\overline{\zeta}^1_{e+1}$ is a small
deformation of $\widetilde{\overline{\xi}}^0_{e+1}$, by
lemma~\ref{lem-opDA} we have that $(\sigma^1)^*\OO_{\Xi^1}(\sigma^1)$
is deformation ample.  So $\overline{\xi}^1_{e+1}$ is positive.
Therefore it is inductable.

\

To define the $1$-morphisms $\overline{\xi}^k_{e+1}:B \rightarrow
\Kgnb{0,1}{X,e+1}$ we repeat the procedure above.  Given
$\overline{\xi}^k_{e+1}$ and a point $b\in B$ whose image is a general
point of $\Delta_{(1,e)}$, say $(\Xi^k_b,\lambda^k_b,g^k_b:\Xi^k_b
\rightarrow X)$ with $\Xi^k_b = L \cup M$, $\lambda^k_b \in L$ and
$g^k_b:L \rightarrow X$ an embedding of a twistable line, we define
$\widetilde{B}^k$ to be the connected sum of $B$ and $L$ with $b\in B$
identified with $\lambda^k_b \in L$.  Then we define
\begin{equation}
\widetilde{\overline{\xi}}^k_{e+1}:\widetilde{B}^k \rightarrow
\Kgnb{0,1}{X,e+1}
\end{equation}
to be the modification of $\overline{\xi}^k_{e+1}$ associated to $b\in
B$, $L\subset \Xi^k_b$ and $\lambda^k_b \in L$.  By the same argument
as above, deformations of $\widetilde{\overline{\xi}}^k_{e+1}$ smooth
the node of $\widetilde{B}^k$, and a small deformation
$\overline{\xi}^{k+1}_{e+1}$ of $\widetilde{\overline{\xi}}^k_{e+1}$
is inductable.

\

It is quite simple to work out the degrees of the pullbacks by
$\widetilde{\overline{\xi}}^0_{e+1}$ of the tautological divisors of
$\Kgnb{0,1}{\PP^N,e+1}$ in terms of the degrees of the pullbacks by
$\zeta_1$ and $\overline{\zeta}_e$ of the tautological divisors.  This
is left to the interested reader.
\end{proof}

\section{Twistable lines} \label{sec-twist}

In this section, we will prove that if $n+1\geq d^2$, then for a
general hypersurface $X_d \subset \PP^n$ of degree $d$ and a general
line $L\subset X$, we have that $L$ is a twistable line on $X$.  To
prove this we introduce some incidence correspondences.  Let $N_d =
\binom{n+d}{n}-1$ and let $\PP^{N_d}$ denote the projective space
parametrizing hypersurfaces $X_d\subset \PP^n$ of degree $d$.  Let
${\mc X} \subset \PP^{N_d} \times \PP^n$ denote the universal family
of degree $d$ hypersurfaces in $\PP^n$.  Let $\GG(1,n)$ denote the
Grassmannian variety of lines in $\PP^n$.  Let $F({\mc X}) \subset
\PP^{N_d} \times \GG(1,n)$ denote the parameter space of pairs
$([X],[L])$ where $X\subset \PP^n$ is a hypersurface of degree $d$,
$L\subset \PP^n$ is a line and $L\subset X$.  Observe that the
projection $F({\mc X}) \rightarrow \GG(1,d)$ is a projective bundle of
relative dimension $N_d - (d+1)$.

\

Let $P(t) = (t+1)^2$ denote the Hilbert polynomial of a quadric
surface in $\PP^3$, and let $U \subset \textit{Hilb}^{P(t)}({\PP^n})$
denote the open subscheme parametrizing subschemes $\Sigma \subset
\PP^n$ which are projectively equivalent to a smooth, quadric surface
in $\PP^3 \subset \PP^n$.  Let $V \subset U \times \GG(1,n)$ denote
the parameter space of pairs $([\Sigma],[L])$ where $\Sigma \subset
\PP^n$ is a smooth quadric surface, where $L\subset \PP^n$ is a line,
and where $L\subset \Sigma$.  The projection map $V\rightarrow U$ has
a Stein factorization $V \rightarrow \tilde{U} \rightarrow U$ where
$\tilde{U} \rightarrow U$ is a finite, \'etale double cover, and where
$V\rightarrow \tilde{U}$ is a $\PP^1$-bundle.  Let $W\subset \PP^{N_d}
\times U \times \GG(1,n)$ denote the parameter space of triples
$([X],[\Sigma],[L])$ where $X\subset \PP^n$ is a hypersurface of
degree $d$, $([\Sigma],[L])$ is a point of $V$ and where $\Sigma
\subset X$.  The projection map $W \rightarrow V$ is a projective
bundle of relative dimension $N_d - (d+1)^2$.

\

Now for a triple $([X],[\Sigma],[L])\in W$, there is a map
(well-defined up to nonzero scalar) $d_X:\CC^{n+1} \rightarrow
H^0(\PP^n,\OO_{\PP^n}(d-1))$ which evaluates the partial derivatives
of a defining equation of $X$.  We may compose this map with the
restriction map $H^0(\PP^n,\OO_{\PP^n}(d)) \rightarrow
H^0(\Sigma,\OO_{\Sigma}(d-1))$. Let this map be denoted by
$d_{X,\Sigma}:\CC^{n+1} \rightarrow H^0(\Sigma,\OO_{\Sigma}(d-1))$.
More precisely, let $E$ be the trivial vector bundle on $W$ of rank
$n+1$, let $G$ be the vector bundle on $U$ whose fiber at a point
$\Sigma$ is just $H^0(\Sigma,\OO_{\Sigma}(d-1))$, and let $F$ be the
vector bundle on $W$ which is
$\text{pr}_1^*(\OO_{\PP^{N_d}}(1))\otimes \text{pr}_2^*G$.  Then there
is a map of vector bundles $d:E \rightarrow F$ whose fiber over
$([X],[\Sigma],[L])$ is the map $d_{X,\Sigma}$ constructed above.  Let
$W^o\subset W$ denote the open subscheme (possibly empty) over which
$d$ is surjective (i.e. the complement of the support of the cokernel
of $d$).

\begin{lem} \label{lem-quadric1}
  For any point $([X],[\Sigma],[L])\in W^o$, we have
\begin{enumerate}

\item $X$ is smooth along $\Sigma$,

\item $H^i\lt( \Sigma, N_{\Sigma/X} \rt)$ is zero for $i>0$,

\item $H^i\lt( \Sigma,N_{\Sigma/X}(-L) \rt)$ is zero for $i>0$,
  
\item $H^i\lt( \Sigma,N_{\Sigma/X}\otimes \OO_{\Sigma}(-1) \rt)$ is
  zero for $i>0$,

\item $H^1\lt( L, N_{L/X} \rt)$ is zero,

\item $H^1\lt( L, N_{L/X}(-1) \rt)$ is zero,
  
\item the projection morphism $W \rightarrow \PP^{N_d}$ given by
  $([X],[\Sigma],[L]) \mapsto [X]$ is smooth at $([X],[\Sigma],[L])$,
  
\item the projection morphism $F({\mc X}) \rightarrow \PP^{N_d}$ given
  by $([X],[L]) \mapsto [X]$ is smooth at $([X],[L])$, and

\item the projection morphism $\pi:W\rightarrow F({\mc X})$ given by
  $([X],[\Sigma],[L]) \mapsto ([X],[L])$ is smooth at
  $([X],[\Sigma],[L])$.  
\end{enumerate}
\end{lem}

\begin{proof}
  Since the partial derivatives of a defining equation of $X$ generate
  $H^0(\Sigma,\OO_{\Sigma}(d-1))$, in particular the locus where they
  all vanish is disjoint from $\Sigma$.  By the Jacobian criterion, we
  conclude that $X$ is smooth along $\Sigma$.

\

For a smooth quadric surface $\Sigma \subset \PP^3 \subset \PP^n$, we
have a short exact sequence:
\begin{equation}
\begin{CD}
0 @>>> N_{\Sigma/\PP^3} @>>> N_{\Sigma/\PP^n} @>>>
N_{\PP^3/\PP^n}|_{\Sigma} @>>> 0
\end{CD}
\end{equation}
Since $N_{\Sigma/\PP^3} \cong \OO_{\Sigma}(2)$ and since
$N_{\PP^3/\PP^n}|_{\Sigma} \cong \OO_{\Sigma}(1)^{\oplus (n-3)}$, we
have a short exact sequence:
\begin{equation}
\begin{CD}
0 @>>> \OO_{\Sigma}(2) @>>> N_{\Sigma/\PP^n} @>>>
\OO_{\Sigma}(1)^{\oplus (n-3)} @>>> 0
\end{CD}
\end{equation}
From this it is easy to compute that $H^i\lt( \Sigma,N_{\Sigma/\PP^n}
\rt)$, $H^i\lt( \Sigma,N_{\Sigma/\PP^n}\otimes \OO_{\Sigma}(-1) \rt)$
and $H^i\lt( \Sigma,N_{\Sigma/\PP^n}\otimes \OO_{\Sigma}(-L) \rt)$ are
all zero for $i>0$.

\

There is a short exact sequence:
\begin{equation}
\begin{CD}
0 @>>> N_{\Sigma/X} @>>> N_{\Sigma/\PP^n} @>>> N_{X/\PP^n}|_\Sigma
@>>> 0
\end{CD}
\end{equation}
Of course $N_{X/\PP^n}|_{\Sigma} \cong \OO_{\Sigma}(d)$, and for ${\mc
  L} = \OO_{\Sigma}$, for ${\mc L} = \OO_{\Sigma}(-L)$, and for ${\mc
  L} = \OO_{\Sigma}(-1)$, we compute that $H^i\lt( \Sigma,
\OO_{\Sigma}(d) \otimes {\mc L} \rt)$ is zero for $i>0$.  The bundle
$N_{\Sigma/\PP^n}$ was computed in the last paragraph.  For the line
bundles ${\mc L} = \OO_\Sigma$, ${\mc L} = \OO_{\Sigma}(-L)$ and ${\mc
  L} = \OO_{\Sigma}(-1)$, we computed that $H^i(\Sigma,
N_{\Sigma/\PP^n}\otimes {\mc L})$ is zero for $i>0$.  It immediately
follows from the long exact sequence in cohomology that $H^2\lt(
\Sigma, N_{\Sigma/X} \otimes {\mc L} \rt)$ is zero.  We also conclude
that $H^1(\Sigma, N_{\Sigma/X}\otimes {\mc L})$ is zero iff the
corresponding map
\begin{equation}
H^0(\Sigma, N_{\Sigma/\PP^n}\otimes {\mc L}) \rightarrow H^0(\Sigma,
N_{X/\PP^n}|_\Sigma) 
\end{equation}
is surjective.  

\

In the case that ${\mc L} = \OO_\Sigma(-1)$, the map from the last
paragraph factors the map
\begin{equation}
H^0(\PP^n,T_{\PP^n}(-1)) \rightarrow H^0(\Sigma,\OO_{\Sigma}(d-1)).
\end{equation}
But this map is precisely the map $d_{X,\Sigma}:\CC^{n+1} \rightarrow
H^0(\Sigma,\OO_{\Sigma}(d-1))$.  Since $d_{X,\Sigma}$ is surjective,
we conclude that $H^1(\Sigma, N_{\Sigma/X}\otimes \OO_{\Sigma}(-1))$
is zero.

\

To see that $H^1(\Sigma,N_{\Sigma/X}\otimes \OO_{\Sigma}(-L))$ is
zero, observe that we have $\OO_{\Sigma}(1) \cong \OO_{\Sigma}(L+L')$
where $L'\subset \Sigma$ is any line in the ruling opposite to the
ruling of $L$.  Thus we have a commutative diagram:
\begin{equation}
\begin{CD}
H^0( \Sigma , N_{\Sigma/\PP^n}(-1) ) \otimes_{\CC}
H^0( \Sigma, \OO_{\Sigma}(L') ) @>>> H^0( \Sigma, N_{\Sigma/\PP^n}(-L) \\
@VVV @VVV \\
H^0( \Sigma , N_{X/\PP^n}|_{\Sigma}(-1) )
\otimes_{\CC} H^0( \Sigma , \OO_{\Sigma}(L') ) @>>>
H^0( \Sigma, N_{X/\PP^n}|_{\Sigma}(-L) ) \\
\end{CD}
\end{equation}
The top vertical horizontal arrow is surjective by the last paragraph.
Moreover, the right vertical arrow is
\begin{equation}
\begin{array}{cc}
H^0\lt( \Sigma,
\OO_{\Sigma} \lt( (d-1)L + (d-1)L' \rt) \rt) \otimes H^0\lt( \Sigma,
\OO_{\Sigma}(L') \rt) \ \ \ \ \ \ \ \ \ \ \ \ \ \\
\ \ \ \ \ \ \ \ \ \ \ \ \rightarrow 
H^0\lt( \Sigma, \OO_{\Sigma}\lt( (d-1)L + dL' \rt) \rt),
\end{array}
\end{equation}
which is clearly surjective.  Therefore we conclude that the bottom
horizontal arrow is also surjective, i.e. $H^1(\Sigma,
N_{\Sigma/X}(-L) )$ is zero.  The proof that $H^1(\Sigma,
N_{\Sigma/X})$ is zero is almost identical to the proof that
$H^1(\Sigma, N_{\Sigma/X}(-L)$ is zero and is left to the reader.

\

To see that $H^1\lt( L, N_{L/X}(-1) \rt)$ is zero, first observe we
have a short exact sequence:
\begin{equation}
\begin{CD}
0 @>>> N_{\Sigma/X}(-1) @>>> N_{\Sigma/X}(-L') @>>> N_{\Sigma/X}|_L(-1) @>>> 0
\end{CD}
\end{equation}
By our computations and the long exact sequence in cohomology, we
conclude that $H^1\lt( L, N_{\Sigma/X}|_L(-1) \rt)$ is zero.  Next
observe that we have a short exact sequence:
\begin{equation}
\begin{CD}
0 @>>> N_{L/\Sigma}(-1) @>>> N_{L/X}(-1) @>>> N_{\Sigma/X}|_L(-1) @>>> 0
\end{CD}
\end{equation}
Of course $N_{L/\Sigma} \cong \OO_L(1)$, so $H^1\lt(
L,N_{L/\Sigma}(-1) \rt)$ is zero.  And we have seen that $H^1\lt(
L,N_{\Sigma/X}|_L(-1) \rt)$ is zero.  Therefore by the long exact
sequence in cohomology, we conclude that $H^1\lt(L, N_{L/X}(-1) \rt)$
is zero. By an almost identical argument, we also conclude that
$H^1\lt( L, N_{L/X} \rt)$ is zero.

\

Now by \cite[proposition I.2.14.2]{K}, the obstruction space for the
relative Hilbert scheme $\textit{Hilb}^{P(t)}({\mc X}/\PP^{N_d})$ at a
point $([X],[\Sigma])$ is contained in $H^1(\Sigma,N_{\Sigma/X})$.
Since the obstruction space vanishes, it follows by \cite[theorem
2.10]{K} that $\textit{Hilb}^{P(t)}({\mc X}/\PP^{N_d}) \rightarrow
\PP^{N_d}$ is smooth at $([X],[\Sigma])$.  As we have seen $W
\rightarrow \textit{Hilb}^{P(t)}({\mc X}/\PP^{N_d})$ is smooth.
Therefore the composition $W \rightarrow \PP^{N_d}$ is smooth along
$W^o$.

\

For basically the same reason as above, we conclude that the
projection map $F({\mc X}) \rightarrow \PP^{N_d}$ is smooth along the
image of $\pi:W^o \rightarrow F({\mc X})$.  Since $W^o \rightarrow
\PP^{N_d}$ is smooth, and since $F({\mc X}) \rightarrow \PP^{N_d}$ is
smooth along the image of $W^o$, to prove that $\pi$ is smooth along
$W^o$, it suffices to check that the derivative map $d\pi:
T_{W^o/\PP^{N_d}} \rightarrow \pi^* T_{F({\mc X})/\PP^{N_d}}$ is
surjective at every point.  This exactly reduces to the statement that
$H^0 \lt( \Sigma, N_{\Sigma/X} \rt) \rightarrow H^0\lt( L,
N_{\Sigma/X}|_L \rt)$ is surjective.  Since the cokernel is contained
in $H^1\lt( \Sigma, N_{\Sigma/X}(-L) \rt)$, which is zero, we conclude
the map is surjective.  Therefore $\pi:W^o \rightarrow F({\mc X})$ is
smooth.
\end{proof}

Now suppose that $([X],[\Sigma],[L])$ is a point in $W^o$.  We
associate to this triple a morphism $\zeta: L \rightarrow
\Kgnb{0,1}{X,1}$ as follows.  Let $\sigma:L \rightarrow \Sigma$ be the
inclusion and let $\text{pr}_L: \Sigma \rightarrow L$ be the unique
projection such that $\sigma$ is a section of $\text{pr}_L$.  Let
$g:\Sigma \rightarrow X$ be the inclusion.  Then we have a morphism
$\zeta:L \rightarrow \Kgnb{0,1}{X,1}$ given by the data:
\begin{equation}
\zeta = \lt(\text{pr}_L:\Sigma \rightarrow L, \sigma:L \rightarrow
\Sigma, g:\Sigma 
\rightarrow X \rt).
\end{equation}

\begin{lem} \label{lem-quadric2}
  If $([X],[\Sigma],[L])$ is a triple in $W^o$ and if $X$ is smooth,
  then the corresponding morphism $\zeta:L \rightarrow
  \Kgnb{0,1}{X,1}$ is flexible.
\end{lem}

\begin{proof}
  We need to check the axioms of definition~\ref{defn-twisting}.
  Since $g\circ \sigma:L \rightarrow X$ is an embedding, in particular
  this map is stable, i.e. axiom $(1)$ is satisfied.  By part $(5)$ of
  lemma~\ref{lem-quadric1}, we conclude that $\Kgnb{0,0}{X,1}$ is
  unobstructed at $[g\circ \sigma: L \rightarrow X]$, i.e. axiom $(2)$
  is satisfied.

\

To check axiom $(3)$, consider the normal bundle ${\mc N}$ of the
regular embedding $(\text{pr}_L, g): \Sigma \rightarrow L \times X$.
This fits into a short exact sequence:
\begin{equation}
\begin{CD}
0 @>>> \text{pr}_L^* T_L @>>> {\mc N} @>>> N_{\Sigma/X} @>>> 0
\end{CD}
\end{equation}
By part $(3)$ of remark~\ref{rmk-twisting}, we need to check that $R^1
\lt(\text{pr}_L\rt)_* {\mc N}(-\sigma)$ is zero.  
It is clear that for each fiber $L'$ of $\text{pr}_L:\Sigma
\rightarrow L$, we have that ${\mc N}(-\sigma)|_{L'}$ is just
$N_{L'/X}(-1)$.  By part $(6)$ of lemma~\ref{lem-quadric1}, we
conclude that $H^1\lt(L',N_{L'/X}(-1)\rt)$ is zero.  Therefore we
conclude that $R^1 \lt( \text{pr}_L \rt)_* {\mc N}(-\sigma)$ is zero.

\

By part $(3)$ of remark~\ref{rmk-twisting}, we have that $\zeta^*
T_{\text{ev}}$ is equivalent to $\lt( \text{pr}_L \rt)_* {\mc
  N}(-\sigma)$.  Twisting the short exact sequence above by
$\OO_{\Sigma}(-L)$ and applying the long exact sequence of higher
direct images, we see that $\lt( \text{pr}_L \rt)_* {\mc N}(-\sigma)$
fits between $\lt( \text{pr}_L \rt)_* \text{pr}_L^* T_L(-L)$ and $\lt(
\text{pr}_L \rt)_* N_{\Sigma/X}(-L)$ with cokernel $R^1 \lt(
\text{pr}_L \rt)_* \text{pr}_L^* T_L(-L)$.  For any fiber $L'$ of
$\text{pr}_L$, we have $T_L(-L)|_{L'}$ is isomorphic to
$\OO_{L'}(-1)$.  Therefore $\lt(\text{pr}_L\rt)_* \text{pr}_L^*
T_L(-1)$ and $R^1 \lt(\text{pr}_L \rt)_* \text{pr}_L^* T_L(-L)$ are
both zero, i.e. $\lt( \text{pr}_L \rt)_* {\mc N}(-L)$ is isomorphic to
$\lt( \text{pr}_L \rt)_* N_{\Sigma/X}(-L)$.

\

To show that axiom $(4)$ holds, we want to prove that for any point
$p'\in L$ with corresponding fiber $L' =
\text{pr}_L^{-1}\lt\{p'\rt\}$, we have that $\lt( \text{pr}_L \rt)_*
N_{\Sigma/X}(-L)\otimes \OO_L(-p')$ has no $H^1$.  Observe that since
$R^1\lt( \text{pr}_L \rt)_* \text{pr}_L^* T_L(-L-L')$ and $R^1 \lt(
\text{pr}_L \rt)_* {\mc N}(-L -L')$ are both zero, it follows from the
long exact sequence of higher direct images that also $R^1 \lt(
\text{pr}_L \rt)_* N_{\Sigma/X}(-L-L')$ is zero.  Therefore by the
Leray spectral sequence, we conclude that $H^1\lt( \Sigma,
N_{\Sigma/X}(-L-L') \rt)$ equals $H^1 \lt( L, \lt( \text{pr}_L
\rt)_*\lt( N_{\Sigma/X}(-L) \rt)(-p') \rt)$.  But by part $(4)$ of
lemma~\ref{lem-quadric1}, this is zero.  Thus axiom $(4)$ is
satisfied.

\

Finally, observe that $\sigma^* \OO_{\Sigma}(\sigma)$ is the trivial
line bundle, so axiom $(5)$ is satisfied.
\end{proof}

By lemma~\ref{lem-quadric2}, for any pair $([X],[L])\in F({\mc X})$,
if we can find a corresponding triple $([X],[\Sigma],[L])\in W^o$,
then it follows that $L$ is a twistable line on $X$.  Now we come to
the main result of this section.

\begin{prop} \label{prop-quadric}
  If $n+1 \geq d^2$ and $d\geq 2$, then $W^o\rightarrow F({\mc X})$ is
  dominant.  Thus for a general pair $([X],[L])\in F({\mc X})$, we
  have that $L$ is a twistable line on $X$.
\end{prop}

\begin{proof}
  By part $(9)$ of lemma~\ref{lem-quadric1}, it suffices to prove that
  $W^o$ is nonempty.  Let $I_d$ be the set of pairs of integers
  $I_d=\{(i,j) :0\leq i,j \leq d-1, i+j \geq 2\}$.  Choose coordinates
  on $\PP^n$ of the form $(Y_0,Y_1,Y_2,Y_3) \cup \lt( X_{i,j}
  \rt)_{(i,j)\in I_d} \cup (Z_m: m = 1,\dots, n+1-d^2)$.  Let $\Sigma
  \subset \PP^n$ be the smooth quadric surface with ideal
\begin{equation}
I_{\Sigma} = \langle Y_0 Y_3 - Y_1 Y_2 \rangle + \langle  X_{i,j} |
(i,j) \in I_d \rangle + \langle Z_m | m= 1, \dots, n+1-d^2 \rangle.
\end{equation}
This is the image of the embedding $f:\PP^1 \times \PP^1 \rightarrow
\PP^n$ given by sending a point $\lt( [U_0:U_1], [V_0:V_1] \rt)\in
\PP^1 \times \PP^1$ to
\begin{equation}
\lt( [U_0:U_1],[V_0:V_1] \rt) \mapsto [U_0V_0: U_0V_1: U_1V_0: U_1V_1:
0 : \dots :0 ].
\end{equation}

\

We make the following convention.  Given a pair $(i,j)$ of integers,
we set $k=\min(i,j)$, we set $i'=i-k$ and we set $j'=j-k$.  Consider
the hypersurface $X\subset \PP^n$ with defining equation
\begin{equation}
F = \lt( Y_0 Y_3 - Y_1 Y_2 \rt) Y_3^{d-2} + \sum_{(i,j) \in I_d} Y_0^k
Y_1^{i'} Y_2^{j'} Y_3 X_{i,j}.
\end{equation}
It is clear that $\Sigma \subset X$.  We claim that the derivative map
$dF:\CC^{n+1} \rightarrow H^0\lt(\Sigma,\OO_{\Sigma}(d-1) \rt)$ is
surjective.  Observe first that we have
\begin{eqnarray}
\frac{\partial F}{\partial Y_0} \mapsto U_1^{d-1}V_1^{d-1},
\frac{\partial F}{\partial Y_1} \mapsto U_1^{d-1}V_0V_1^{d-2}, \\
\frac{\partial F}{\partial Y_2} \mapsto U_0U_1^{d-2}V_1^{d-1},
\frac{\partial F}{\partial Y_3} \mapsto U_0U_1^{d-2}V_0V_1^{d-2}.
\end{eqnarray}
And observe that for $(i,j)\in I_d$, we have that
\begin{equation}
\frac{\partial F}{\partial X_{i,j}} \mapsto
U_0^iU_1^{d-1-i}V_0^jV_1^{d-1-j}.
\end{equation}
Since the partial derivatives of the form $\frac{\partial F}{\partial
  Y_i}$ give the terms $U_0^iU_1^{d-1-i}V_0^jV_1^{d-1-j}$ with $(i,j)
= (0,0), (0,1), (1,0),$ and $(1,1)$, and since these are precisely the
pairs $(i,j)$ not contained in $I_d$, we conclude that $dF$ is
surjective.  Thus, for any line $L\subset \Sigma$, we have that
$([X],[\Sigma],[L])$ is in $W^o$.
\end{proof}

\section{A very positive, very twisting family of lines}

In the last section, we proved that if $n+1 \geq d^2$, and $d\geq 2$,
then for a general hypersurface $X_d\subset \PP^n$ of degree $d$, a
general line $L\subset X$ is twistable, in other words
hypothesis~\ref{hyp-2} holds.  In this section, we will prove that if
$n\geq d^2 + d+2$, and $d\geq 3$ then there exists a morphism
$\overline{\zeta}_1:\PP^1 \rightarrow \Kgnb{0,1}{X,1}$ which is both
very twistable and very positive.  This provides the base case for the
induction argument of section~\ref{sec-induct}.

\

The arguments in this section are very similar to those of the last
section.  In that section, the key result was that for $([X],[L])$
general, there is a quadric surface $\Sigma$ with $L\subset \Sigma
\subset X$ such that $H^i\lt(\Sigma,N_{\Sigma/X}(-1)\rt)$ is zero for
$i>0$.  This result in turn reduced to finding a single degree $d$
polynomial $F$ on $\PP^n$, vanishing on some quadric surface $\Sigma$,
and such that
\begin{equation}
d_{F,\Sigma}:\CC^{n+1} \rightarrow H^0\lt(\Sigma,
\OO_{\Sigma}(d-1)\rt)
\end{equation}
is surjective.  

\

In this section, the role of $L\subset X$ will be replaced by a
rational normal curve $C_0\subset X$ of some degree $k\leq n$ (in the
end we will only need the case $k=2d-3$).  The role of the quadric
surface will be replace by a rational normal scroll $\Sigma$ of degree
$2k-1$ such that $C_0 \subset \Sigma \subset X$.  The cohomology
vanishing result of the last section will be replaced by the vanishing
of $H^i\lt(\Sigma, N_{\Sigma/X}(-C_0 - 2L) \rt)$ for $i>0$, where $L$
is any line of ruling of $\Sigma$.  The computation in this section
will be to find a single degree $d$ polynomial $F$ on $\PP^n$,
vanishing on $\Sigma$, and such that the image of the map,
\begin{equation}
d_{F,\Sigma}: \CC^{n+1} \rightarrow H^0\lt( \Sigma,\OO_{\Sigma}(d-1)
\rt),
\end{equation}
let's call it $W\subset H^0\lt( \Sigma, \OO_{\Sigma}(d-1) \rt)$, has
the property that the induced map
\begin{equation}
W \otimes H^0 \lt(\Sigma,\OO_{\Sigma}\lt( (k-3)L \rt) \rt) \rightarrow
H^0 \lt( \Sigma, \OO_{\Sigma}(d-1) \otimes \OO_{\Sigma}\lt( (k-3) L
\rt) \rt)
\end{equation}
is surjective.  A similar polynomial $F$ to that of the last section
satisfies this condition.

\subsection{Generating Linear Systems on ${\mathbb F}_1$}
In the last section, the relevant surface was the Hirzebruch surface
${\mathbb F}_0 = \PP^1 \times \PP^1$ embedded as a quadric surface.
In this section, the relevant surface is the Hirzebruch surface
${\mathbb F}_1$ embedded as a rational normal scroll of degree $2k-1$.
We will perform our computations using the projective model of
${\mathbb F}_1$:
\begin{equation}
{\mathbb F}_1 = \lt\{ ([T_0:T_1],[T_0U:T_1U:V]) \in \PP^1 \times \PP^2
| T_0 (T_1U) = T_1 (T_0U) \rt\}.
\end{equation}
In the equation above, $T_0U$ and $T_1U$ are simply variables on
$\PP^2$.  We denote the projection maps by $\text{pr}_1:{\mathbb F}_1
\rightarrow \PP^1$ and $\text{pr}_2: {\mathbb F}_1 \rightarrow \PP^2$.
We denote by $\OO_{{\mathbb F}_1}(F)$ the line bundle $\text{pr}_1^*
\OO_{\PP^1}$ and by $\OO_{{\mathbb F}_1}(E+F)$ the line bundle
$\text{pr}_2^* \OO_{\PP^2}$.  Here $\OO_{{\mathbb F}_1}(E)$ is the
divisor class of the directrix $E\subset {\mathbb F}_1$.  This
explains our terminology $T_0U$ and $T_1U$: $U$ is a nonzero element
of $H^0\lt({\mathbb F}_1,\OO_{{\mathbb F}_1}(E) \rt)$, and $T_0U$ and
$T_1U$ are the products of $U$ with the two global sections $T_0$ and
$T_1$ of $H^0\lt({\mathbb F}_1, \OO_{{\mathbb F}_1}(F) \rt)$.

\

We note that the line bundles $\OO_{{\mathbb F}_1}(E+F)$ and
$\OO_{{\mathbb F}_1}(F)$ generate the Picard group of ${\mathbb F}_1$.
Thus we adopt the terminology for line bundles on ${\mathbb F}_1$:
\begin{equation}
\OO(a,b) := \OO_{{\mathbb F}_1}\lt( a(E+F) + bF \rt).
\end{equation}
Note that $E+F$ and $F$ are each NEF, but not ample.  We conclude that
these two line bundles generate the NEF cone.  Thus every NEF line
bundle on ${\mathbb F}_1$ is of the form $\OO(a,b)$ for some
nonnegative integers $a,b$.

\

Now suppose that $\OO(a,b)$ is some NEF line bundle and $W\subset
H^0\lt( {\mathbb F}_1, \OO(a,b) \rt)$ is a linear series.

\begin{defn} \label{defn-cgen}
  For an integer $c\geq 0$, we say that $W$ is a $c$-\emph{generating
    linear system}, if the associated map
\begin{equation}
\mu_{W,c}: W \otimes H^0\lt({\mathbb F}_1, \OO(0,c)\rt)
\rightarrow H^0\lt( {\mathbb F}_1, \OO(a,b+c) \rt)
\end{equation}
is surjective.
\end{defn}

The question we want to answer is, when is $W$ a $c$-generating linear
system.  In particular, what is the minimal necessary dimension of a
$c$-generating linear system?  To simplify the answer, we write $b-1 =
\beta_d(c+1) + \beta_r$ where $\beta_d,\beta_r$ are integers with $0
\leq \beta_r < c+1$, and we write $a+b-1 = \alpha_d(c+1) + \alpha_r$
where $\alpha_d, \alpha_r$ are integers with $0\leq \alpha_r < c+1$.

\

\begin{lem} \label{lem-comput}
  Define the functions
\begin{eqnarray}
M(a,b,c) = a^2 + \lt(2b+3(c+1) \rt)a + 2b +4(c+1) + 2, \\
E(a,b,c) = \lt(\beta_r^2 - (c+1)\beta_r\rt) - 
\lt(\alpha_r^2 - (c-1)\alpha_r\rt) 
\end{eqnarray}
The minimal necessary dimension for a $c$-generating linear system
\begin{equation}
W\subset H^0\lt( {\mathbb F}_1, \OO(a,b)
\rt)
\end{equation}
is $\text{dim}(W) = \frac{1}{2(c+1)}\lt( M(a,b,c) + E(a,b,c) \rt)$.
\end{lem}

\begin{proof}
  This is just a computation.  For any nonnegative integers $a',b'$
  there is a decreasing filtration on $H^0\lt( {\mathbb F}_1,
  \OO(a',b') \rt)$ given by
\begin{equation}
F^iH^0\lt( {\mathbb F}_1, \OO(a',b')\rt) = H^0\lt( {\mathbb F}_1,
 \OO(a',b')(-iE)  \rt). 
\end{equation}
For any linear system $W\subset H^0\lt( {\mathbb F}_1, \OO_(a,b)
\rt)$, there is an induced filtration $F^iW = F^i\cap W$.  And the
multiplication map $\mu_{W,c}$ respects the filtrations on $W$ and on
$H^0\lt( {\mathbb F}_1, \OO(a,b+c) \rt)$.  If $\mu_{W,c}$ is
surjective, then the associated graded pieces
\begin{equation}
\text{gr}^i\mu_{W,c}: \text{gr}^iW \otimes H^0\lt( {\mathbb F}_1,
\OO(0,c) \rt) \rightarrow \text{gr}^i H^0\lt( {\mathbb
F}_1, \OO(a,b+c) \rt)
\end{equation}
are all surjective.  As the dimension of $W$ is the sum of the
dimensions of the pieces $\text{gr}^i W$, we should compute the
minimum possible dimension of a vector subspace $W^i \subset
\text{gr}^i H^0\lt( {\mathbb F}_1, \OO(a,b) \rt)$ such that
$\text{gr}^i\mu_{W^i,c}$ is surjective (where $\text{gr}^i\mu_{W^i,c}$
is the obvious map).

\

Of course the associated graded parts for $\OO(a',b')$ are just
\begin{equation}
\text{gr}^i H^0\lt( {\mathbb F}_1, \OO(a',b') \rt)
       = \lt\{ \begin{array}{ll}
                  H^0\lt( E, \OO_E(b'+i) \rt), & 0\leq i \leq a' \\
                  \{0\}, & i > a' 
                \end{array} \rt.
\end{equation}
So we are looking for a subset $W^i \subset H^0\lt( E, \OO_E(b+i)
\rt)$ such that the multiplication map
\begin{equation}
\text{gr}^i \mu_{W^i,c}: W^i \otimes H^0\lt( E, \OO_E(c) \rt)
\rightarrow H^0\lt( E, \OO_E(b+c+i) \rt)
\end{equation}
is surjective.  Counting the dimensions of the spaces on the left and
the right, we have
\begin{equation}
\text{dim}(W^i)\lt(c+1\rt) \geq \lt(b+c+i+1\rt),
\end{equation}
in other words,
\begin{equation}
\text{dim}(W^i) \geq \lt\lfloor\frac{b+i-1}{c+1} \rt\rfloor +2.
\end{equation}
On the other hand, we can acheive this bound: simply take $W^i$ to be
generated by the set of monomials:
\begin{equation}
\{ U^i V^{a-i} T_0^{(b+i)-j(c+1)}T_1^{j(c+1)}| j=1,\dots,r\}\cup  
\{ U^i V^{a-i} T_1^{b+i}\}
\end{equation}
where $r = \lt\lfloor\frac{b+i-1}{c+1} \rt\rfloor$.  Thus we have that
the minimum necessary dimension for a $c$-generating linear series is
\begin{equation}
\text{dim}(W) = 2a+2 + \sum_{i=0}^a \lt\lfloor \frac{b+i-1}{c+1}
\rt\rfloor.
\end{equation}
If we write $b-1 = \beta_d(c+1) + \beta_r$ with $0\leq \beta_r < c+1$
and if we write $a+b-1 = \alpha_d(c+1) + \alpha_r$ with $0 \leq
\alpha_r < c+1$, then the sum above is just
\begin{equation}
2a+2 + \sum_{i=0}^a \lt\lfloor \frac{b+i-1}{c+1} \rt\rfloor =
\frac{1}{2(c+1)}\lt( 
M(a,b,c) + E(a,b,c) \rt)
\end{equation}
where $M(a,b,c)$ and $E(a,b,c)$ are as above.
\end{proof}

For the next result, we introduce the Cox homogeneous coordinate ring:
\begin{equation}
S = S({\mathbb F}_1) := \CC[T_0,T_1,U,V] = \oplus_{(a,b) \in \ZZ^2}
H^0\lt( {\mathbb F}_1, \OO(a,b) \rt). 
\end{equation}
This is a $\ZZ^2$-graded ring, graded by $\text{deg}(T_0) =
\text{deg}(T_1) = (0,1)$, $\text{deg}(V) = (1,0)$ and $\text{deg}(U) =
(1,-1)$.  For any multi-degree $(a,b)\in \ZZ^2$, we have that
$S_{(a,b)} = H^0\lt( {\mathbb F}_1, \OO(a,b) \rt)$.  We put a graded
lexicographical monomial order on $S$ where the grading is by total
degree ($\text{deg(a,b)} = a+b$), and where $U>V>T_0>T_1$.

\

In the proof of the lemma above, a special role was played by the
linear system
\begin{eqnarray}
W_0(a,b,c) = \text{span}\lt\{ U^i V^{a-i}
T_0^{(b+i)-j(c+1)}T_1^{j(c+1)}|i=0,\dots, a,
j=1,\dots,r(i) \rt\}  \\
+ \text{span}\lt\{U^i V^{a-i}
T_1^{b+i}\rt|i=0,\dots,a\}.
\end{eqnarray} 
Here $r(i) = \lt\lfloor\frac{b+i-1}{c+1} \rt\rfloor$.

\begin{lem}\label{lem-comput2}
  Suppose $W\subset H^0\lt( {\mathbb F}_1, \OO(a,b) \rt)$ is a linear
  system.  If the vector space of initial terms $\text{IN}(W)$
  contains $W_0(a,b,c)$, then $W$ is a $c$-generating linear system.
\end{lem}

\begin{proof}
  Clearly we have that the vector space of initial terms of
  $\text{image}(\mu_{W,c})$ satisfies
\begin{equation}
\text{IN}(W)\cdot S_{(0,c)} \subset \text{IN}\lt(\text{image}(\mu_{W,c})
\rt).
\end{equation}
So if $\text{IN}(W)$ contains $W_0(a,b,c)$, then we have that
\begin{equation}
W_0(a,b,c)\cdot S_{(0,c)} \subset \text{IN}\lt( \text{image}(\mu_{W,c})
\rt).
\end{equation}
By construction, $W_0(a,b,c)\cdot S_{(0,c)} = S_{(a,b+c)}$.  So we
have $\text{IN}\lt( \text{image}(\mu_{W,c}) \rt) = S_{(a,b+c)}$, and
therefore $\text{image}(\mu_{W,c}) = S_{(a,b+c)}$.  Therefore $W$ is a
$c$-generating linear system.
\end{proof}

An important special case for us is when $a=d-1, b=(d-1)(k-1)$ and
$c=k-3$ for some positive integers $d$ and $k$ (here $d$ will be the
degree of the hypersurface $X\subset \PP^n$, and $k$ will be the
degree of the curve $C_0\subset X$).  In particular, if $k\geq 2d-3$,
then we have $b-1 = (d-1)(k-2) + d-2$ and $a+b-1 = (d-1)(k-2) + 2d-3$.
So the formulas above reduce to
\begin{equation}
\begin{array}{rcl}
M(d-1,(d-1)(k-1),k-3) & =  & 
2(k-2)(d-1)^2 + \\
5(k-2)(d-1) + 
4(k-2) & + &
\lt(3(d-1)^2 + 2(d-1) -2 \rt), \\
E(d-1,(d-1)(k-2),k-3) & = & 
(k-2)(d-1) - \\
\lt( 3(d-1)^2 + 2(d-1) - 2 \rt)
\end{array}
\end{equation}
So, if $k\geq 2d-3$, we have
\begin{equation}
\frac{1}{2(k-2)}(M+E) = (d-1)^2 + 3(d-1) + 2 = d^2 + d.
\end{equation}

\subsection{Cohomology Results}
We introduce some incidence correspondences, analogous to those
introduced in section~\ref{sec-twist}.  Let $N_d = \binom{n+d}{n} -1$
and let $\PP^{N_d}$ denote the projective space parametrizing
hypersurfaces $X_d \subset \PP^n$ of degree $d$.  Let ${\mc X} \subset
\PP^{N_d} \times \PP^n$ denote the universal family of degree $d$
hypersurfaces in $\PP^n$.  Let $k$ be any integer with $1\leq k \leq
\frac{n}{2}$ (later we will only need the case that $k=2d-3$).  Let
${\mc R}^k(\PP^n) \subset \textit{Hilb}^{kt+1}(\PP^n)$ denote the open
subscheme parametrizing curves $C_0\subset \PP^n$ which are
projectively equivalent to a degree $k$ rational normal curve $C_0
\subset \PP^k \subset \PP^n$.  Observe that ${\mc R}^k\lt(\PP^n\rt)$
is a homogeneous space of $\text{PGL}_{n+1}$, and therefore is smooth
and connected.  Let ${\mc R}^k({\mc X}) \subset \PP^{N_d} \times {\mc
  R}^k(\PP^n)$ denote the parameter space for pairs $([X],[C_0])$
where $C_0\subset X$.  Observe that the projection ${\mc R}^k({\mc X})
\rightarrow {\mc R}^k(\PP^n)$ is a projective bundle of relative
dimension $N_d - (kd+1)$.

\

Let $Q(t) = \frac{1}{2}(t+1)((2k-1)t+2)$ denote the Hilbert polynomial
of a rational normal scroll of degree $2k-1$ in $\PP^{2k}$.  Let ${\mc
  U}\subset \textit{Hilb}^{Q(t)}({\PP^n})$ denote the open subscheme
parametrizing subschemes $\Sigma \subset \PP^n$ which are projectively
equivalent to a rational normal scroll of degree $2k-1$ in
$\PP^{2k}\subset \PP^n$ which is abstractly isomorphic to ${\mathbb
  F}_1$.  Let ${\mc V}\subset {\mc U} \times {\mc R}^k(\PP^n)$ denote
the parameter space of pairs $([\Sigma],[C_0])$ where $C_0\subset
\Sigma$ and such that, using the isomorphism of $\Sigma$ and ${\mathbb
  F}_1$, the line bundle of $C_0$ is $\OO(1,0)$.  The projection map
${\mc V} \rightarrow {\mc U}$ factors as an open subset (with nonempty
fibers) of a projective bundle over ${\mc U}$ of relative dimension
$2$ (actually each fiber is isomorphic to the $\AAA^2$ of irreducible
curves in the linear system $|\OO(1,0)|$).

\

Let ${\mc W} \subset \PP^{N_d}\times {\mc U} \times {\mc R}^k(\PP^n)$
denote the parameter space for triples $([X],[\Sigma],[C_0])$ where
$([\Sigma],[C_0])$ is in ${\mc V}$ and where $\Sigma \subset X$.  The
projection map ${\mc W} \rightarrow {\mc V}$ is a projective bundle of
relative dimension $N_d - Q(d)$.

\

Now for a triple $([X],[\Sigma],[C_0]) \in {\mc W}$, we define
$d_{X,\Sigma}: \CC^{n+1} \rightarrow H^0\lt( \Sigma, \OO_{\Sigma}(d-1)
\rt)$ as in section~\ref{sec-twist}.  More precisely, let ${\mc E}$ be
the trivial vector bundle on ${\mc W}$ of rank $n+1$, let ${\mc G}$ be
the vector on ${\mc U}$ whose fiber at a point $\Sigma$ is just
$H^0\lt(\Sigma, \OO_{\Sigma}(d-1)\rt)$, and let ${\mc F}$ be the
vector bundle on ${\mc W}$ which is
$\text{pr}_1^*(\OO_{\PP^{N_d}}(1))\otimes \text{pr}_2^* {\mc G}$.
Then there is a map of vector bundles $d:{\mc E} \rightarrow {\mc F}$
whose fiber over $([X],[\Sigma],[C_0])$ is the map $d_{X,\Sigma}$
constructed above.  Let ${\mc W}^o\subset {\mc W}$ be the open
subscheme (possibly empty) parametrizing pairs $([X],[\Sigma],[C_0])$
such that the image of $d_{X,\Sigma}$ is a $(k-3)$-generating linear
series in $H^0\lt( \Sigma, \OO_{\PP^n}(d-3)|_\Sigma \rt)$.

\begin{lem} \label{lem-scroll1}
  Let $f:{\mathbb F}_1 \rightarrow \Sigma$ be an isomorphism to a
  degree $2k-1$ rational normal scroll $\Sigma \subset \PP^{2k}\subset
  \PP^n$.  For each pair of integers $a,b \geq 0$, consider the bundle
\begin{equation}
N(a,b) = f^*\lt( N_{\Sigma/\PP^n(-1}) \rt) \otimes 
  \OO(a,b) 
\end{equation}
and the subbundle
\begin{equation}
N'(a,b) = f^*\lt( N_{\Sigma/\PP^{2k} } \rt) \otimes
  \OO(a,b).
\end{equation}
Then we have the following:
\begin{enumerate}

\item $N'(0,0)$ is generated by global sections and satisfies $H^i\lt(
  {\mathbb F}_1, N'(0,0) \rt)$ is zero for $i>0$,

\item $N(0,0)$ is generated by global sections and satisfies $H^i\lt(
  {\mathbb F}_1, N(0,0) \rt)$ is zero for $i>0$,

\item if ${\mathcal F}$ is any coherent sheaf on ${\mathbb F}_1$ which
  is generated by global sections and satisfies $H^i\lt( {\mathbb
  F}_1, {\mathcal F} \rt)$ is zero for $i>0$, then for every pair of
  nonnegative integers $(a,b)$ we have that ${\mathcal F}(a,b) := {\mathcal
  F}\otimes \OO(a,b)$ is generated by global
  sections and satisfes $H^i\lt( {\mathbb F}_1, {\mathcal F}(a,b)
  \rt)$ is zero for $i>0$.

\end{enumerate}
In particular, we conclude that for any pair of nonnegative integers
$(a,b)$, we have that $N(a,b)$ (resp. $N'(a,b)$) is generated by
global sections and satisfies $H^i\lt({\mathbb F}_1, N(a,b)\rt)$ is
zero for $i>0$ (resp. $H^i\lt( {\mathbb F}_1, N'(a,b) \rt)$ is zero
for $i>0$).
\end{lem}

\begin{proof}
  Recall that $\text{pr}_1:{\mathbb F}_1 \rightarrow \PP^1$ is the
  projection morphism such that $\text{pr}_1^*\OO_{\PP^1}(1) =
  \OO_{{\mathbb F}_1}(f)$.  Via $\text{pr}_1$, ${\mathbb F}_1$ is
  isomorphic as a $\PP^1$-scheme to the total space of the projective
  bundle:
\begin{equation}
{\mathbb F}_1 \cong \PP\lt( \OO_{\PP^1}(-k) \oplus \OO_{\PP^1}(-(k-1))
\rt).
\end{equation}
Under this isomorphism $\OO(1,k-1)$ corresponds to the tautological
quotient bundle $\OO(1)$ on $\PP\lt( \OO_{\PP^1}(-k) \oplus
\OO_{\PP^1}(-(k-1)) \rt)$.  In other words, up to projective
equivalence, the map $f:{\mathbb F}_1 \rightarrow \PP^{2k}$
corresponds to the complete linear system of $\OO(1)$, i.e.
$f^*\OO_{\PP^{2k}}(1) \cong \OO(1)$.  Using this identification, it is
easy to see that we have a short exact sequence of vector bundles on
${\mathbb F}_1$:
\begin{equation}
0 \rightarrow \text{pr}_1^* T_{\PP^1} \rightarrow \text{pr}_1^*\lt(
\OO_{\PP^1}(1)^{\oplus (2k-1)} \rt)\otimes f^*\OO_{\PP^{2k}}(1) \rightarrow
f^* N_{\Sigma/\PP^{2k}} \rightarrow 0
\end{equation}
Twisting by $f^*\OO_{\PP^{2k}}(-1)$, we observe that $N'(0,0)$ is a
quotient of $\text{pr}_1^*\lt(\OO_{\PP^1}(1)^{\oplus (2k-1)} \rt)$,
and so is generated by global sections.  Observe that we have the
vanishing
\begin{equation}
\lt(\text{pr}_1\rt)_*\lt(\text{pr}_1^*T_{\PP^1}\otimes
f^*\OO_{\PP^{2k}}(-1) \rt) =
R^1\lt(\text{pr}_1\rt)_*\lt(\text{pr}_1^*T_{\PP^1}\otimes 
f^*\OO_{\PP^{2k}}(-1) \rt) = \{0\}.
\end{equation}
Applying the long exact sequence of higher direct images to our short
exact sequence (after twisting by $f^*\OO_{\PP^{2k}}(-1)$), we have
that $R^1\lt( \text{pr}_1 \rt)_* \lt(f^*N_{\Sigma/\PP^{2k}}(-1) \rt)$
is zero, and $\lt( \text{pr}_1 \rt)_* \lt(f^*N_{\Sigma/\PP^{2k}}(-1)
\rt)$ is $\OO_{\PP^1}(1)^{\oplus (2k-1)}$.  From this and the Leray
spectral sequence associated to $\text{pr}_1:{\mathbb F}_1 \rightarrow
\PP^1$, we conclude that $H^i\lt({\mathbb F}_1, N'(0,0) \rt)$ is zero
for $i>0$.  This proves item $(1)$ of the lemma.

\

Also observe that we have a short exact sequence:
\begin{equation}
\begin{CD}
0 @>>> N'(0,0) @>>> N(0,0) @>>> \OO(0,0)^{\oplus (n-2k)}
@>>> 0
\end{CD}
\end{equation}
Since $H^i\lt( {\mathbb F}_1, \OO_{{\mathbb F}_1} \rt)$ is zero for
$i>0$, we conclude that $N(0,0)$ is generated by global sections and
satisfies $H^i\lt({\mathbb F}_1, N(0,0) \rt)$ is zero for $i>0$.
This proves item $(2)$ of the lemma.

\

Now suppose that ${\mathcal F}$ is a coherent sheaf on ${\mathbb F}_1$
such that ${\mathcal F}$ is generated by global sections and such that
$H^i\lt( {\mathbb F}_1, {\mathcal F} \rt)$ is zero for $i>0$.  We will
prove by double induction on $(a,b)$ that the same is true for
${\mathcal F}(a,b) := {\mathcal F}\otimes \OO(a,b)$.

\

First we prove the result when $b=0$.  We proceed by induction on $a$.
If $a=0$, the result is tautological.  Thus suppose that $a>0$ and
suppose the result is proved for $a-1$.  Let $D\subset {\mathbb F}_1$
be a general member of the linear system $|\OO(1,0)|$.  Then $D$ is a
smooth curve isomorphic to $\PP^1$.  Since $D$ is general, we have a
short exact sequence:
\begin{equation}
\begin{CD}
0 @>>> {\mathcal F}(a-1,0) @>>> {\mathcal F}(a,0) @>>> {\mathcal F}(a,0)|_D
@>>> 0
\end{CD}
\end{equation}
Now $\mathcal{F}|_D$ is generated by global sections.  And
$\OO_{{\mathbb F}_1}(a(e+f))|_D$ is isomorphic to $\OO_{\PP^1}(a)$.
Thus we conclude that also ${\mathcal F}(a,0)|_D$ is generated by
global sections.  By the induction assumption, $H^1\lt(({\mathbb
  F}_1,{\mathcal F}(a-1,0) \rt)$ is zero, we conclude by the long
exact sequence of cohomology associated to the short exact sequence
above, that all the global sections of ${\mathcal F}(a,0)|_D$ lift to
global sections of ${\mathcal F}(a,0)$.  Therefore ${\mathcal F}(a,0)$
is generated by global sections.  Also, a coherent sheaf on $\PP^1$
which is generated by global sections has no higher cohomology.
Combining this with the induction assumption and using the long exact
sequence in cohomology associated to the short exact sequence above,
we conclude that $H^i\lt( {\mathbb F}_1, {\mathcal F}(a,0) \rt)$ is
zero for $i>0$.  Therefore we conclude by induction that for all
$a>0$, ${\mathcal F}(a,0)$ is generated by global sections and
$H^i\lt( {\mathbb F}_1, {\mathcal F}(a,0) \rt)$ is zero for $i>0$.

\

Now we prove the result with $b$ arbitrary.  We proceed by induction
on $b$.  If $b=0$, the result was proved in the last paragraph.  Thus
suppose that $b>0$ and suppose the result is proved for $b-1$.  Let $L
\subset {\mathbb F}_1$ be a general fiber of $\text{pr}_1$.  Then $L$
is smooth and isomorphic to $\PP^1$.  Since $L$ is general, we have a
short exact sequence:
\begin{equation}
\begin{CD}
0 @>>> {\mathcal F}(a,b-1) @>>> {\mathcal F}(a,b) @>>> {\mathcal
  F}(a,b)|_L @>>> 0
\end{CD}
\end{equation}
Now $\OO_{{\mathbb F}_1}(a(e+f)+bf)|_L$ is isomorphic to
$\OO_{\PP^1}(a)$.  By a similar analysis to that in the last
paragraph, we conclude that ${\mathcal F}(a,b)$ is generated by global
sections and that $H^i\lt( {\mathbb F}_1, {\mathcal F}(a,b) \rt)$ is
zero for $i>0$.  So item $(3)$ is proved by induction on $b$.
\end{proof}

\begin{lem} \label{lem-scroll2}
  Let $([X],[\Sigma],[C_0]) \in {\mc W}^o$ be any triple, and let
  $f:{\mathbb F}_1 \rightarrow \Sigma$ be some fixed isomorphism.  Let
  $N'(a,b)$ and $N(a,b)$ be as in lemma~\ref{lem-scroll1}.  Also let
  us denote
\begin{equation}
N_X(a,b) f^*\lt( N_{\Sigma/X}(-1) \rt)\otimes \OO(a,b)
\end{equation}  
Then we have the following:
\begin{enumerate}

\item $X$ is smooth along $\Sigma$
  
\item for each pair of nonnegative integers $(a,b)$, we have $H^i\lt(
  {\mathbb F}_1, N_X(a,b+k-3) \rt)$ is zero for $i>0$,
  
\item for any line of ruling $L\subset \Sigma$, $H^1\lt( L,
  N_{L/X}(a-1) \rt)$ is zero for $a\geq 0$ any integer,
  
\item $H^1\lt( C_0, N_{C_0/X}(a-2) \rt)$ is zero for $a\geq 0$ any
  integer,
  
\item the projection morphism ${\mc W} \rightarrow \PP^{N_d}$ given by
  $([X],[\Sigma],[C_0]) \mapsto [X]$ is smooth at
  $([X],[\Sigma],[C_0])$,
  
\item for any line $L\subset \Sigma$, the projection morphism $F({\mc
    X}) \rightarrow \PP^{N_d}$ given by $([X],[L]) \mapsto [X]$ is
  smooth at $([X],[L])$,
  
\item the projection morphism ${\mc R}^k({\mc X}) \rightarrow
  \PP^{N_d}$ given by $([X],[C_0]) \mapsto [X]$ is smooth at
  $([X],[C_0])$, and
  
\item the projection morphism $\pi: {\mc W} \rightarrow {\mc R}^k({\mc
    X})$ given by $([X],[\Sigma],[C_0]) \mapsto ([X],[C_0])$ is smooth
  at $([X],[\Sigma],[C_0])$.
\end{enumerate}

\end{lem} 

\begin{proof}
  Since the partial derivatives of a defining equation of $X$ give a
  $c$-generating linear series, in particular they generate the sheaf
  $\OO_{\Sigma}(d-1)$.  Thus, there is no point of $\Sigma$ at which
  all the partial derivatives vanish.  By the Jacobian criterion, we
  conclude that $X$ is smooth along $\Sigma$.  This proves item $(1)$.

\

For item $2$, we observe that we have a short exact sequence:
\begin{equation}
\begin{CD}
0 @>>> N_{\Sigma/X} @>>> N_{\Sigma/\PP^n} @>>> N_{X/\PP^n}|_\Sigma
@>>> 0
\end{CD}
\end{equation}
For ease of notation, define $\alpha = a + (d-1)$ and $\beta =
b+(d-1)(k-1)$.  We have a short exact sequence:
\begin{equation}
\begin{CD}
0 @>>> N_X(a,b) @>>> N(a,b) @>>> \OO(\alpha,\beta) @>>> 0 
\end{CD}
\end{equation}
When $a,b \geq 0$, it follows by lemma~\ref{lem-scroll1} that $H^i\lt(
{\mathbb F}_1, N(a,b) \rt)$ is zero for $i\geq 0$.  By a simple
calculation, we also see that
\begin{equation}
H^i\lt( {\mathbb F}_1,\OO(\alpha,\beta) \rt) = \{0\}, i\geq 0.
\end{equation}
So we conclude that $H^2\lt( {\mathbb F}_1,N_X(a,b) \rt)$ is zero, and
$H^1\lt( {\mathbb F}_1, N_X(a,b) \rt)$ is zero iff the following map
is surjective:
\begin{equation}
H^0\lt( {\mathbb F}_1, N(a,b) \rt) \rightarrow H^0\lt( {\mathbb F}_1,
\OO(\alpha,\beta) \rt).
\end{equation}
We have a commutative diagram:
\begin{equation}
\begin{CD}
H^0\lt( {\mathbb F}_1, N(a,b) \rt) \otimes H^0\lt( {\mathbb F}_1,
\OO(a',b') \rt) @>>> H^0\lt( {\mathbb F}_1, N(a+a',b+b') \rt) \\
@VVV @VVV \\
H^0\lt( {\mathbb F}_1, \OO\lt(\alpha,\beta
\rt) \rt) \otimes H^0\lt( {\mathbb F}_1, \OO(a',b') \rt) @>>>
H^0\lt( {\mathbb F}_1,
\OO\lt( \alpha+a', \beta+b' \rt) \rt)
\end{CD}
\end{equation}
We conclude that if $H^1\lt( {\mathbb F}_1, N_X(a,b) \rt)$ is zero and
if $a',b' \geq 0$, then we also have that $H^1\lt( {\mathbb F}_1,
N_X(a+a',b+b') \rt)$.  So to prove item $(2)$, we are reduced to the
case $a=0, b=k-3$.  But then the commutative diagram above factors as
\begin{equation}
\begin{CD}
H^0\lt( {\mathbb F}_1, T_{\PP^n} \rt) \otimes H^0\lt( {\mathbb F}_1,
\OO(0,k-3) \rt) @>>> H^0\lt( {\mathbb F}_1, N(0,k-3) \rt) \\
@VVV @VVV \\
H^0\lt( {\mathbb F}_1, \OO\lt(d-1,(d-1)(k-1)
\rt) \rt) \otimes H^0\lt( {\mathbb F}_1, \OO(0,k-3) \rt) @>>>  
H^0\lt( {\mathbb F}_1, \OO\lt(\alpha,\beta \rt) \rt)
\end{CD}
\end{equation}
By definition, the composition
\begin{equation}
H^0\lt( {\mathbb F}_1 T_{\PP^n} \rt) \otimes
H^0\lt( {\mathbb F}_1, \OO(0,k-3) \rt) \rightarrow H^0\lt( {\mathbb
  F}_1, \OO\lt(\alpha, \beta \rt) \rt)
\end{equation}
is surjective iff the triple $([X],[\Sigma],[C_0])$ is in ${\mathcal
  W}^o$.
So if $([X],[\Sigma],[C_0])$ is in ${\mathcal W}^o$, then $H^1\lt(
  {\mathbb F}_1, N_X(0,k-3) \rt)$ is zero.  Thus item $(2)$ holds.

\

For item $(3)$, observe that we have a short exact sequence:
\begin{equation}
\begin{CD}
0 @>>> N_{L/\Sigma}(a-1) @>>> N_{L/X}(a-1) @>>> N_{\Sigma/X}|_L(a-1) @>>> 0
\end{CD}
\end{equation}
Of course $N_{L/\Sigma} \cong \OO_L$, thus $H^1\lt( L,
N_{L/\Sigma}(a-1) \rt)$ is zero for $a\geq 0$.  So to prove item
$(3)$, it suffices to prove that $H^1\lt(L, N_{\Sigma/X}|_L(a-1) \rt)$
is zero.  Since $\OO(a-1,b)|_L \cong \OO_L(a-1)$, we have a short
exact sequence:
\begin{equation}
\begin{CD}
0 @>>> N_X(a,k-3) @>>> N_X(a,k-2) @>>> N_{\Sigma/X}|_L(a-1) @>>> 0
\end{CD}
\end{equation}
By item $(2)$, for $a\geq 0$ the higher cohomology of the first two
terms vanishes.  Thus by the long exact sequence in cohomology
associated to this short exact sequence, we have that $H^1\lt( L,
N_{\Sigma/X}|_L(a-1) \rt)$ is zero for $a\geq 0$.  This proves item
$(3)$.

\

Item $(4)$ is almost identical to item $(3)$ and is left as an
exercise to the reader.

\

By \cite[proposition 2.14.2]{K}, the obstruction space for the
relative Hilbert scheme $\textit{Hilb}^{Q(t)}({\mc X}/\PP^{N_d})$ 
at a point $([X],[\Sigma])$ is contained in $H^1\lt(\Sigma,
N_{\Sigma/X} \rt)$.  For a triple $([X],[\Sigma],[C_0])$ in ${\mc
W}^o$, it follows from item $(2)$ that 
\begin{equation}
H^1\lt( \Sigma, N_{\Sigma/X} \rt) = H^1\lt( {\mathbb F}_1, N_X(1,k-1)
\rt)
\end{equation} 
vanishes.  It follows by \cite[theorem 2.10]{K} that
$\textit{Hilb}^{Q(t)}({\mc X}/\PP^{N_d}) \rightarrow \PP^{N_d}$ is
smooth at $([X],[\Sigma])$.  Of course the projection ${\mc W}^o
\rightarrow \textit{Hilb}^{Q(t)}({\mc X}/\PP^{N_d})$ is an open subset
of a projective bundle, and so is smooth.  Thus we conclude that the
composite morphism ${\mc W}^o \rightarrow \PP^{N_d}$ is smooth.  This
proves item $(5)$.

\

Item $(6)$ is similar to item $(5)$ and follows from item $(3)$ which
shows that $H^1\lt( L, N_{L/X} \rt)$ is zero.

\

Item $(7)$ is similar to item $(5)$ and follows from item $(4)$ which
shows that $H^1\lt( C_0, N_{C_0/X} \rt)$ is zero.

\

Since ${\mc W}^o \rightarrow \PP^{N_d}$ is smooth at
$([X],[\Sigma],[C_0])$ and since ${\mc R}^k \rightarrow \PP^{N_d}$ is
smooth at $([X],[C_0])$, to prove that $\pi:{\mc W}^0 \rightarrow {\mc
  R}^k({\mc X})$ is smooth at $([X],[\Sigma],[C_0])$, it suffices to
check that the derivative map $d\pi: T_{{\mc W}^0/\PP^{N_d}}
\rightarrow \pi^* T_{{\mc R}^k({\mc X})/\PP^{N_d}}$ is surjective at
$([X],[\Sigma],[C_0])$.  This exactly reduces to the statement that
$H^0 \lt( \Sigma, N_{\Sigma/X} \rt) \rightarrow H^0 \lt( C_0,
N_{\Sigma/X}|_{C_0} \rt)$ is surjective.  Since the cokernel is
contained in $H^1\lt( {\mathbb F}_1, N_X(0,k-1) \rt)$, which is zero
by item $(3)$, we conclude the derivative $d\pi$ is surjective.
Therefore $\pi:{\mc W}^o \rightarrow {\mc R}^k({\mc X})$ is smooth.
This proves item $(8)$.
\end{proof}

Now suppose that $([X],[\Sigma],[C_0])$ is a point in ${\mc W}^o$.
Let $\sigma:C_0 \rightarrow \Sigma$ be the inclusion and let
$\text{pr}_{C_0}: \Sigma \rightarrow C_0$ be the unique projection
such that $\sigma$ is a section of $\text{pr}_L$ (in the model of
$\Sigma$ as ${\mathbb F}_1$, $\text{pr}_{C_0}$ is simply
$\text{pr}_1$).  Let $g:\Sigma \rightarrow X$ be the inclusion.  Then
we have an induced morphism $\zeta:C_0 \rightarrow \Kgnb{0,1}{X,1}$
given by the diagram:
\begin{equation}
\begin{CD}
\Sigma @> g >> X \\
@VV \text{pr}_{C_0} V \\
C_0
\end{CD}
\end{equation}

\begin{lem} \label{lem-scroll3}
  If $([X],[\Sigma],[C_0])$ is a triple in ${\mc W}^o$, then the
  corresponding morphism $\zeta:C_0 \rightarrow \Kgnb{0,1}{X,1}$ is
  very twisting and very positive.
\end{lem}

\begin{proof}
  First we check the axioms of definition~\ref{defn-verytwisting}.
  Since $g\circ \sigma:C_0 \rightarrow X$ is an embedding, in
  particular this map is stable, i.e. axiom $(1)$ is satisfied.  By
  item $(7)$ of lemma~\ref{lem-scroll2}, we conclude that
  $\Kgnb{0,0}{X,k}$ is unobstructed at $[g\circ \sigma: C_0
  \rightarrow X]$, i.e. axiom $(2)$ is satisfied.

\

The argument that axiom $(3)$ holds is exactly the same as the
argument that axiom $(3)$ holds in the proof of
lemma~\ref{lem-quadric2}, where item $(6)$ of lemma~\ref{lem-quadric1}
is replaced by item $(3)$ of lemma~\ref{lem-scroll2}.

\

As in the proof of lemma~\ref{lem-quadric2}, we have that
$\zeta^*T_{\text{ev}}$ is isomorphic to $\lt( \text{pr}_{C_0} \rt)_*
N_X(0,k-1)$.  Of course $\zeta^*T_{\text{ev}}$ is ample iff $H^1\lt(
C_0, \zeta^*T_{\text{ev}}(-2) \rt)$ is zero.  By the isomorphism above
and a Leray spectral sequence argument analogous to the one in the
proof of lemma~\ref{lem-quadric2}, this cohomology group equals
$H^1\lt( {\mathbb F}_1, N_X(0,k-3) \rt)$.  By item $(2)$ of
lemma~\ref{lem-scroll2}, this group is zero.  Therefore
$\zeta^*T_{\text{ev}}$ is an ample bundle.  So axiom $(4)$ is
satisfied.

\

Finally, observe that $\sigma^* \OO_{\Sigma}(\sigma)$ is
$\OO_{C_0}(1)$ and so is ample.  So axiom $(5)$ is satisfied.  Thus
$\zeta$ is a very twisting family.

\

Next we consider the axioms in definition~\ref{defn-pos}.  Axioms
$(1)$, $(2)$ and $(3)$ follow immediately from axioms $(1)$, $(2)$ and
$(3)$ of definition~\ref{defn-twisting} proved above.  To see that
axiom $(4)$ holds, observe that we have a short exact sequence:
\begin{equation}
\begin{CD}
0 @>>> \zeta^* T_{\text{ev}} @>>> \zeta^* T_{\Kgnb{0,1}{X,1}} @>>>
\lt( g\circ \sigma\rt)^* T_X @>>> 0
\end{CD}
\end{equation}
We have seen above that $\zeta^* T_{\text{ev}}$ is ample.  Moreover,
it follows by item $(4)$ of lemma~\ref{lem-scroll2} that $N_{C_0/X}$
is an ample vector bundle.  Since also $T_{C_0}$ is an ample line
bundle, we conclude that $T_X|_{C_0}$ is an ample vector bundle.
Since the first and third terms in the short exact sequence above are
ample, we conclude that $\zeta^* T_{\Kgnb{0,1}{X,1}}$ is an ample
vector bundle.  Since $\zeta^* \text{pr}^* T_{\Kgnb{0,0}{X,1}}$ is a
quotient bundle of $\zeta^* T_{\Kgnb{0,1}{X,1}}$, we also have that
$\zeta^* \text{pr}^* T_{\Kgnb{0,0}{X,1}}$ is an ample vector bundle.
So axiom $(4)$ holds.

\

Finally, we have seen above that $\sigma^*\OO_{\Sigma}(\sigma)$ equals
$\OO_{C_0}(1)$, which is ample.  Thus axiom $(5)$ holds.  So $\zeta$
is a very positive family.
\end{proof} 

\begin{prop} \label{prop-scroll}
If $n \geq d^2 + d + 2$ and if $d\geq 3$, for $k=2d-3$, we have that
${\mc W}^o \rightarrow {\mc R}^k({\mc X})$ is dominant, and ${\mc
R}^k({\mc X}) \rightarrow \PP^{N_d}$ is dominant.  So for $[X]\in
\PP^{N_d}$ general, there exists a very twisting, very positive family
$\zeta: C_0 \rightarrow \Kgnb{0,1}{X,1}$.
\end{prop}

\begin{proof}
  By item $(8)$ of lemma~\ref{lem-scroll2}, it suffices to prove that
  ${\mc W}^o$ is nonempty.  We have to find a pair $([X],[\Sigma])$
  such that for $a=d-1$, $b=(d-1)(k-1)$ and for $c=k-3$, we have that
  the image of the derivative map
\begin{equation}
d_{X,\Sigma}:H^0\lt( (\PP^n,T_{\PP^n}(-1) \rt) \rightarrow H^0\lt(
{\mathbb F}_1, \OO(a,b) \rt)
\end{equation}
is a $c$-generating linear system.

\

Recall that $S = \CC[T_0,T_1,U,V]$ is the $\ZZ^2$-graded Cox
homogeneous coordinate ring of ${\mathbb F}_1$.  Let $A_d$ denote the
set of monomials in the vector subspace $S_{(d-1,(d-1)(k-1))}$ which
occur in the linear system $W_0(a,b,c)$, i.e.
\begin{equation}
\begin{array}{c}
A_d = \\
\lt\{ U^i V^{d-1-i}
T_0^{((d-1)(k-1)+i)-j(k-2)}T_1^{j(k-2)}|i=0,\dots, d-1,
j=1,\dots,r(i) \rt\}  \\
\cup \lt\{U^i V^{d-1-i}
T_1^{(d-1)(k-1)+i}\rt|i=0,\dots,d-1\}
\end{array}
\end{equation}
where $r(i) = d-1+\lt\lfloor\frac{d-2+i}{k-2} \rt\rfloor$.
Let $B_d$ denote the set of
monomials of the form 
\begin{equation}
\begin{array}{c}
B_d = \\
\lt\{U^{d-1} T_0^{(d-1)k - i(k-2)} T_1^{i(k-2)} | i=0,\dots, d-2
\rt\} \cup \\
\lt\{ U^{d-2} V T_0^{(d-1)k-1-i(k-2)} T_1^{i(k-2)} | i=0,\dots, d-2 \rt\}
\cup \\
\lt\{ U^{d-3} V^2 T_0^{(d-1)k-2-(j+2)(k-2)} T_1^{(j+2)(k-2)} |
j=0,\dots, d-4 \rt\} \cup \\
\lt\{ U^{d-4} V^3 T_0^{(d-1)k-3-(j+2)(k-2)} T_1^{(j+2)(k-2)} | j=0,\dots,
d-4 \rt\}.
\end{array}
\end{equation}
Let $C_d$ denote the set of monomials $C_d = A_d - B_d$.  Now $A_d$
contains $d^2 + d$ monomials, and $B_d$ contains $4d-8$ monomials.
Choose homogeneous coordinates on $\PP^n$ of the form
\begin{equation}
\lt\{Y_0,\dots,Y_{2k} \rt\} \cup \lt\{ X_M | M \in C_d \rt\} \cup
\lt\{Z_l | l=1,\dots, n-(d^2+d+2) \rt\}.
\end{equation}

\

Let $f:{\mathbb F}_1 \rightarrow \PP^{2k} \subset \PP^n$ be the map
defined by sending a point $([T_0:T_1],[T_0U:T_1U:V])$ to the point in
$\PP^n$ with coordinates $X_m = 0, m\in C_d$, with $Z_l=0, l=1,\dots,
n-(d^2+d)$, and with
\begin{eqnarray}
Y_0 = T_0^kU,\dots, Y_i = T_0^{k-i}T_1^iU, \dots, Y_k = T_1^kU, \\
Y_{k+1} = T_0^{k-1}V,\dots, Y_{k+1+j} = T_0^{k-1-j}T_1^j V,\dots,
Y_{2k} = T_1^{k-1}V.
\end{eqnarray}
This is an embedding whose image $\Sigma = f({\mathbb F}_1)$ is a
rational normal scroll of degree $2k-1$.  

\

Now the pullback map $H^0\lt( (\PP^{2k},\OO_{\PP^{2k}}(1) \rt)
\rightarrow H^0\lt( {\mathbb F}_1, \OO(1,k-1) \rt)$ is surjective by
construction.  And the natural map
\begin{equation}
\text{Sym}^{d-1}H^0 \lt( {\mathbb F}_1, \OO(1,k-1) \rt) \rightarrow
H^0 \lt( {\mathbb F}_1, \OO\lt(d-1,(d-1)(k-1)\rt) \rt)
\end{equation} 
is surjective.  Therefore the pullback map 
\begin{equation}
H^0\lt( \PP^{2k},\OO_{\PP^{2k}}(d-1) \rt) \rightarrow H^0\lt( {\mathbb
F}_1, \OO\lt(d-1,(d-1)(k-1) \rt) \rt)
\end{equation}
is surjective.  For each monomial $M\in C_d$, choose a polynomial
$G_M(Y_0,\dots,Y_{2k})$ such that $f^*G_M = M$.  

\

Consider the hypersurface $X\subset \PP^n$ with defining equation
\begin{eqnarray}
F = \sum_{i=0}^{d-2} \lt(Y_iY_{k+2+i} - Y_{i+1}Y_{k+1+i}\rt)
Y_0^{d-2-i}Y_{k-3}^i + \\
\sum_{j=0}^{d-4}\lt(Y_{d-1+j}Y_{k+d+1+j} - Y_{d+j}Y_{k+d+j}
\rt)Y_0^{d-4-j}Y_{k-3}^jY_{2k-3}Y_{2k+1-d} + \\ 
\sum_{M\in C_d}G_M(Y_0,\dots,Y_{2k})X_M.
\end{eqnarray}
The corresponding derivative map $d_{X,\Sigma}$ acts on the partial
derivatives $\frac{\partial F}{\partial Y_i}$ by
\begin{equation}
\begin{array}{lcl}
\frac{\partial F}{\partial Y_0} & \mapsto & U^{d-2}VT_0^{(d-1)k-2}; \\
\frac{\partial F}{\partial Y_{i+1}} & \mapsto &
-U^{d-2} V T_0^{(d-1)k-1-i(k-2)} T_1^{i(k-2)} + \\
& & U^{d-2} V T_0^{(d-1)k-2-(i+1)(k-2)} T_1^{(i+1)(k-2)}, \\
& & i=0,\dots,d-3,
\\
\frac{\partial F}{\partial Y_{d-1}} & \mapsto &
-U^{d-2} V T_0^{(d-1)k-1-(d-2)(k-2)} T_1^{(d-2)(k-2)} + \\
& & U^{d-4} V^3 T_0^{(d-1)k-4-2(k-2)} T_1^{2(k-2)+1}, \\
\frac{\partial F}{\partial Y_{d+j}} & \mapsto &
-U^{d-4} V^3 T_0^{(d-1)k-3-(j+2)(k-2)} T_1^{(j+2)(k-2)} + \\
& & U^{d-4} V^3 T_0^{(d-1)k-4-(j+3)(k-2)} T_1^{(j+3)(k-2)},\\ 
& & j=0,\dots,
d-5, \\
\frac{\partial F}{\partial Y_{k-1}} & \mapsto & 
-U^{d-4} V^3
T_0^{(d-1)k-3-(d-2)(k-2)} T_1^{(d-2)(k-2)}, \\
\frac{\partial F}{\partial Y_k} & \mapsto & 0, \\
\frac{\partial F}{\partial Y_{k+1}} & \mapsto & -U^{d-1} T_0^{(d-1)k-1}T_1, \\
\frac{\partial F}{\partial Y_{k+2+i}} & \mapsto & U^{d-1} T_0^{(d-1)k -
  i(k-2)} T_1^{i(k-2)} - \\
& & U^{d-1} T_0^{(d-1)k-1-(i+1)(k-2)} T_1^{(i+1)(k-2)+1},\\
& & i=0,\dots, d-3,
\\
\frac{\partial F}{\partial Y_{k+d}} & \mapsto & U^{d-1}
T_0^{(d-1)k-(d-2)(k-2)} T_1^{(d-2)(k-2)} - \\
& & U^{d-3} V^2 T_0^{(d-1)k - 3 -2(k-2)} T_1^{2(k+2) + 1}, \\
\frac{\partial F}{\partial Y_{k+d+1+j}} & \mapsto & U^{d-3} V^2
T_0^{(d-1)k-2-(j+2)(k-2)} T_1^{(j+2)(k-2)} - \\
& & U^{d-3} V^2 T_0^{(d-1)k-3-(j+3)(k-2)} T_1^{(j+3)(k-2)+1},\\
& & j=0,\dots,
d-5, \\
\frac{\partial F}{\partial Y_{2k}} & \mapsto & U^{d-3} V^2
T_0^{(d-1)k-2-(d-2)(k-2)} T_1^{(d-2)(k-2)}.
\end{array}
\end{equation}

\

In the list above, every partial derivative is a sum of at most two
monomials, i.e. each partial derivative is a binomial, and the first
term listed is the initial term with respect to our monomial order.
From this list, it is clear that every monomial in $B_d$ occurs as the
initial term of some partial derivative.  Also, for each $M\in C_d$,
the partial derivative $\frac{\partial F}{\partial X_M}$ simply maps
to $M$.  Thus every monomial in $C_d$ occurs as the initial term of
some partial derivative.  Thus every monomial in $A_d$ occurs as the
initial term of some partial derivative.  So the vector space of
initial terms of $\text{image}(d_{X,\Sigma})$ contains $W_0(a,b,c)$.
Therefore, by lemma~\ref{lem-comput2}, we conclude that
$\text{image}(d_{X,\Sigma})$ is a $c$-generating linear system.  So,
for any irreducible curve $C_0$ in the linear system of $|\OO(1,0)|$,
we have that $([X],[\Sigma],[C_0])$ is in ${\mc W}^o$.
\end{proof}

\section{Proof of the main theorem}

In this section we prove theorem~\ref{thm-thm1}.  As explained at the
end of section~\ref{sec-results}, if $d< \frac{n+1}{2}$, then for a
general hypersurface $X_d\subset \PP^n$, hypothesis~\ref{hyp-1},
hypothesis~\ref{hyp-1.5}, and hypothesis~\ref{hyp-1.75} are satisfied.
By proposition~\ref{prop-quadric}, if $d\geq 2$ and $d^2 \leq n+1$,
then for $X_d\subset \PP^n$ general we have that
hypothesis~\ref{hyp-2} is satisfied.  Finally, if $d\geq 3$ and if
$d^2+d+2 \leq n$, then for $X_d\subset \PP^n$ general there exists a
very twisting, very positive family $\zeta:C_0 \rightarrow
\Kgnb{0,1}{X,1}$.  Thus $(\zeta,\zeta)$ is an inducting pair.  Now by
theorem~\ref{thm-induction}, we conclude that for every $e \geq 1$
there exists an inducting pair $(\zeta_1,\overline{\zeta}_e)$.  In
particular, there exists a very positive $1$-morphism
$\overline{\zeta}_e:C \rightarrow \Kgnb{0,1}{X,e}$.  As seen in the
proof of theorem~\ref{thm-induction}, we may assume that $C$ is
smooth, i.e. $C$ is equal to $\PP^1$, and that the image of
$\text{pr}\circ \overline{\zeta}_e:\PP^1 \rightarrow \Kgnb{0,0}{X,e}$
is contained in the smooth part of the fine moduli locus.  So, by item
$(2)$ of remark~\ref{rmk-pos}, we conclude that $\overline{\zeta}_e$
is a \emph{very free} morphism.  And by \cite[proposition 7.4]{HRS2},
we have that $\Kgnb{0,0}{X,e}$ is an irreducible variety.  Therefore
by \cite[theorem IV.3.7]{K}, we conclude that $\Kgnb{0,0}{X,e}$ is
rationally connected.

\bibliography{my}
\bibliographystyle{abbrv}

\end{document}